\documentclass{article}

\usepackage{amsfonts,pstricks}
\usepackage{amsmath,amsthm, amssymb,  mathrsfs, epsfig}
\usepackage{hyperref}
\usepackage{graphicx}
\usepackage{tikz}
\usepackage{caption}
\usepackage{authblk}
\usepackage{bm}
\usepackage{tensor}
\usepackage{tikz-cd}
\usetikzlibrary{decorations.markings}
\usepackage{multirow}
\usepackage{enumitem}
\usepackage{url}

\hypersetup{linktocpage = true}

\def\R{\mathbb{R}}
\def\Z{\mathbb{Z}}
\def\C{\mathbb{C}}
\def\Cat{\mathcal{C}}
\def\D{\mathcal{D}}
\def\T{\mathcal{T}}
\def\G{\mathcal{G}}
\def\Mod{\text{Mod}}

\def\cross#1#2{{#1}_{#2}^{\times}}
\def\L#1{L(#1)}
\def\Ltil#1{\tilde{L}(#1)}
\def\lsup#1#2{\tensor[^{#1}]{{#2}}{}}
\def\lsuprsub#1#2#3{\tensor[^{#1}]{{#2}}{_{#3}}}

\def\unit{\mathbf{1}}
\def\lrangle{\langle\;,\;\rangle}
\def\Sym{\text{S}}
\def\ColorS{\mathcal{S}}
\def\Vect{\mathcal{V}ect}
\def\CY{CY}
\def\Y{Y}
\def\DW{DW}
\def\Aut#1{\text{Aut}_{\otimes}(#1)}
\def\uG#1{\underline{#1}}

\def\GBSFC{G\text{-BSFC}}

\def\TQFT{\mathrm{TQFT}}
\def\hc{\circ_h}

\newcommand*\circled[1]{\tikz[baseline=(char.base)]{
            \node[shape=circle,draw,inner sep=1pt] (char) {#1};}}
\newcommand{\insertimage}[2]{\vcenter{\hbox{\includegraphics[scale=#1]{#2}}}}

\def\CatName#1{#1-crossed braided spherical fusion}
\def\CatNameTu#1{#1-ribbon}
\def\CatNameKi#1{#1-equivariant}

\def\CatNameShort#1{#1-BSFC}
\def\Sphere{\mathbb{S}}

\DeclareMathOperator{\Hom}{Hom}

\DeclareMathOperator{\Tr}{Tr}

\def \Yshape #1#2#3#4#5
{
\draw (#1,-1+ #2)--(#1,#2) -- (-1+ #1,1 + #2);
\draw (#1,#2) -- (1+ #1,1+ #2);
\draw (-1.5+#1,1+ #2) node{{\tiny ${#3}$}};
\draw (1.5+#1, 1+ #2) node{{\tiny ${#4}$}};
\draw (-0.5+#1,-1+ #2) node{{\tiny ${#5}$}};
}
\def \YshapeNormalFont #1#2#3#4#5
{
\draw (#1,-1+ #2)--(#1,#2) -- (-1+ #1,1 + #2);
\draw (#1,#2) -- (1+ #1,1+ #2);
\draw (-1.5+#1,1+ #2) node{ ${#3}$};
\draw (1.5+#1, 1+ #2) node{ ${#4}$};
\draw (-0.5+#1,-1+ #2) node{ ${#5}$};
}
\def \InverseYshape #1#2#3#4#5
{
\draw (#1,1+ #2)--(#1,#2) -- (-1+ #1,-1 + #2);
\draw (#1,#2) -- (1+ #1,-1+ #2);
\draw (-1.5+#1,-1+ #2) node{{\tiny ${#3}$}};
\draw (1.5+#1, -1+ #2) node{{\tiny ${#4}$}};
\draw (-0.5+#1,1+ #2) node{{\tiny ${#5}$}};
}

\def \Phishape #1#2#3#4#5#6
{
\Yshape{#1}{#2}{#3}{#4}{#5}
\begin{scope}[yshift = 2cm]
\draw (#1,1+ #2)--(#1,#2) -- (-1+ #1,-1 + #2);
\draw (#1,#2) -- (1+ #1,-1+ #2);
\draw (-0.5+#1,1+ #2) node{{\tiny ${#6}$}};
\end{scope}
}

\def \Ishape #1#2#3#4#5#6#7
{
\Yshape{#1}{2+#2}{#5}{#6}{#7}
\draw (#1,1+ #2)--(#1,#2) -- (-1+ #1,-1 + #2);
\draw (#1,#2) -- (1+ #1,-1+ #2);
\draw (-1.5+#1,-1+ #2) node{{\tiny ${#3}$}};
\draw (1.5+#1, -1+ #2) node{{\tiny ${#4}$}};
}

\def \ShortIshape #1#2#3#4#5#6#7
{
\draw (#1,1+ #2)--(#1,#2) -- (-1+ #1,-1 + #2);
\draw (#1,#2) -- (1+ #1,-1+ #2);
\draw (-1+#1, 2+#2) -- (#1, 1+#2) -- (1+#1, 2+#2);
\draw (-1.5+#1,-1+ #2) node{{\tiny ${#3}$}};
\draw (1.5+#1, -1+ #2) node{{\tiny ${#4}$}};
\draw (-1.5+#1, 2+ #2) node{{\tiny ${#5}$}};
\draw (1.5+#1, 2+ #2) node{{\tiny ${#6}$}};
\draw (-0.5+#1, 0.5+ #2) node{{\tiny ${#7}$}};
}

\def \IshapeL #1#2#3#4#5
{
\draw (#1,1+ #2)--(#1,#2) -- (-1+ #1,-1 + #2);
\draw (#1,#2) -- (1+ #1,-1+ #2);
\draw (-1+#1, 2+#2) -- (#1, 1+#2) -- (1+#1, 2+#2);
\draw (-0.8+#1, 0.5+ #2) node{{\tiny ${#3}$}};
\draw (#1,-0.3+#2) node{{\tiny $#4$}};
\draw (#1,1.3+#2) node{{\tiny $#5$}};
}

\def \IshapeR #1#2#3#4#5
{
\draw (#1,1+ #2)--(#1,#2) -- (-1+ #1,-1 + #2);
\draw (#1,#2) -- (1+ #1,-1+ #2);
\draw (-1+#1, 2+#2) -- (#1, 1+#2) -- (1+#1, 2+#2);
\draw (0.8+#1, 0.5+ #2) node{{\tiny ${#3}$}};
\draw (#1,-0.3+#2) node{{\tiny $#4$}};
\draw (#1,1.3+#2) node{{\tiny $#5$}};
}

\def \IshapeRP #1#2#3#4#5#6
{
\draw (#1,1+ #2)--(#1,#2) -- (-1+ #1,-1 + #2);
\draw (#1,#2) -- (1+ #1,-1+ #2);
\draw (-1+#1, 2+#2) -- (#1, 1+#2) -- (1+#1, 2+#2);
\draw (0.6+#1, 0.5+ #2) node{{\tiny ${#3}$}};
\draw [dotted] plot [smooth] coordinates {(-0.6+#1,-0.3+#2) (-0.9+#1,0.5+#2) (-0.6+#1,1.3+#2)};
\draw [dotted] plot [smooth] coordinates {(0.9+#1,-0.3+#2) (1.2+#1,0.5+#2) (0.9+#1,1.3+#2)};
\draw (-1.1+#1,1.1+#2) node{{\tiny $#6$}};
\draw (#1,-0.3+#2) node{{\tiny $#4$}};
\draw (#1,1.3+#2) node{{\tiny $#5$}};
}

\def \Vline #1#2#3#4
{
\draw (#1,#2) -- (#1,#2 + #3);
\draw (-0.5+ #1,#2) node{\tiny ${#4}$};
}

\def \VlineDot #1#2#3
{
\draw (#1, #2) -- (#1, #2 + #3);
}

\def \Braid#1#2
{
\draw (2+#1,#2) -- (#1,2+#2);
\draw (0.8+#1,0.8+#2)[color = white, line width = 2mm] -- (1.2+#1,1.2+#2);
\draw (#1,#2) -- (2+#1,2+#2);
}
\def \InvBraid#1#2
{
\draw (#1,#2) -- (2+#1,2+#2);
\draw (2+#1,#2)[color = white, line width = 2mm] -- (#1,2+#2);
\draw (2+#1,#2) -- (#1,2+#2);
}

\def \SmallBraid#1#2#3 
{
\draw (#3+#1,#2) -- (#1,#3+#2);
\draw (0.4*#3+#1,0.4*#3+#2)[color = white, line width = 2mm] -- (0.6*#3+#1,0.6*#3+#2);
\draw (#1,#2) -- (#3+#1,#3+#2);
}

\def \SmallInvBraid#1#2#3 
{
\draw (#1,#2) -- (#3+#1,#3+#2);
\draw (0.6*#3+#1,0.4*#3+#2)[color = white, line width = 2mm] -- (0.4*#3+#1,0.6*#3+#2);
\draw (#3+#1,#2) -- (#1,#3+#2);
}

\def\RecMor#1#2#3#4#5#6 
{
\draw[postaction={decorate}] (#1,0.5+#2) -- (#1,#2) node[left]{{\tiny $#3$}};
\draw[postaction={decorate}] (1+#1,0.5+#2) -- (1+#1,#2)node[right]{{\tiny $#4$}};

\draw[postaction={decorate}] (#1,2+#2)node[left]{{\tiny $#5$}} -- (#1,1.5+#2);
\draw[postaction={decorate}] (1+#1,2+#2)node[right]{{\tiny $#6$}} -- (1+#1,1.5+#2);

\draw (-0.5+#1,0.5+#2) rectangle (1.5+#1,1.5+#2);
}

\def \Segment#1#2#3#4#5
{
\tikz[scale = 0.5, baseline = {(current bounding box.center)}]{
 \draw (#1,#2) -- (2+ #1,#2);
 \draw (-0.5+#1,#2) node{{\tiny $#3$}};
 \draw (2.5+#1,#2) node{{\tiny $#4$}};
 \draw (1+#1,0.5+#2) node{{\tiny $#5$}};
}
}
\def \Triangle#1#2#3#4#5#6#7#8#9
{
\tikz[scale = 0.5, baseline= {(current bounding box.center)}]{
 \draw (#1,#2) -- (2+#1,#2) -- (1+#1,1.73+#2) -- (#1,#2);
 \draw (-0.5+#1,#2) node{{\tiny $#3$}};
 \draw (2.5+#1,#2) node{{\tiny $#4$}};
 \draw (1+#1,2.23+#2) node{{\tiny $#5$}};
 \draw (1+#1,-0.3+#2) node{{\tiny $#6$}};
 \draw (2+#1,1+#2) node{{\tiny $#7$}};
 \draw (0+#1,1+#2) node{{\tiny $#8$}};
 \draw (1+#1,0.5+#2) node{{\tiny $#9$}};
}
}
\title{Higher Categories and Topological Quantum Field Theories}

\author{Shawn X. Cui%
       \thanks{Email: \texttt{xingshan@vt.edu}}}
\affil{Stanford Institute for Theoretical Physics,\\ Stanford University, Stanford, California 94305, USA}
\date{ }

\theoremstyle{plain}
\newtheorem{theorem}{Theorem}[section]

\newtheorem{lemma}[theorem]{Lemma}
\newtheorem{proposition}[theorem]{Proposition}

\theoremstyle{definition}
\newtheorem{definition}[theorem]{Definition}

\theoremstyle{remark}
\newtheorem{remark} [theorem]{Remark}

\begin{document}
\maketitle
\begin{abstract}
We construct a state-sum type invariant of smooth closed oriented $4$-manifolds out of a $G$-crossed braided spherical fusion category ($G$-BSFC) for $G$ a finite group. The construction can be extended to obtain a $(3+1)$-dimensional topological quantum field theory (TQFT). The invariant of $4$-manifolds generalizes several known invariants in literature such as the Crane-Yetter invariant from a ribbon fusion category and Yetter's invariant from homotopy $2$-types. If the $G$-BSFC is concentrated only at the sector indexed by the trivial group element, a cohomology class in $H^4(G,U(1))$ can be introduced to produce a different invariant, which reduces to the twisted Dijkgraaf-Witten theory in a special case. Although not proven, it is believed that our invariants are strictly different from other known invariants. It remains to be seen if the invariants are sensitive to smooth structures. It is expected that the most general input to the state-sum type construction of $(3+1)$-TQFTs is a spherical fusion $2$-category. We show that a $G$-BSFC corresponds to a monoidal $2$-category with certain extra structure, but that structure does not satisfy all the axioms of a spherical fusion $2$-category given by M. Mackaay. Thus the question of what axioms properly define a spherical fusion $2$-category is open.
\end{abstract}

\tableofcontents 

\section{Introduction}

The notion of a topological quantum field theory ($\TQFT$) was invented by E. Witten based on path integrals in physics \cite{witten1988topological} and was given a mathematical formulation in terms of axioms by M. Atiyah \cite{atiyah1988topological} in the $1980$s. Since then there has been a vast study of $\TQFT$s both from the physics side and mathematics side. Throughout the paper, we work in the category of smooth oriented manifolds. Roughly speaking, for every positive integer $d$, a $(d+1)$-dimensional $\TQFT$ ($(d+1)$-$\TQFT$ for short) associates to every closed $d$-manifold a finite dimensional Hilbert space and to every $(d+1)$-manifold a vector in the Hilbert space corresponding to its boundary. These assignments should satisfy certain compatibility properties as specified by the axioms. The empty set is considered as a special closed $d$-manifold and the Hilbert space associated to it is required to be $\C$. Then a $(d+1)$-$\TQFT$ produces a complex scalar, called the partition function, for each closed $(d+1)$-manifold, and the scalar is an invariant of closed $(d+1)$-manifolds. This is an important application of $\TQFT$ to topology. 

The study of $\TQFT$s is closely related to higher category theory \cite{lurie2009higher, Baez:1995xq}. In general, a $(d+1)$-$\TQFT$ is to be described by the data of a $d$-category. On one hand, strict $d$-categories are well-defined for any $d$ \cite{eilenberg1966closed, crans1995combinatorial}, but this is insufficient for the purpose of $\TQFT$s since many important $d$-categories are not strict and can not be strictified either \cite{Baez1997introduction}. On the other hand, weak $d$-categories are only rigorously defined for small $d$ (such as $d = 1,2,3$)\footnote{A definition of weak 4-categories was given in \cite{trimble1995notes}, which was further clarified in \cite{hoffnung2011spans}}, and it is still controversial what should be the right notion of weak $d$-categories for higher $d$, although there have been many efforts in this direction \cite{street1987algebra, baez1998higher, tamsamani1996equivalence}. In the following, by $d$-categories we always mean weak $d$-categories. Special $d$-categories can be obtained from $k$-categories with certain extra structure for $k < n$. For instance, a monoidal $1$-category is a $2$-category and a braided monoidal $1$-category is a $3$-category \cite{kapanov1994category, gordon1995coherence, gurski2006}. Higher categories are natural resources for $\TQFT$s as shown below. We first give a brief overview of some categorical constructions of $(2+1)$- and $(3+1)$-$\TQFT$s.

There has been a fundamental achievement in $(2+1)$-$\TQFT$s which builds a nontrivial connection between monoidal categories, Hopf algebras, and $3$-manifolds. N. Reshetikhin and V. Turaev constructed an invariant of $3$-manifolds using modular tensor categories, which is believed to be the mathematical realization of E. Witten's $\TQFT$ from nonAbelian Chern-Simons theories \cite{reshetikhin1991invariants}. V. Turaev and O. Viro gave a state-sum invariant of $3$-manifolds (Turaev-Viro invariant) from the quantum $6j$ symbols of $U_q(sl_2)$ for $q$ a certain root of unity \cite{turaev1992state}. Later J. Barrett and B. Westbury generalized this construction (Turaev-Viro-Barrett-Westbury invariant, or TVBW for short) by using any spherical fusion category \cite{barrett1996invariants}. These invariants can all be extended to define a $(2+1)$-$\TQFT$ \footnote{The Reshetikhin-Turaev $\TQFT$ has an anomaly, but this is not the concern of this paper.}. Apart from the categorical constructions, another approach is by using Hopf algebras, among which the Kuperberg invariant \cite{kuperberg1996noninvolutory} and the Hennings invariant \cite{kauffman1995invariants, hennings1996invariants} are nonsemisimple generalizations of the TVBW invariant and the Reshetikhin-Turaev invariant, respectively. A special case of the Kuperberg invariant (and also the TVBW invariant) reduces to the Dijkgraaf-Witten theory \cite{dijkgraaf1990topological}. The study of $(2+1)$-$\TQFT$s has led to applications in quantum groups, $3d$ topology, and topological quantum computing. For example, the Turaev-Viro invariant can distinguish certain $3$-manifolds which are homotopy equivalent.

In one dimension higher, the theory of $(3+1)$-$\TQFT$s, however, is not understood as well as its counterpart in $(2+1)$ dimensions. The Dijkgraaf-Witten invariant \cite{dijkgraaf1990topological} in $(3+1)$ dimensions, as well as in other dimensions, measures the number of group morphisms from the fundamental group of a closed manifold to a given finite group. L. Crane and I. Frankel constructed a $4$-manifold invariant out of some algebraic structure, called Hopf categories \cite{crane1994four}. In \cite{crane1993categorical}, L. Crane and D. Yetter gave a state-sum invariant (Crane-Yetter invariant) using a semisimple subquotient of the category of representations of $U_q(sl_2)$ for $q$ a certain principal root of unity. L. Crane, L. Kauffman and D. Yetter generalized the construction to ribbon fusion categories \cite{crane1997state}\footnote{In \cite{crane1997state}, they were called semisimple tortile categories. The generalized invariant is still called the Crane-Yetter invariant.}, which reduces to the (untwisted) Dijkgraaf-Witten invariant for the category of finite dimensional representations of a finite group \cite{barenz2016dichromatic}. 
 The Walker-Wang model \cite{walker2012top} is believed to be a Hamiltonian realization of the Crane-Yetter invariant.
  The modular Crane-Yetter invariant, which is obtained from a modular tensor category, turns out to be a function of the Euler characteristics and the signature \cite{Crane93evaluatingthe}, and thus is a classical invariant. From a different direction, D. Yetter gave a construction of $(3+1)$-$\TQFT$ from homotopy $2$-types \cite{yetter1993tqft}, which is equivalent to a crossed module or a categorical group. Along a similar line, A. Kapustin \cite{kapustin2013higher} and M. Mackaay \cite{mackaay2000finite} obtained $4$-manifold invariants from $2$-groups with some additional structures. More recently, R. Kashaev produced a $(3+1)$-$\TQFT$ out of a cyclic group $\Z_N$ \cite{kashaev2014asimple}.

Topology in $4d$ is different from all other dimensions. $\R^4$ is the only Euclidean space accepting more than one smooth structure (infinitely many actually)\cite{gompf1985infinite}. After the work of M. Freedman \cite{freedman1990topology}, topological $4$-manifolds are fairly well understood. However, classifying smooth $4$-manifolds remains one of the most difficult open problems. Nontrivial $(3+1)$-$\TQFT$s are rare since any such $\TQFT$ gives rise to an invariant of smooth $4$-manifolds. The categorical constructions of invariants mentioned above are not known to be sensitive to smooth structures.

The main result of this paper is a state-sum construction of an invariant of $4$-manifolds from a \CatName{$G$} category (\CatNameShort{$G$}) where $G$ is a finite group. \CatNameShort{$G$}s were introduced by V. Turaev \cite{turaev2000homotopy} in the context of homotopy quantum field theories where they were called ribbon crossed $G$-categories. A. Kirillov later studied \CatNameShort{$G$}s in the context of conformal field theories where they were called $G$-equivariant categories \cite{Kirillov_ong-equivariant}.
\CatNameShort{$G$}s were also studied in \cite{etingof2010fusion, drinfeld2010braided, cui2015gauging} and have applications in condensed matter physics \cite{maissam2014symmetry}. Roughly, a \CatNameShort{$G$} $\cross{\Cat}{G}$ consists of the following structures (see Section \ref{sec:GBSFC} for a detailed definition).

\begin{enumerate}
 \item $\cross{\Cat}{G}$ is a spherical fusion category.
 \item $\cross{\Cat}{G} = \bigoplus\limits_{g \in G} \Cat_g$ is graded by $G$ where $\Cat_g$ is a full subcategory and $\Cat_{g_1} \otimes \Cat_{g_2} \subset \Cat_{g_1g_2}$, $g, g_1, g_2 \in G$.
 \item There is a $G$-action on $\cross{\Cat}{G}$ by pivotal functors such that for any $g, g' \in G$, the action of $g$ maps $\Cat_{g'}$ to $\Cat_{gg'g^{-1}}$.
 \item There is a $G$-crossed braiding
 \begin{equation*}
  \{c_{X, Y}: X \otimes Y \longrightarrow \lsup{g}{Y} \otimes X \ | \  X \in \Cat_g, Y \in \cross{\Cat}{G}\},
 \end{equation*}
 which satisfies certain compatibility conditions.
\end{enumerate}

Given a \CatNameShort{$G$} $\cross{\Cat}{G}$ and a closed oriented $4$-manifold $M$, the procedure of constructing the invariant of $M$ from $\cross{\Cat}{G}$ goes as follows. Take a triangulation $\T$ of $M$. A coloring $\ColorS$ assigns an element of $G$ to each $1$-simplex and a simple object of $\cross{\Cat}{G}$ to each $2$- and $3$-simplex. These assignments are subject to certain constraints. To avoid distraction from technical details here, we simply present the invariant in the form,
\begin{equation}
\label{equ:formal_partition}
Z_{\cross{\Cat}{G}}(M,\T) = \sum\limits_{\mathcal{S}} \frac{\prod\limits_{\Delta_{4} \in \T^4} Z_{\mathcal{S}}(\Delta_{4})\prod\limits_{\Delta_{2}\in \T^2}Z_{\mathcal{S}}(\Delta_{2})\prod\limits_{\Delta_{0}\in \T^0}Z_{\mathcal{S}}(\Delta_{0})}
{\prod\limits_{\Delta_{3}\in \T^3}Z_{\mathcal{S}}(\Delta_{3})\prod\limits_{\Delta_{1}\in \T^1}Z_{\mathcal{S}}(\Delta_{1})},
\end{equation}
where $\T^i$ is the set of $i$-simplices and $Z_{\mathcal{S}}(\Delta_i)$ is some factor associated with the $i$-simplex $\Delta_i$. Thus for each coloring $\mathcal{S}$, every simplex in $\T$ contributes a factor to the invariant. Usually, the contributions from the top dimensional simplices are the most important ones. See Section \ref{sec:partition} for more details. The following is the main theorem of the paper.

\begin{theorem}[Main, informal]
The formula for $Z_{\cross{\Cat}{G}}(M, \T)$ is independent of the choice of the triangulation $\T$ and thus $Z_{\cross{\Cat}{G}}(M):= Z_{\cross{\Cat}{G}}(M, \T)$ is an invariant of closed smooth oriented $4$-manifolds.
\end{theorem}

Moreover, the construction can be extended to obtain a $(3+1)$-$\TQFT$. The following theorem shows that the $4$-manifold invariant constructed here is rather rich and it  generalizes several known categorical invariants in literature.

\begin{theorem}
\begin{itemize}
\item If $G$ is the trivial group, then $\cross{\Cat}{G}$ is a ribbon fusion category and $Z_{\cross{\Cat}{G}}(M)$ is equal to the Crane-Yetter invariant of $M$ from $\cross{\Cat}{G}$.
\item If all the simple objects of $\cross{\Cat}{G}$ are invertible, then $\cross{\Cat}{G}$ corresponds to a crossed module (or equivalently a categorical group), and $Z_{\cross{\Cat}{G}}(M)$ is equal to Yetter's invariant from the corresponding categorical group.
\item If the $G$-action is trivial and $\Cat_g = 0$ for all $g \neq e \in G$, then again $\cross{\Cat}{G}$ is a ribbon fusion category, and $Z_{\cross{\Cat}{G}}(M)$ is the product of the untwisted Dijkgraaf-Witten invariant from $G$ and the Crane-Yetter invariant from $\cross{\Cat}{G}$.
\end{itemize}
\end{theorem}

Further generalizations of our invariant are also possible. In particular, we show that if $\cross{\Cat}{G} = \bigoplus\limits_{g \in G}\Cat_g$ is a \CatNameShort{$G$} such that $\Cat_g = 0$ for all $g \neq e$, then we can introduce a $4$-cohomology class $\omega \in H^4(G, U(1))$ and the new partition function takes the same form as Equation \ref{equ:formal_partition} except for each $4$-simplex $\Delta = (ijklm)$, $Z_{\mathcal{S}}(\Delta_{4})$ is replaced with,
\begin{align*}
Z_{\mathcal{S}}(\Delta_{4})' &= Z_{\mathcal{S}}(\Delta_{4})\,\omega\left(\mathcal{S}(ij),\mathcal{S}(jk),\mathcal{S}(kl),\mathcal{S}(lm)\right).
\end{align*}
Denote the new partition function by $Z_{\cross{\Cat}{G}, \omega}(M,\T)$.

\begin{theorem}
$Z_{\cross{\Cat}{G}, \omega}(M,\T)$ is again independent of the choice of the triangulation $\T$, and is thus an invariant of closed smooth oriented $4$-manifolds.
\end{theorem}

As a special case, we have the the following proposition.
\begin{proposition}
If the $G$-action is trivial and $\Cat_g = 0$ for all $g \neq e \in G$, then $Z_{\cross{\Cat}{G},\omega}(M,\T)$ is the product of the $\omega$-twisted Dijkgraaf-Witten invariant from $G$ and the Crane-Yetter invariant from $\cross{\Cat}{G}$.
\end{proposition}

From a different perspective, in $(2+1)$ dimensions, the typical state-sum model (Turaev-Viro-Barrett-Westbury invariant) involves a spherical fusion ($1$-)category. Thus in $(3+1)$ dimensions, one would expect there to be a notion of a \lq spherical fusion $2$-category' which results in the most general construction of state-sum $(3+1)$-$\TQFT$s and all known categorical invariants $(\TQFT \text{s})$ of state-sum type should fit into this framework. Monoidal $2$-categories are defined in \cite{kapanov1994category} and monoidal $2$-categories with duals are given in \cite{Baez2003705}. In \cite{mackaay1999spherical}, M. Mackaay proposed a definition of spherical fusion $2$-categories\footnote{In \cite{mackaay1999spherical} they are called nondegenerate finitely semisimple semistrict spherical 2-categories of nonzero dimension.} as monoidal $2$-categories with some additional structures. Based on his definition, he formally defined a $4$-manifold invariant. However, as explained later, his definition is too restrictive and excludes many interesting examples. In particular, his construction does not cover the invariant from a \CatNameShort{$G$} due to the following proposition.

\begin{proposition}[Informal]
A monoidal $2$-category $\D(\cross{\Cat}{G})$ with duals can be constructed from a \CatNameShort{$G$} $\cross{\Cat}{G}$, but $\D(\cross{\Cat}{G})$ does not satisfy the axioms of a spherical fusion $2$-category according to the definition in \cite{mackaay1999spherical}.
\end{proposition}

On the other hand, in \cite{barrett2012gray} spherical Gray categories are defined where a Gray category is a semistrict $3$-category (tricategory). A semistrict monoidal $2$-category can be viewed as a Gray category with one object. It will be shown that $\D(\cross{\Cat}{G})$ does become a spherical $2$-category (or rather, a spherical Gray category) if we adapt the definition of sphericity for Gray categories to monoidal $2$-categories. It is not clear though that a spherical $2$-category in the sense of \cite{barrett2012gray} leads to a $(3+1)$-$\TQFT$. We leave this as a future direction of study.

Lastly, $\GBSFC$s are not rare. In \cite{drinfeld2010braided, kirillov2002modular, kirillov2002modular2}, it has been proved that equivalence classes of $\GBSFC$s are in one-to-one correspondence, by equivariantization and de-equivariantization, with equivalence classes of spherical braided fusion categories containing $\text{Rep}(G)$ as a subcategory. Also, given a group morphism from $G$ to the group of automorphisms of a unitary braided fusion category $\Cat$, if certain obstructions vanish, then $\Cat$ can be extended to a unitary $G$-crossed braided fusion category, which is also a $\GBSFC$, with $\Cat$ as the sector indexed by the trivial group element \cite{etingof2010fusion}.

The structure of the paper is organized as follows. In Section \ref{sec:GBSFC}, we give a review of \CatNameShort{$G$}s. A \CatNameShort{$G$} can be understood either by embedding it into a strict \CatNameShort{$G$} (Section \ref{subsec:stratification}) or by extracting from it a set of discrete data satisfying certain equations (\ref{subsec:skeletonization}). Section \ref{sec:partition} is the core of the paper where three equivalent definitions of the invariant are given and the main theorem is stated. In Section \ref{sec:example}, we give several examples of the invariants and also introduce a variation of the invariants. Section \ref{sec:proofmain} contains the proof of the main theorem. In Section \ref{sec:twocat} we show that a monoidal $2$-category with certain extra structure can be constructed from a \CatNameShort{$G$}. Finally in Section \ref{sec:future} we provide some open questions for research.

\section{$G$-crossed Braided Spherical Fusion Categories (\CatNameShort{$G$}s)}
\label{sec:GBSFC}
We assume the readers are familiar with the concepts tensor (monoidal) categories and their various specializations such as fusion categories, spherical categories, and ribbon categories. Some knowledge of functors and natural transformations is also expected. We use the word \lq monoidal' and \lq tensor' interchangeably. There are a number of excellent references developing these concepts. See for instance \cite{kassel2012quantum, Etingof_finitetensor, bojko2001lectures, Etingof_onfusion, Wang2010topological}. Since our main object to use is a \CatName{$G$} category (\CatNameShort{$G$}) where $G$ is a finite group, we give a review of such categories. See \cite{Kirillov_ong-equivariant, turaev2000homotopy} for more detailed discussions. But notice that these categories are called \CatNameTu{$G$} categories in \cite{turaev2000homotopy} and \CatNameKi{$G$} categories in \cite{Kirillov_ong-equivariant}. Some conventions from this paper are also different from the references.

If $\Cat$ is a category, denote the set of objects by $\Cat^0$ and the set of morphisms by $\Cat^1$. We follow the convention that compositions in a category are read from right to left, namely, if $f \in \Hom(X, Y), g \in \Hom(Y, Z),$ then $g \circ f \in \Hom(X, Z)$. The identity map on an object $X$ is denoted by $I_X$, or $Id_X$. The subscript will often be dropped. Throughout the paper, $G$ denotes a finite group.

\subsection{Definition of \CatNameShort{$G$}s}
\label{subsec:GBSFC_Def}

Let $(\Cat, \otimes, \unit, a, l, r)$ be a tensor category, where $a,\  l,\  r$ are the structure isomorphisms,
\begin{align*}
\begin{array}{rrcl}
a_{X,Y,Z}\ \colon& (X \otimes Y) \otimes Z & \overset{\simeq}{\longrightarrow} & X \otimes (Y \otimes Z), \\
l_X\ \colon &\unit \otimes X & \overset{\simeq}{\longrightarrow}&  X,\\
r_X\ \colon & X \otimes \unit& \overset{\simeq}{\longrightarrow}& X,
\end{array}
\end{align*}
which satisfy the Pentagon Identity and Triangle Identity. When no confusion arises, we often drop the structure symbols and claim $\Cat$ as a tensor category.

Denote by $\Aut{\Cat}$ the tensor category where objects are tensor auto-equivalences of $\Cat$, morphisms are natural transformations, and the tensor product of two tensor equivalences is the composition of functors. Denote by $\uG{G}$ the tensor category where the objects are elements of $G$ and the tensor product is given by group multiplication. There is only one morphism, the identity map, from an object to itself, and no morphism between different objects.

\begin{definition}
\label{def:G-action}
For a tensor category $\Cat$, a $G$-action on $\Cat$ is a tensor functor $ \uG{G} \longrightarrow \Aut{\Cat}$.
\end{definition}

The above is a compact way to describe a $G$-action. To be more clear, we unpack the definition into a set of specific axioms. Let
\begin{align*}
\begin{array}{rrcl}
(F, \eta, \epsilon)\ \colon &\uG{G} &\longrightarrow & \Aut{\Cat}
\end{array}
\end{align*}
be a tensor functor. For $g \in G,\  X \in \Cat^0, \ f \in \Cat^1$, we often write $F(g)(X)$ and $F(g)(f)$ as $\lsup{g}{X}$ and $\lsup{g}{f}$, respectively. For each object $g \in \uG{G}^0$, i.e., an element $g \in G$, $F(g)$ is a tensor auto-equivalence and hence it is endowed with natural isomorphisms,
\begin{align*}
\begin{array}{rrcl}
\gamma_{g; X,Y}\ \colon& \lsup{g}{(X \otimes Y)}& \overset{\simeq}{\longrightarrow}& \lsup{g}{X} \otimes \lsup{g}{Y}, \\
\sigma_{g}\ \colon& \lsup{g}{\unit}& \overset{\simeq}{\longrightarrow}& \unit,
\end{array}
\end{align*}
which are compatible with the tensor structure of $\Cat$. Since $(F, \eta, \epsilon)$ is a tensor functor, there are natural isomorphisms $\eta_{g,h}\colon F(gh) \longrightarrow F(g) \circ F(h)$ and  $\epsilon\colon F(e) \longrightarrow Id_{\Cat}$. Or more specifically, for each $X \in \Cat^0$, there are isomorphisms $\eta_{g,h;X}\colon \lsup{gh}{X} \longrightarrow \lsup{g}{(\lsup{h}{X})}$ and $\epsilon_{X}\colon \lsup{e}{X} \longrightarrow X$, such that the following diagrams commute.
\begin{equation}
\insertimage{1}{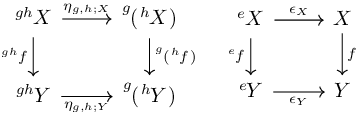}
\end{equation}

\begin{equation}
\insertimage{1}{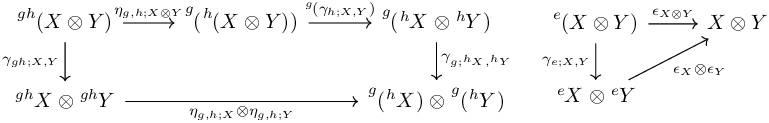}
\end{equation}

\begin{equation}
\insertimage{1}{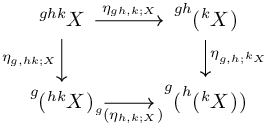}
\end{equation}

\begin{equation}
\insertimage{1}{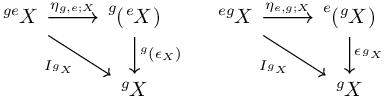}
\end{equation}


\begin{definition}
A \CatName{$G$} category (\CatNameShort{$G$}) is a $G$-graded spherical fusion category $\cross{\Cat}{G} = \bigoplus\limits_{g \in G} \Cat_g$ together with a $G$-action and a family of natural isomorphisms
\begin{align}
\{c_{X, Y}\colon X \otimes Y \longrightarrow \lsup{g}{Y} \otimes X \ | \  X \in \Cat_g, Y \in \cross{\Cat}{G}\},
\end{align}
such that the following conditions are satisfied.
\begin{enumerate}
  \item Each $\Cat_g$, called the $g$-sector, is a full subcategory. The only morphism between two objects from different sectors is the zero morphism. Moreover, $\Cat_g \otimes \Cat_{g'} \subset \Cat_{gg'}$ and $\lsup{g}{(\Cat_{g'})} \subset \Cat_{gg'g^{-1}}$.
  \item The $c_{X, Y}\,'$s, called $G$-crossed braidings, are self consistent, namely, for $X \in \Cat_g,\  Y \in \Cat_h,\  Z \in \Cat_{k}$, the following diagrams commute,
    \begin{equation}
    \insertimage{1}{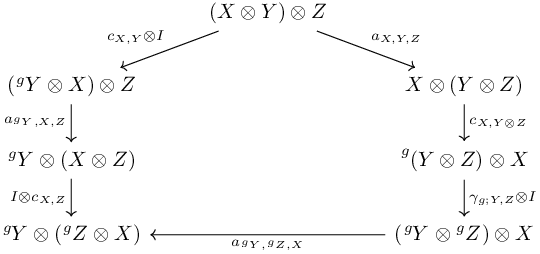}
\end{equation}
\begin{equation}
\insertimage{1}{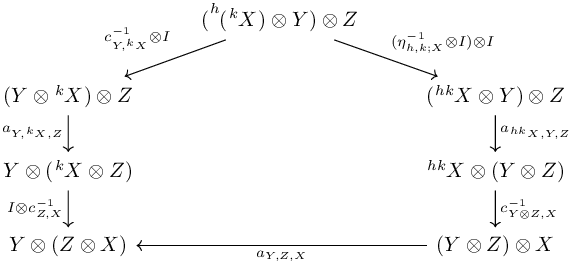}
\end{equation}
 \item The isomorphism $c_{X, Y}$ is natural with respect to both $X$ and $Y$. That is, for $X, \ X' \in \Cat_g,\  Y,\  Y' \in \Cat_{g'},\ \phi \in \Hom(X, X'),\ \psi \in \Hom(Y, Y'),$ we have,
 \begin{equation}
 \insertimage{1}{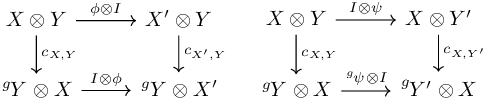}
\end{equation}

  \item The $G$-action is consistent with the $G$-crossed braiding. Namely, for $X \in \Cat_g, \ Y \in \cross{\Cat}{G}$, the following diagram commutes.
  \begin{equation}
  \insertimage{1}{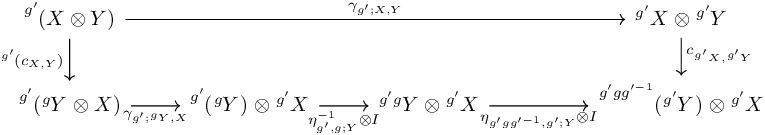}
\end{equation}
\item The $G$-action is consistent with the pivotal structure. Thus, if $\delta_{X}: X \longrightarrow X^{**}$ is the pivotal structure, then the following diagram commutes.
    \begin{equation}
    \insertimage{1}{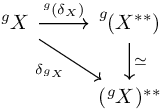}
    \end{equation}
    The vertical arrow above represents the canonical isomorphism induced by the action of $g$.
\end{enumerate}

\end{definition}

Note that if $G$ is trivial, then the definition above is the same as that of spherical braided fusion categories, so an \CatNameShort{$\{e\}$} is simply a spherical braided fusion category, or a ribbon fusion category. It should be noted that in general a \CatNameShort{$G$} is not a braided tensor category. However, the sector $\Cat_e$ is indeed always a ribbon fusion category. We also do not require the grading to be faithful. For instance, one can take an arbitrary ribbon fusion category $\Cat$ and set $\Cat_e = \Cat, \  \Cat_g = 0,\  g \neq e$, then $\cross{\Cat}{G} = \bigoplus\limits_{g \in G} \Cat_g = \Cat$ is a \CatNameShort{$G$} with the trivial $G$-action. More interesting examples are provided in Section \ref{sec:example}.

There are two opposite directions to study a \CatNameShort{$G$}. One direction is strictifying a \CatNameShort{$G$}, where one shows that a \CatNameShort{$G$} is equivalent, in some properly defined sense, to a \CatNameShort{$G$} in which the structure isomorphisms $a_{X, Y, Z},\, l_X,\, r_X,\, \gamma_{g;X,Y},\, \sigma_{g},\, \eta_{g,h;X},\, \epsilon_{X}$ are all identity maps. A \CatNameShort{$G$} is called strict if it satisfies the above properties. The advantage of using strict \CatNameShort{$G$}s is that it is easy to deduce identities of morphisms and moreover, identities which hold in a strict category also hold in an equivalent but nonstrict category after appropriately inserting certain structure isomorphisms. The other direction is skeletonizing a \CatNameShort{$G$}, where a representative for each isomorphism class of simple objects is chosen, and the category is described by a set of discrete data which satisfy some equations involving the representatives. This method is useful when one needs to perform calculations in a category. For instance, in Section \ref{sec:partition}, the discrete data will be used to compute the invariant of a $4$-manifold. We elaborate these two concepts in Section \ref{subsec:stratification} and Section \ref{subsec:skeletonization}, respectively.

\subsection{Strictifying \CatNameShort{$G$}s}
\label{subsec:stratification}

A \CatNameShort{$G$} is called strict, if it is a strict spherical fusion category and the structure isomorphisms $\gamma_{g;X,Y}, \,\sigma_{g},\, \eta_{g,h;X}, \,\epsilon_{X}$ are the identity maps. Every \CatNameShort{$G$} is equivalent to a strict one as indicated in the following theorem. Roughly speaking, an equivalence of \CatNameShort{$G$}s is an equivalence as a tensor functor that preserves all additional structures, e.g., crossed braiding, $G$-action, etc. In particular, we require such an equivalence to preserve the $G$-grading. That is, it maps the $g$-sector to the $g$-sector. For a rigorous definition, see \cite{muger2010onthe}.

\begin{theorem}\cite{muger2010onthe}
Let $\cross{\Cat}{G}$ be a \CatNameShort{$G$}, then there exists a strict \CatNameShort{$G$} $\D$ and an equivalence $F\colon \cross{\Cat}{G} \longrightarrow \D$ of \CatNameShort{$G$}s.
\end{theorem}

In a strict \CatNameShort{$G$}, it is convenient to represent morphisms with colored graph diagrams. We only list the basic rules for this representation. For detailed treatment, see \cite{turaev1994quantum}. But note that we will follow a slightly different convention.

A graph diagram is a collection of rectangles \footnote{In \cite{turaev1994quantum}, they are called coupons.}, immersed segments, and immersed circles in $\R \times [0,1]$. The segments and circles are called $1$-strata of the diagram. The following conditions are required to be satisfied.
\begin{itemize}
\item All rectangles are disjoint from each other and lie in $\R \times (0,1)$. We assume that  the longer sides of each rectangle are parallel to $\R \times \{0\}$ (i.e., horizontal) and the shorter sides vertical to $\R \times \{0\}$.
\item The $1$-strata end either on $\R \times \{0,1\}$ or on the horizontal sides of a rectangle. The interior of the $1$-strata lies in $\R \times (0,1)$ and does not intersect with any rectangles.
\item The $1$-strata are directed and may only have double crossings in $\R \times (0,1)$ with overcrossing/undercrossing data.
\end{itemize}

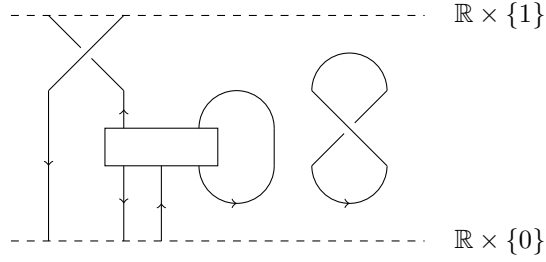
\begin{figure}
\centering
\begin{tikzpicture}[scale = 0.5]
\begin{scope}[decoration={
    markings,
    mark=at position 0.5 with {\arrow{>}}}]
 \draw [postaction={decorate}] (0,4) -- (0,0);
 \draw [postaction={decorate}] (2,2) -- (2,0);
  \draw [postaction={decorate}] (3,0) -- (3,2);
   \draw (1.5,2) rectangle (4.5, 3) ;
   \draw [postaction={decorate}] (4,2) arc(-180:0:1cm) ;
   \draw (4,3) arc(180:0:1cm);
   \draw (6,2) -- (6,3);
   \draw [postaction={decorate}](2,3) -- (2,4);
   \Braid{0}{4}
   \draw [postaction={decorate}](7, 2) arc(-180:0:1cm);
   \InvBraid{7}{2}
   \draw (7,4) arc(180:0:1cm);
   \draw[dashed] (-1,0) -- (10,0);
   \draw (12, 0)  node{$\R \times \{0\}$};
   \draw[dashed] (-1, 6) -- (10,6) ;
   \draw (12, 6) node{$\R \times \{1\}$};
 \end{scope}
\end{tikzpicture}
\caption{Example of a graph diagram}\label{fig:graphexample}
\end{figure}

\noindent See Figure \ref{fig:graphexample} for an example of a graph diagram. Usually we will not draw the dashed lines representing $\R \times \{0,1\}$ explicitly and assume that the bottom (resp. the top) of a graph diagram is bounded by the line $\R \times \{0\}$ (resp. $\R \times \{1\}$).

Given a graph diagram $\G$, we think of the $1$-strata broken at the undercrossings, namely, the $1$-strata decomposing into a collection of arcs, where each arc starts and ends either at an undercrossing or at one of the end points of the $1$-strata. For instance, the graph diagram in Figure \ref{fig:graphexample} contains seven arcs. Denote by $\G^0,\ \G^1,$ and $\G^2$ the set of crossings, arcs, and rectangles, respectively.

If $A \in \Cat_g$ is an object, denote $|A| = g$. A $\cross{\Cat}{G}$-coloring of a graph diagram $\G$ is an assignment $\left(\{\psi_{\alpha}\colon\ \alpha  \in \G^0 \}, \ 
\{A_{\beta}\colon\  \beta \in \G^1 \}, \ \{f_{\gamma}\colon\  \gamma \in \G^2\}\right)$ where the $\psi_{\alpha}\,'$s and $f_{\gamma}\,'$s are morphisms in $\cross{\Cat}{G}$ and the $A_{\beta}\,'$s are homogeneous objects of $\cross{\Cat}{G}$, such that the following conditions are satisfied.
\begin{itemize}
\item For each rectangle $\gamma$, denote by $\beta_1, \cdots, \beta_{m}$ the set of arcs incident to the bottom of $\gamma$ and by $\beta^1, \cdots, \beta^n$ the set of arcs incident to the top of $\gamma$, both enumerated from left to right. Note that some of the $\beta_i\,'$s and $\beta_j\,'$s could be the same arc. For each $\beta_i$ (resp. $\beta^j$), define $\epsilon_i$ (resp. $\epsilon^j$) to be $+1$ if $\beta_i$ (resp. $\beta^j$) is directed downwards near $\gamma$, and $-1$ otherwise. Let
\begin{align*}
A_{\gamma}:= \bigotimes\limits_{i=1}^{m} A_{\beta_i}^{\epsilon_i}, \qquad A^{\gamma}:= \bigotimes\limits_{j=1}^{n} A_{\beta^j}^{\epsilon^j},
\end{align*}
where for an object $B$, $B^{+1}:= B$ and $B^{-1}:= B^*$. If $m=0$, define $A_{\gamma}:= \unit$. Similarly if $n=0$, define $A^{\gamma} := \unit.$ Then we require $f_{\gamma}$ to be a morphism in $\Hom(A_{\gamma}, A^{\gamma})$. For instance, in Figure \ref{fig:rectangle} $(a)$, we represent the coloring by placing an object beside each arc and a morphism inside each rectangle. Then $f \in \Hom(A_1 \otimes A_2 \otimes A_3^*,\, B_1^* \otimes B_2)$.
\item There are two types of crossings, the positive crossing (Figure \ref{fig:rectangle} $(b)$) and the negative crossing (Figure \ref{fig:rectangle} $(c)$). Again, we represent the coloring of a crossing by placing a morphism beside it. Let the arc corresponding to the overcrossing, the arc entering the undercrossing, and the arc leaving the undercrossing be colored by the objects $A$, $B_1$, and $B_2$, respectively. In the case of a positive crossing, we require $\psi$ to be an isomorphism $\lsup{|A|}{B_2} \overset{\thicksim}{\longrightarrow} B_1$, and in the case of a negative crossing, we require $\psi$ to be an isomorphism $\lsup{|A|}{B_1} \overset{\thicksim}{\longrightarrow} B_2 $. We use the convention that if $\psi = I$, then we drop it from the diagram.
\end{itemize}

\begin{figure}
 \centering
\begin{tikzpicture}[scale = 0.5]
 \begin{scope}[decoration={
    markings,
    mark=at position 0.5 with {\arrow{>}}}]
    \draw [postaction={decorate}] (0,1) -- (0,0) node[left]{$A_1$};
    \draw [postaction={decorate}] (1.5,1) -- (1.5,0) node[left]{$A_{2}$};
    \draw [postaction={decorate}] (3,0)node[right]{$A_3$} -- (3,1) ;
    \draw [postaction={decorate}] (0,2) -- (0,3) node[left]{$B_1$};
    \draw [postaction={decorate}] (3,3)node[right]{$B_2$} -- (3,2) ;
    \draw (-0.5,1) rectangle (3.5,2)node[pos = 0.5]{$f$};

    \draw (1.5,-1) node{$(a)$};
    \end{scope}

    \begin{scope}[xshift = 7cm, decoration={
    markings,
    mark=at position 0.75 with {\arrow{>}}}]
    \draw [postaction={decorate}] (0,2)node[left]{$B_1$} -- (2,0)node[right]{$B_2$};
    \draw (2,2)[color = white, line width = 2mm] -- (0,0);
    \draw [postaction={decorate}] (2,2) -- (0,0)node[left]{$A$};
    \draw (1.6,1.1) node{$\psi$};
    \draw (1,-1) node{$(b)$};
    \end{scope}

     \begin{scope}[xshift = 13cm, decoration={
    markings,
    mark=at position 0.75 with {\arrow{>}}}]
    \draw [postaction={decorate}] (2,2)node[right]{$B_1$} -- (0,0)node[left]{$B_2$};
    \draw (0,2)[color = white, line width = 2mm] -- (2,0);
    \draw [postaction={decorate}] (0,2) -- (2,0)node[right]{$A$};
    \draw (1.6,1.1) node{$\psi$};
    \draw (1,-1) node{$(c)$};
    \end{scope}

 \end{tikzpicture}
 \caption{Colorings of a graph diagram}\label{fig:rectangle}
\end{figure}
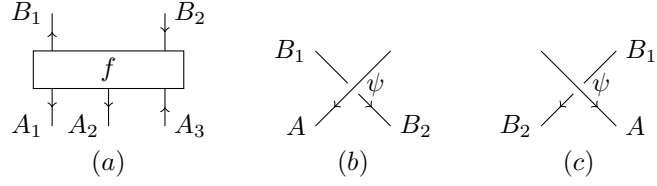
\begin{remark}
The motivation for labeling crossings is as follows. Consider the positive crossing in Figure \ref{fig:rectangle} $(b)$. As shall be illustrated below, the crossing is to be interpreted as the $G$-crossed braiding $c_{A, B_2}$,
\begin{align*}
\begin{array}{rrcl}
c_{A, B_2}\ \colon & A \otimes B_2 & \longrightarrow & \lsup{|A|}{B_2} \otimes A.
\end{array}
\end{align*}
That means we need to have $\lsup{|A|}{B_2} = B_1$. However, as a common principle in categories, it is more natural to require two objects to be {\it isomorphic} rather than {\it equal} on the nose. Thus we endow the crossing with an isomorphism $\psi\colon \lsup{|A|}{B_2} \longrightarrow B_1$.
\end{remark}

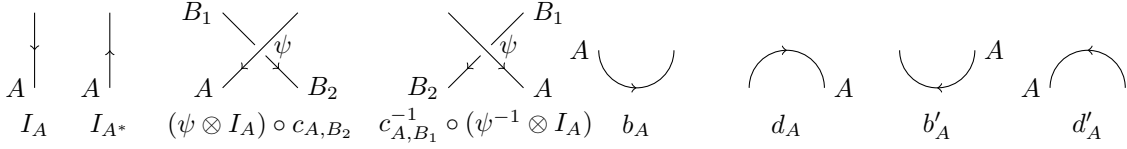
\begin{figure}
 \centering
\begin{tikzpicture}[scale = 0.5]
 \begin{scope}[decoration={
    markings,
    mark=at position 0.5 with {\arrow{>}}}]
   \draw [postaction={decorate}] (0,2) -- (0,0) node[left]{$A$};
   \draw (0,-1) node{$I_A$};

   \draw [postaction={decorate}] (2,0)node[left]{$A$} -- (2,2) ;
   \draw (2,-1) node{$I_{A^*}$};
    \end{scope}

    \begin{scope}[xshift = 5cm, decoration={
    markings,
    mark=at position 0.75 with {\arrow{>}}}]
    \draw [postaction={decorate}] (0,2)node[left]{$B_1$} -- (2,0)node[right]{$B_2$};
    \draw (2,2)[color = white, line width = 2mm] -- (0,0);
    \draw [postaction={decorate}] (2,2) -- (0,0)node[left]{$A$};
    \draw (1.6,1.1) node{$\psi$};
    \draw (1,-1) node{$(\psi \otimes I_A)\circ c_{A, B_2}$};
    \end{scope}

     \begin{scope}[xshift = 11cm, decoration={
    markings,
    mark=at position 0.75 with {\arrow{>}}}]
    \draw [postaction={decorate}] (2,2)node[right]{$B_1$} -- (0,0)node[left]{$B_2$};
    \draw (0,2)[color = white, line width = 2mm] -- (2,0);
    \draw [postaction={decorate}] (0,2) -- (2,0)node[right]{$A$};
    \draw (1.6,1.1) node{$\psi$};
    \draw (1,-1) node{$c_{A, B_1}^{-1}\circ (\psi^{-1} \otimes I_A)$};
    \end{scope}

    \begin{scope}[xshift = 15cm, decoration={
    markings,
    mark=at position 0.5 with {\arrow{>}}}]
 \draw [postaction={decorate}] (0,1) node[left]{$A$} arc(-180:0:1cm) ;
 \draw (1,-1) node{$b_A$};
 \end{scope}

  \begin{scope}[xshift = 19cm, decoration={
    markings,
    mark=at position 0.5 with {\arrow{>}}}]
 \draw [postaction={decorate}]   (0,0)  arc(180:0:1cm) node[right]{$A$};
 \draw (1,-1) node{$d_A$};
 \end{scope}

 \begin{scope}[xshift = 23cm, decoration={
    markings,
    mark=at position 0.5 with {\arrow{>}}}]
 \draw [postaction={decorate}]   (2,1) node[right]{$A$} arc(0:-180:1cm);
 \draw (1,-1) node{$b_A'$};
 \end{scope}

  \begin{scope}[xshift = 27cm, decoration={
    markings,
    mark=at position 0.5 with {\arrow{>}}}]
 \draw [postaction={decorate}]   (2,0)  arc(0:180:1cm) node[left]{$A$};
 \draw (1,-1) node{$d_A'$};
 \end{scope}

 \end{tikzpicture}
 \caption{Interpretations of a graph diagram}\label{fig:interpretation}
\end{figure}

The second condition on the coloring implies that $|A| \, |B_2| = |B_1| \, |A|$ for a positive crossing and $|B_2| \, |A| = |A| \, |B_1|$ for a negative crossing. We will think of the diagram $\G$ as living in $\R \times [0,1] \times \R$ with the identification $\R \times [0,1]$ with $\R \times [0,1] \times \{0\}$, where the positive $z$-axis (the third axis) points outward of the plane $\R \times [0,1]$. We then push the $1$-strata of $\G$ near an undercrossing slightly into the plane.  Denote this deformed graph by $\tilde{\G}$. It is not hard to see that a coloring of $\G$ determines a group morphism from $\pi_1(\R \times [0,1] \times \R \setminus \tilde{\G}, pt)$ to $G$ where the base point $pt$ can be chosen to be any point with a large $z$ coordinate. Given a coloring of $\G$, let $\beta_1, \cdots, \beta_m$ (resp. $\beta^1, \cdots, \beta^n $) be the set of arcs intersecting with $\R \times \{0\}$ (resp. $\R \times \{1\}$) listed from left to right, and similarly define $\epsilon_i, \epsilon^j$ for each arc $\beta_i, \beta^j$ as in the definition of a coloring. Let

\begin{align*}
A_{\G}:= \bigotimes\limits_{i=1}^{m} A_{\beta_i}^{\epsilon_i}, \qquad A^{\G}:= \bigotimes\limits_{j=1}^{n} A_{\beta^j}^{\epsilon^j}.
\end{align*}
And again let $A_{\G} = \unit$ if $m=0$ and $A^{\G} = \unit$ if $n=0$. By \cite{turaev1994quantum}, a $\cross{\Cat}{G}$-colored graph diagram $\G$ can be interpreted as a morphism of $\cross{\Cat}{G}$ in $\Hom(A_{\G}, A^{\G})$. The rules of the interpretation are as follows.
\begin{itemize}
 \item If $\G$ is one of the graph diagrams listed in Figure \ref{fig:interpretation}, then it is interpreted as the morphism shown below the corresponding diagram.
 \item If $\G$ contains a single rectangle and a set of arcs each of which is either a vertical segment connecting $\R \times \{0\}$ and the lower side of the rectangle or a vertical segment connecting $\R \times \{1\}$ and the upper side of the rectangle, (see Figure \ref{fig:rectangle} $(a)$ for instance), then it is interpreted as the coloring labeling the rectangle.
 \item Stacking a graph diagram on top of another corresponds to the composition of the morphisms they represent, and juxtaposition of diagrams corresponds to the tensor product of morphisms.
\end{itemize}

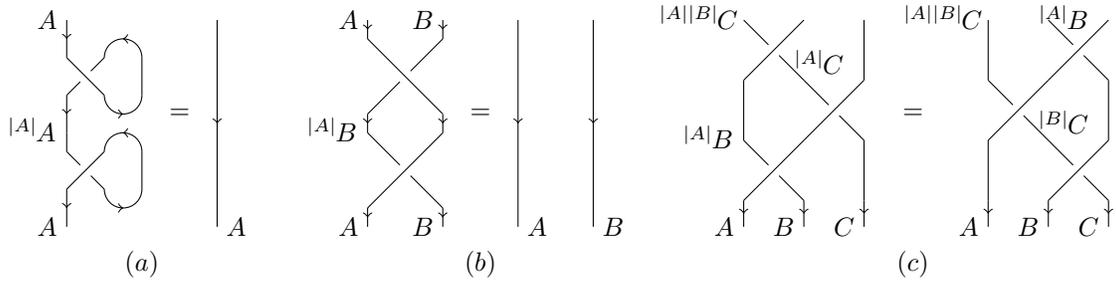
\begin{figure}
 \centering
\begin{tikzpicture}[scale = 0.5]
 \begin{scope}[decoration={
    markings,
    mark=at position 0.5 with {\arrow{>}}}]
   \draw [postaction={decorate}](0,1) -- (0,0) node[left]{$A$};
   \SmallBraid{0}{1}{1}
   \draw [postaction={decorate}](1,1) arc(-180:0:0.5cm);
   \draw [postaction={decorate}](2,2) arc(0:180:0.5cm);
   \draw (2,1) -- (2,2);
   \draw (0,2) -- (0,2.5);

   \begin{scope}[yshift = 2.5cm]
   \draw [postaction={decorate}](0,1) -- (0,0) node[left]{$\lsup{|A|}{A}$};
   \SmallInvBraid{0}{1}{1}
   \draw [postaction={decorate}](1,1) arc(-180:0:0.5cm);
   \draw [postaction={decorate}](2,2) arc(0:180:0.5cm);
   \draw (2,1) -- (2,2);
   \draw [postaction={decorate}](0,3)node[left]{$A$} -- (0,2);
   \end{scope}

   \draw (3, 3) node{$=$};

   \draw [postaction={decorate}](4, 5.5) -- (4,0) node[right]{$A$};

   \draw (2, -1) node{$(a)$};
 \end{scope}

 \begin{scope}[xshift = 8cm, decoration={
    markings,
    mark=at position 0.5 with {\arrow{>}}}]
     \draw [postaction={decorate}] (0, 0.5) -- (0,0) node[left]{$A$};
      \draw [postaction={decorate}] (2, 0.5) -- (2,0) node[left]{$B$};
    \Braid{0}{0.5}
    \draw [postaction={decorate}] (0, 3) -- (0,2.5) node[left]{$\lsup{|A|}{B}$};
     \draw [postaction={decorate}] (2, 3) -- (2,2.5);
     \InvBraid{0}{3}
     \draw [postaction={decorate}] (0, 5.5)node[left]{$A$} -- (0,5) ;
      \draw [postaction={decorate}] (2, 5.5)node[left]{$B$} -- (2,5) ;
      \draw (3,3) node{$=$};
      \draw [postaction={decorate}](4, 5.5) -- (4,0) node[right]{$A$};
      \draw [postaction={decorate}](6, 5.5) -- (6,0) node[right]{$B$};

   \draw (3, -1) node{$(b)$};
 \end{scope}

 \begin{scope}[xshift = 18cm, decoration={
    markings,
    mark=at position 0.5 with {\arrow{>}}}]
     \draw [postaction={decorate}] (0, 0.7) -- (0,0) node[left]{$A$};
      \draw [postaction={decorate}] (1.6, 0.7) -- (1.6,0) node[left]{$B$};
       \draw [postaction={decorate}] (3.2, 0.7) -- (3.2,0) node[left]{$C$};
       \SmallBraid{0}{0.7}{1.6}
       \SmallBraid{1.6}{2.3}{1.6}
       \SmallBraid{0}{3.9}{1.6}
       \draw (3.2, 2.3) -- (3.2, 0.7);
       \draw (0, 3.9) -- (0, 2.3) node[left]{$\lsup{|A|}{B}$};
       \draw (3.2, 5.5) -- (3.2, 3.9);
       \draw (2, 4.3) node{$\lsup{|A|}{C}$};
       \draw (-1.2, 5.5) node{$\lsup{|A||B|}{C}$};

       \draw (4.5,3) node{$=$};

       \begin{scope}[xshift = 6.5cm]
       \draw [postaction={decorate}] (0, 0.7) -- (0,0) node[left]{$A$};
      \draw [postaction={decorate}] (1.6, 0.7) -- (1.6,0) node[left]{$B$};
       \draw [postaction={decorate}] (3.2, 0.7) -- (3.2,0) node[left]{$C$};
       \SmallBraid{0}{2.3}{1.6}
       \SmallBraid{1.6}{0.7}{1.6}
       \SmallBraid{1.6}{3.9}{1.6}
       \draw (0, 2.3) -- (0, 0.7);
       \draw (3.2, 3.9) -- (3.2, 2.3);
       \draw (0, 5.5) -- (0, 3.9);
       \draw (2, 2.7) node{$\lsup{|B|}{C}$};
       \draw (-1.2, 5.5) node{$\lsup{|A||B|}{C}$};
       \draw (2, 5.5) node{$\lsup{|A|}{B}$};
       \end{scope}

       \draw (4.5, -1) node{$(c)$};

 \end{scope}

 \end{tikzpicture}
 \caption{Regular isotopies of a graph diagram (I)}\label{fig:regular}
\end{figure}

\begin{figure}
 \centering
\begin{tikzpicture}[scale = 0.5]
 \begin{scope}[decoration={
    markings,
    mark=at position 0.5 with {\arrow{>}}}]
    \draw [postaction={decorate}] (0,3) -- (0,0)node[left]{$A$};
    \draw [postaction={decorate}] (2,1) -- (2,0)node[right]{$B$};
    \draw [postaction={decorate}] (2,3)node[right]{$C$} -- (2,2);
     \draw (1.5,1) rectangle (2.5,2)node[pos = 0.5]{$f$};
     \Braid{0}{3}
     \draw (-1,5) node{$\lsup{|A|}{C}$};

     \draw (4, 2.5) node{$=$};
     \begin{scope}[xshift = 7cm]
      \draw [postaction={decorate}] (0,5)node[left]{$\lsup{|A|}{C}$} -- (0,4);
    \draw [postaction={decorate}] (0,3) -- (0,2)node[left]{$\lsup{|A|}{B}$};
    \draw [postaction={decorate}] (2,5) -- (2,2);
     \draw (-0.7,3) rectangle (0.7,4)node[pos = 0.5]{$\lsup{{\tiny |A|}}{f}$};
     \Braid{0}{0}
     \draw (-0.5,0) node{$A$};
     \draw (2.5,0) node{$B$};
     \end{scope}

     \draw (4, -1) node{$(a)$};
\end{scope}

\begin{scope}[xshift = 13cm, decoration={
    markings,
    mark=at position 0.5 with {\arrow{>}}}]
    \draw [postaction={decorate}] (2,3) -- (2,0)node[right]{$C$};
    \draw [postaction={decorate}] (0,1) -- (0,0)node[left]{$A$};
    \draw [postaction={decorate}] (0,3)node[left]{$B$} -- (0,2);
     \draw (-0.5,1) rectangle (0.5,2)node[pos = 0.5]{$f$};
     \Braid{0}{3}
     \draw (-1,5) node{$\lsup{|B|}{C}$};

     \draw (4, 2.5) node{$=$};
     \begin{scope}[xshift = 7cm]
      \draw [postaction={decorate}] (2,5)node[right]{$B$} -- (2,4);
    \draw [postaction={decorate}] (0,5)node[left]{$\lsup{|A|}{C}$} -- (0,2);
    \draw [postaction={decorate}] (2,3) -- (2,2);
     \draw (1.5,3) rectangle (2.5,4)node[pos = 0.5]{$f$};
     \Braid{0}{0}
     \draw (-0.5,0) node{$A$};
     \draw (2.5,0) node{$C$};
     \end{scope}

     \draw (4, -1) node{$(b)$};
\end{scope}

 \end{tikzpicture}
 \caption{Regular isotopies of a graph diagram (II). Note that in $(b)$, $|A| = |B|$.}\label{fig:regular2}
\end{figure}
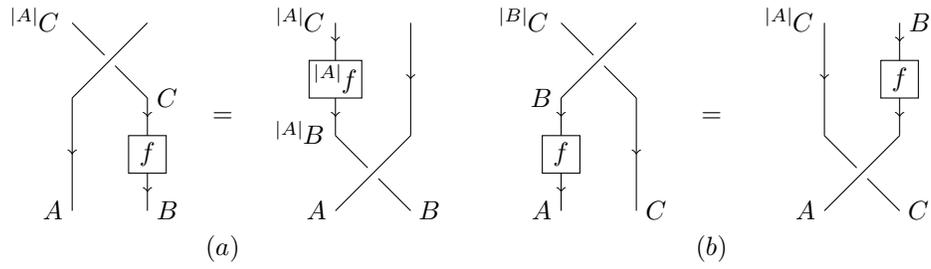

Moreover, the morphism represented by a colored diagram is invariant under regular isotopies of the diagram, some of which are drawn in Figures \ref{fig:regular} and \ref{fig:regular2}, where for simplicity we assume the coloring at each crossing is the identity map. We refer the readers to \cite{turaev1994quantum} for the treatment of more general cases.

%
%

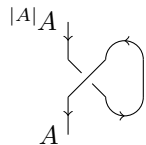
\begin{figure}
 \centering
\begin{tikzpicture}[scale = 0.5]
 \begin{scope}[decoration={
    markings,
    mark=at position 0.5 with {\arrow{>}}}]
    \draw [postaction={decorate}] (0,1)  -- (0,0) node[left]{$A$};
   \draw [postaction={decorate}] (1,1) arc(-180:0:0.5cm);
   \draw [postaction={decorate}] (2,2) arc(0:180:0.5cm);
   \SmallBraid{0}{1}{1}
   \draw (2,1) -- (2,2);
   \draw [postaction={decorate}](0,3)node[left]{$\lsup{|A|}{A}$} -- (0,2) ;
 \end{scope}
\end{tikzpicture}
\caption{The twist $\theta_A$} \label{fig:twist crossed}
\end{figure}

Similar to the case of braided spherical fusion categories, here we can also define the twist $\theta_A\colon A \longrightarrow \lsup{|A|}{A}$ for a homogeneous object $A$ by the diagram in Figure \ref{fig:twist crossed}.

The following proposition states properties of the twist parallel to those in a spherical braided category.
\begin{proposition}\cite{Kirillov_ong-equivariant}
\label{prop:twist_prop}
In a \CatNameShort{$G$}, the twist $\theta_{(\cdot)}$ satisfies$\colon$
\begin{enumerate}
 \item $\theta_{A \otimes B} = (\theta_{\lsup{gh\bar{g}}{A}} \otimes \theta_{\lsup{g}{B}})\circ c_{\lsup{g}{B},A} \circ c_{A,B}$, for $A \in \Cat_g,\ B \in \Cat_h$, and $\bar{g}:= g^{-1}$;
 \item $\theta_{A^{*}} = (\theta_{A})^*$;
 \item $\theta_{\unit} = Id$;
 \item $\lsup{h}{(\theta_{A})} = \theta_{\lsup{h}{A}}$.
\end{enumerate}
\end{proposition}

\subsection{Skeletonizing \CatNameShort{$G$}s}
\label{subsec:skeletonization}
Let $\cross{\Cat}{G}$ be a \CatNameShort{$G$}. We extract a set of discrete data to characterize $\cross{\Cat}{G}$. See \cite{maissam2014symmetry} for a similar treatment. Let $\tilde{L} = \Ltil{\cross{\Cat}{G}}$ be a complete set of representatives of simple objects, namely, $\Ltil{\cross{\Cat}{G}}$ contains a representative for each isomorphism class of simple objects. We assume $\cross{\Cat}{G}$ to have the property that $\tilde{L}$ can be chosen to be closed under the $G$-action and taking duals\footnote{Any \CatNameShort{$G$} is equivalent to one with such a property.}. For simplicity, we also assume $\cross{\Cat}{G}$ is multiplicity-free, but the method presented here can be easily adjusted to the general case. A convention within this subsection is that all the variables in a summation are implicitly assumed to be in the range of $\tilde{L}$ unless otherwise stated. {\it An undirected graph diagram means that all the segments are directed downwards.} A \CatNameShort{$G$} is described by the data $\left(\tilde{L},\ \overline{(\cdot)},\ |\cdot|,\  \lsup{(\cdot)}{(\cdot)} \right)$ together with the data $\left(N_{ab}^c,\ F^{abc}_{d;nm},\ t_a,\  U_g(a,b;c),\ \eta_a(g,h), \ R_c^{ab}\right)$ defined as follows, where $a,b,c,d,m,n \in \tilde{L}, \, g,h \in G$.

\textbf{Simple objects.} As defined above, $\tilde{L}$ is a complete set of representatives of simple objects. There is a special element, the unit object $\unit \in \tilde{L}$.

\textbf{The dual.} $\overline{(\cdot)}\colon \tilde{L} \longrightarrow \tilde{L}$ is an involution giving the dual of an object. We write the dual $\bar{a}$ and $a^*$ interchangeably. The unit object is self-dual, i.e., $\unit^* = \unit$.

\textbf{$G$-grading.}  $|\cdot|\colon \tilde{L} \longrightarrow G$ denotes the grading map, namely, if $a \in \tilde{L}$ is an object from the $g$-sector, then $|a| := g$. In particular, $|\unit|$ is the identity element of $G$.

\textbf{$G$-action.} $\lsup{(\cdot)}{(\cdot)}\colon G \times \tilde{L} \longrightarrow \tilde{L}$ denotes the $G$-action on $\tilde{L}$ which commutes with the map $\overline{(\cdot)}$ of taking the dual. The action of the group element $g$ on an object $a \in \tilde{L}$ is then given by $\lsup{g}{a}$.

\begin{figure}
 \begin{tikzpicture}[scale = 0.5]
 \begin{scope}
  \Yshape{0}{2}{a}{b}{c}
  \InverseYshape{6}{2}{a}{b}{c}

  \draw (0,-1) node{$ {\tiny B_{c}^{ab}}$};
  \draw (6,-1) node{${\tiny B_{ab}^c}$};
  \end{scope}


  \begin{scope}[xshift = 13cm]
  \Phishape{0}{0}{a}{b}{c}{c'}
  \draw (2,1) node{{\tiny $=$}};
  \draw (3,1) node{{\tiny $\delta_{c,c'}$}};
  \Vline{4}{-1}{4}{c}
  \end{scope}

  \begin{scope}[xshift = 20cm]
  \draw (0,1) node {$\sum\limits_{{\tiny c}}$};
  \Ishape{2}{0}{a}{b}{a}{b}{c}
  \draw (4,1) node {{\tiny $=$}};
  \Vline{5}{-1}{4}{a}
  \Vline{7}{-1}{4}{b}
  \end{scope}

 \end{tikzpicture}
 \caption{Graphical definition and normalizations of $B_{ab}^c$ and $B_{c}^{ab}$}\label{fig:identities}
\end{figure}
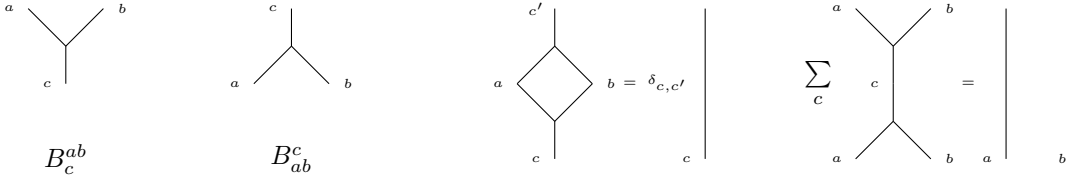

\textbf{Fusion rule.}
\begin{equation}
a \otimes b \ \simeq \ \bigoplus\limits_{c} N_{ab}^c \,c\, , \quad N_{ab}^c = 0,1.
\end{equation}
A triple $(a,b,c)$ is called {\it admissible} if $N_{ab}^{c} = 1$. In this case the morphism spaces $\Hom(c, a \otimes b)$ and $\Hom(a \otimes b, c)$ are both one-dimensional and we choose a basis element $B_{c}^{ab} \in \Hom(c, a \otimes b)$, $B_{ab}^c \in \Hom(a \otimes b, c)$ so that they satisfy the following normalization conditions. \footnote{The normalization conditions here are different from those in the physics literature where they take $B_{ab}^{c'} \circ B_{c}^{ab} = \delta_{c,c'}\sqrt{\frac{d_ad_b}{d_c}}Id_{c}$. But the data to be described below will be the same under the two normalizations.}
\begin{equation}
B_{ab}^{c'} \circ B_{c}^{ab} = \delta_{c,c'} Id_c \quad \text{and} \quad \sum\limits_{c} B_{c}^{ab}\circ B_{ab}^c = Id_{a \otimes b},
\end{equation}
where $B_{ab}^c$ and $B_{c}^{ab}$ are defined to be zero whenever $N_{ab}^c = 0$. See Figure \ref{fig:identities} for their graphical representations.

\textbf{\boldmath{$F$}-symbol (\boldmath{$6j$}-symbol).} When $\Hom(d, a \otimes b \otimes c) \neq 0$, it has two sets of bases,
    \begin{align*}
    \mathcal{B}_1 &= \{(B_{m}^{ab} \otimes I)\circ B_{d}^{mc}\ |\ m \in \tilde{L}, (a,b,m), (m,c,d) \text{ admissible}\}, \\
    \mathcal{B}_2 &= \{(I \otimes B_{n}^{bc})\circ B_{d}^{an}\ |\ n \in \tilde{L}, (b,c,n), (a,n,d) \text{ admissible}\}.
    \end{align*}
\noindent These two bases are represented graphically in Figure \ref{fig:twobases}. The matrix relating the two bases is called an $F$-matrix and the matrix elements are called $F$-symbols or $6j$-symbols. Explicitly, given a $6$-tuple $(a,b,c,d,m,n) \in \tilde{L}^6$, if $(a,b,m),\ (m,c,d),\ (b,c,n), $ and $ (a,n,d)$ are all admissible, (we say the $6$-tuple is admissible) then the $F^{abc}_{d;nm}$ and $(F^{-1})^{abc}_{d;mn}$ are defined as shown in Figures \ref{fig:FMatrix} and \ref{fig:FMatrixInv}, respectively. Otherwise, $F^{abc}_{d;nm} := (F^{-1})^{abc}_{d;mn} := 0$. By definition, if the $6$-tuple $(a,b,c,d,m,n)$ is admissible, we have,
\begin{align*}
    \sum\limits_{k} F^{abc}_{d;nk}(F^{-1})^{abc}_{d;km} &= \delta_{n,m}.
\end{align*}
     \begin{figure}
     \centering
\setlength{\unitlength}{0.030in}
\begin{picture}(105,50)(0,-10)
 \put(20,10){\line(0,-1){10}}
 \put(20,10){\line(1,1){20}}
 \put(20,10){\line(-1,1){20}}
\put(10,20){\line(1,1){10}}

 \put(2,30){$a$}
 \put(22,30){$b$}
 \put(42,30){$c$}
 \put(16,16){$m$}
 \put(22,2){$d$}

 \put(20,-10){$\mathcal{B}_1$}

 \put(80,10){\line(0,-1){10}}
 \put(80,10){\line(1,1){20}}
 \put(80,10){\line(-1,1){20}}
\put(90,20){\line(-1,1){10}}

 \put(62,30){$a$}
 \put(82,30){$b$}
 \put(102,30){$c$}
 \put(82,15){$n$}
 \put(82,2){$d$}

 \put(80,-10){$\mathcal{B}_2$}
\end{picture}
\caption{Two bases $\mathcal{B}_1$ and $\mathcal{B}_2$ }\label{fig:twobases}
\end{figure}
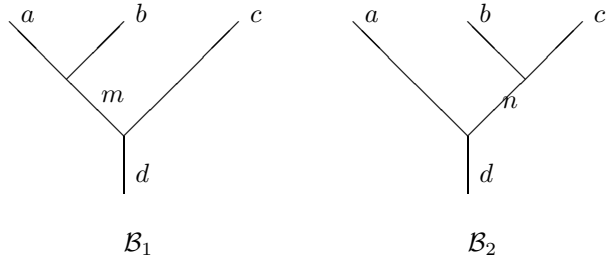

\begin{figure}
     \centering
\setlength{\unitlength}{0.030in}
\begin{picture}(105,40)(0,0)
 \put(20,10){\line(0,-1){10}}
 \put(20,10){\line(1,1){20}}
 \put(20,10){\line(-1,1){20}}
\put(10,20){\line(1,1){10}}

 \put(2,30){$a$}
 \put(22,30){$b$}
 \put(42,30){$c$}
 \put(16,16){$m$}
 \put(22,2){$d$}

  \put(45,15){=}
 \put(50,15){$\sum\limits_{n} F_{d;nm}^{abc}$}

 \put(80,10){\line(0,-1){10}}
 \put(80,10){\line(1,1){20}}
 \put(80,10){\line(-1,1){20}}
\put(90,20){\line(-1,1){10}}

 \put(62,30){$a$}
 \put(82,30){$b$}
 \put(102,30){$c$}
 \put(82,15){$n$}
 \put(82,2){$d$}
\end{picture}

\caption{Definition of the $F$-symbol}\label{fig:FMatrix}
\end{figure}

\begin{figure}
     \centering
\setlength{\unitlength}{0.030in}
\begin{picture}(105,40)(0,0)

 \put(20,10){\line(0,-1){10}}
 \put(20,10){\line(1,1){20}}
 \put(20,10){\line(-1,1){20}}
 \put(30,20){\line(-1,1){10}}

 \put(2,30){$a$}
 \put(22,30){$b$}
 \put(42,30){$c$}
 \put(22,15){$n$}
 \put(22,2){$d$}

  \put(45,15){=}
 \put(50,15){$\sum\limits_{n} (F^{-1})_{d;mn}^{abc}$}

 \put(90,10){\line(0,-1){10}}
 \put(90,10){\line(1,1){20}}
 \put(90,10){\line(-1,1){20}}
\put(80,20){\line(1,1){10}}

 \put(72,30){$a$}
 \put(92,30){$b$}
 \put(112,30){$c$}
 \put(86,16){$m$}
 \put(92,2){$d$}
\end{picture}
\caption{Definition of the $F^{-1}$-symbol}\label{fig:FMatrixInv}
\end{figure}

\textbf{\boldmath{$t$}-symbol.} For $a \in \tilde{L}$, $t_a$ is a nonzero scalar such that $t_a Id\colon a \longrightarrow \bar{\bar{a}}$ corresponds to the pivotal structure of $\cross{\Cat}{G}$ (see Section $4.2$ of \cite{Wang2010topological}).

\textbf{\boldmath{$U$}-symbol.} For any $g \in G, \, a,b,c \in \tilde{L}$ such that $(a,b,c)$ is admissible, $\lsup{g}{(B_{c}^{ab})}$ is a scalar multiple of $B_{\lsup{g}{c}}^{\lsup{g}{a}\lsup{g}{b}}$. Denote this scalar by $U_g(a,b;c)$. Graphically, this is represented by Figure \ref{fig:u-symbol}.
     \begin{figure}
     \centering
     \begin{tikzpicture}[scale = 0.5]
     \YshapeNormalFont{0}{0}{a}{b}{c}
     \draw plot[smooth] coordinates{(-1.5,1.5)(-2,0)(-1.5,-1.5)};
     \draw plot[smooth] coordinates{(1.5,1.5)(2,0)(1.5,-1.5)};
     \draw (-2,1) node{{\tiny $g$}};

     \draw (4,0) node{$=$};
     \draw (7,0) node{$U_{g}(a,b;c)$};

     \YshapeNormalFont{10}{0}{\lsup{g}{a}}{\lsup{g}{b}}{\lsup{g}{c}}
     \end{tikzpicture}
     \caption{Definition of $U$-symbol}\label{fig:u-symbol}
     \end{figure}

\textbf{\boldmath{$\eta$}-symbol.} For $g,h \in G,\, a \in \tilde{L}$, $\eta_a(g,h)$ is the scalar such that $\eta_a(g,h)Id$ represents the isomorphism $\eta(g,h;a)\colon \lsup{gh}{a} \longrightarrow \lsup{g}{(\lsup{h}{a})}$.

 \textbf{\boldmath{$R$}-symbol.} If $(a,b,c)$ is admissible, then $R_{c}^{ab}$ is defined to be the scalar such that $c_{a,b} \circ B_{c}^{ab} = R_{c}^{ab} B_{c}^{\lsup{|a|}{b}a}$, or equivalently $c_{a,b}^{-1} \circ B_{c}^{\lsup{|a|}{b}a} = (R_{c}^{ab})^{-1}B_{c}^{ab}$.  See Figure \ref{fig:RMatrix} for a graphical definition of the $R$-symbol.

\begin{figure}
\centering
\begin{tikzpicture}[scale = 0.5]
\begin{scope}
\YshapeNormalFont{0}{0}{a}{b}{c}
\Braid{-1}{1}
\draw (-1.5, 3) node{$\lsup{|a|}{b}$};
\draw (1.5,3) node{$a$};
\end{scope}
\begin{scope}[xshift = 3.5cm]
\draw (0,1) node{$=$  $R_{c}^{ab}$};
\end{scope}
\begin{scope}[xshift = 6cm]
\draw (0,-1) -- (0,1) -- (-1, 3);
\draw (0,1) -- (1,3);
\draw (-0.5, -1) node{$c$};
\draw (-1.5, 3) node{$\lsup{|a|}{b}$};
\draw (1.5, 3) node{$a$};
\end{scope}
\begin{scope}[xshift = 12cm]
\begin{scope}
\YshapeNormalFont{0}{0}{\lsup{|a|}{b}}{a}{c}
\InvBraid{-1}{1}
\draw (-1.5, 3) node{$a$};
\draw (1.5,3) node{$b$};
\end{scope}
\begin{scope}[xshift = 3.5cm]
\draw (0,1) node{$=$  $(R_{c}^{ab})^{-1}$};
\end{scope}
\begin{scope}[xshift = 6cm]
\draw (0,-1) -- (0,1) -- (-1, 3);
\draw (0,1) -- (1,3);
\draw (-0.5, -1) node{$c$};
\draw (-1.5, 3) node{$a$};
\draw (1.5, 3) node{$b$};
\end{scope}
\end{scope}
\end{tikzpicture}
\caption{Definition of the $R$-symbol $R^{ba}_c$}\label{fig:RMatrix}
\end{figure}

Translating the axioms of a \CatNameShort{$G$} to the data $(N_{ab}^c,\ F^{abc}_{d;nm},\ t_a,\ U_g(a,b;c),\ \eta_a(g,h),\ R_c^{ab})$, we arrive at Equations \ref{equ:fusion} - \ref{equ:RU}, where $a,b,c,d,m,n,l,p,q \in \Ltil{\cross{\Cat}{G}}, \, g,h,k \in G$.
\begin{align}
\label{equ:fusion}
N_{ab}^{c} = N_{\bar{a}c}^{b} = N_{c\bar{b}}^{a} = N_{\lsup{|a|}{b}a}^c =N_{\lsup{g}{a}\lsup{g}{b}}^{\lsup{g}{c}}, \nonumber\\
N_{ab}^{\unit} = \delta_{a, \bar{b}},  \qquad N_{ab}^c     = 0 \text{ if } |a||b| \neq |c| ,\nonumber\\
\sum\limits_{m} N_{ab}^m N_{mc}^d = \sum\limits_{n} N_{bc}^n N_{an}^d .             
\end{align}

\begin{align}
\label{equ:pentagon}
F^{a \unit b}_{c;ba} = 1 \text{ whenever } (a,b,c) \text{ is admissible}, \nonumber \\
F^{a \bar{a} a}_{a; \unit\unit} \neq 0, \nonumber \\
F^{mcd}_{f;qn}F^{abq}_{f;pm} = \sum\limits_{l} F^{abc}_{n;lm}F^{ald}_{f;pn}F^{bcd}_{p;ql}.
\end{align}

\begin{align}
\label{equ:pivotal}
t_{\unit} = 1, \qquad t_{\bar{a}} = t_a^{-1}, \quad t_{\lsup{g}{a}} = t_a, \nonumber \\
t_a^{-1}t_b^{-1}t_c = F^{ab\bar{c}}_{\unit;\bar{a}c}F^{b\bar{c}a}_{\unit;\bar{b}\bar{a}}F^{\bar{c}ab}_{\unit;c\bar{b}} \text{ whenever } (a,b,c) \text{ is admissible.}
\end{align}

\begin{align}
\label{equ:FU}
F^{\lsup{g}{a}\lsup{g}{b}\lsup{g}{c}}_{\lsup{g}{d};\lsup{g}{n}\lsup{g}{m}} &= F^{abc}_{d;nm} \frac{U_g(b,c;n)U_g(a,n;d)}{U_g(a,b;m)U_g(m,c;d)}
\end{align}

\begin{align}
\label{equ:eta}
\eta_a(gh,k)\eta_{\lsup{k}{a}}(g,h) &= \eta_a(g,hk) \eta_a(h,k)
\end{align}

\begin{align}
\label{equ:Ueta}
\frac{\eta_a(g,h)\eta_b(g,h)}{\eta_c(g,h)} &= \frac{U_g(\lsup{h}{a},\lsup{h}{b};\lsup{h}{c})U_h(a,b;c)}{U_{gh}(a,b;c)}
\end{align}

\begin{align}
\label{equ:hexagon1}
R_{m}^{ab}F^{\lsup{|a|}{b}ac}_{d;nm}R_{n}^{ac} &= \sum\limits_{l} F^{abc}_{d;lm}U_{|a|}(b,c;l)R_{d}^{al}F^{\lsup{|a|}{b}\lsup{|a|}{c}a}_{d;n\lsup{|a|}{l}}
\end{align}

\begin{align}
\label{equ:hexagon2}
(R^{b\lsup{|c|}{a}}_{m})^{-1} F^{b \lsup{|c|}{a} c}_{d;n,m}(R^{ca}_{n})^{-1} &= \sum\limits_{l} \eta_a(|b|,|c|)^{-1} F^{\lsup{|b||c|}{a}bc}_{d;lm}(R^{la}_d)^{-1}F^{bca}_{d;nl}
\end{align}

\begin{align}
\label{equ:RU}
R^{\lsup{g}{a}\lsup{g}{b}}_{\lsup{g}{c}} &= R^{ab}_{c} \frac{U_g(\lsup{|a|}{b},a;c )\eta_b(g |a| g^{-1}, g)}{U_g(a,b;c)\eta_b(g, |a|)}
\end{align}

There is a gauge freedom in relabeling the elements of $\tilde{L}$ as well as modifying the choice of the basis elements $B_{ab}^c$ and $B_{c}^{ab}$. For instance, one can also choose $\tilde{B}_{ab}^{c}:= \lambda(a,b;c) B_{ab}^c$ and $\tilde{B}_{c}^{ab}:= \lambda(a,b;c)^{-1}B_{c}^{ab}$ for arbitrary nonzero scalars $\lambda(a,b;c)$. Then the discrete data set, namely, the $F$-symbols, $R$-symbols, etc, will also be changed correspondingly. We call two such data sets gauge equivalent. 

Conversely, given a discrete data set $\left(\tilde{L},\ \overline{(\cdot)},\  |\cdot|,\  \lsup{(\cdot)}{(\cdot)} \right)$ and $\left(N_{ab}^c,\  F^{abc}_{d;nm},\  t_a, \  U_g(a,b;c),\  \eta_a(g,h),\  R_c^{ab}\right)$ satisfying  Equations \ref{equ:fusion} - \ref{equ:RU}, there is a canonical way to construct a \CatNameShort{$G$} with the property that each isomorphism class contains exactly one object. A category with such a property is called {\it skeletal}. The construction is similar to that in \cite{davidovich2013arithmetic, yamagami2002polygonal} except here we also need to define the $G$-action. Since this is straightforward and it will not be used for the rest of the paper, we omit the details.

In Section \ref{subsec:partition3}, we will give a formula for the invariant of a closed $4$-manifold in terms of the data mentioned above.

\section{Partition Function}
\label{sec:partition}
In this section we introduce the key construction of the $4$-manifold invariant from a \CatNameShort{$G$}. The invariant, also called the partition function, is a state-sum model and a $4d$ analogue of the Turaev-Viro-Barrett-Westbury invariant. Although we will only describe the invariant for oriented closed $4$-manifolds, there should be no conceptual difficulty in extending the construction to produce a $(3+1)$-$\TQFT$. See \cite{williamson2016hamiltonian} for a Hamiltonian realization of this $\TQFT$. By $4$-manifolds we always mean closed oriented smooth $4$-manifolds unless otherwise stated. In Sections \ref{subsec:partition1}, \ref{subsec:partition2}, and \ref{subsec:partition3}, three equivalent definitions of the partition function based on ordered triangulations will be given. To show that this indeed defines an invariant of $4$-manifolds, one needs to prove that the partition function is independent of the choice of ordered triangulations. Considering the proof is quite technical and lengthy, we defer it to Section \ref{sec:proofmain}.

\subsection{Definition of Partition Function I}
\label{subsec:partition1}
An ordered triangulation $\T$ of a $4$-manifold $M$ is a triangulation of $M$ with an ordering of its vertices by $0,1,2,\cdots.$ For $k= 0,1,2,3,4$, let $\T^k$ be the set of $k$-simplices. The restriction of the ordering on each $k$-simplex $\sigma \in \T^k$ induces a relative ordering of the vertices of $\sigma$. Under this relative ordering, we write $\sigma$ as $(0\ 1\ \cdots\ k-1 )$ and any $n$-face of $\sigma$ as $(i_1\ i_2\ \cdots \ i_n)$, $0 \leq i_1 < \cdots < i_n \leq k-1$.  The definition of the partition function only depends on the relative ordering as shall be seen below. For each $\sigma \in \T^4$, define the sign, $\epsilon(\sigma)$, of $\sigma$ to be \lq $+$' if the orientation on $\sigma$ induced from $M$ coincides with the one determined by the relative ordering of its vertices; $\epsilon(\sigma)$ is defined to be \lq $-$' otherwise. Let $\cross{\Cat}{G} = \bigoplus\limits_{g\in G}\Cat_g$ be a \CatNameShort{$G$}, and let $\L{\cross{\Cat}{G}}$ be the set of isomorphism classes of simple objects.

\begin{definition}
\label{def:coloring}
 A $\cross{\Cat}{G}$-coloring of $(M, \T)$ is a pair of maps $F = (g,f)$, $g\colon \T^1 \rightarrow G$, $f\colon \T^2 \rightarrow \L{\cross{\Cat}{G}}$, such that for any simplex $\beta = (012) \in \T^2$ with the induced ordering on its vertices,
$$f(012) \in \Cat_{g(02)^{-1}g(01)g(12)}.$$
\end{definition}

Given a $\cross{\Cat}{G}$-coloring $F = (g,f)$, for each $\beta = (012) \in \T^2$, we arbitrarily choose a representative in the class $f(012)$ and denote the representative by $f_{012}$ or just $012$ when no confusion arises. For an edge $(ij)$, we also denote $g(ij)$ by $g_{ij}$ or $ij$ and denote $g(ij)^{-1}$ by $\overline{g}_{ij}$ or $\overline{ij}$.



Assume now a representative for the color of each $2$-simplex has been chosen. Let $\tau = (0123) \in \T^3$ be any $3$-simplex with the induced ordering on its vertices, and consider the boundary map,
$$\partial(0123) = (123) - (023) + (013) - (012).$$
We assign to $\tau = (0123)$ two vector spaces,
$$V^{+}_{F}(0123) := \Hom(f_{023} \otimes \lsuprsub{\overline{g}_{23}}{f}{012}, f_{013} \otimes f_{123}),$$
$$V^{-}_{F}(0123) := \Hom(f_{013} \otimes f_{123}, f_{023} \otimes \lsuprsub{\overline{g}_{23}}{f}{012}).$$
Note that by the definition of a coloring, the objects $f_{023} \otimes \lsuprsub{\overline{g}_{23}}{f}{012}$ and $f_{013} \otimes f_{123}$ are both in the sector $\Cat_{\overline{g}_{03}g_{01}g_{12}g_{23}}$. This is a necessary requirement since otherwise the spaces $V^{+}_{F}(0123)$ and $V^{-}_{F}(0123)$ would equal $0$.

For any two objects $X,Y$ in a spherical fusion category, a pairing on $\Hom(X,Y) \times \Hom(Y,X)$ is defined as,
\begin{align}
\label{equ:pairingdef}
\begin{array}{rccc}
\langle\;,\;\rangle\ \colon & \Hom(X,Y) \times \Hom(Y,X) &\longrightarrow  & \C \nonumber \\
                            &(\phi, \psi)                &\mapsto &           \Tr(\phi\psi)
\end{array}
\end{align}
By \cite{Etingof_onfusion}, for each simple object $a$, $\Tr(id_a)$, namely, the quantum dimension of $a$, is nonzero. It is straightforward to show that the pairing $\langle\;,\;\rangle$ is nondegenerate, and thus induces a natural isomorphism between $V^{+}_{F}(0123)$ and $V^{-}_{F}(0123)^{*}$, as well as $V^{-}_{F}(0123)$ and $V^{+}_{F}(0123)^{*}$.

\begin{figure}
\centerline{\includegraphics{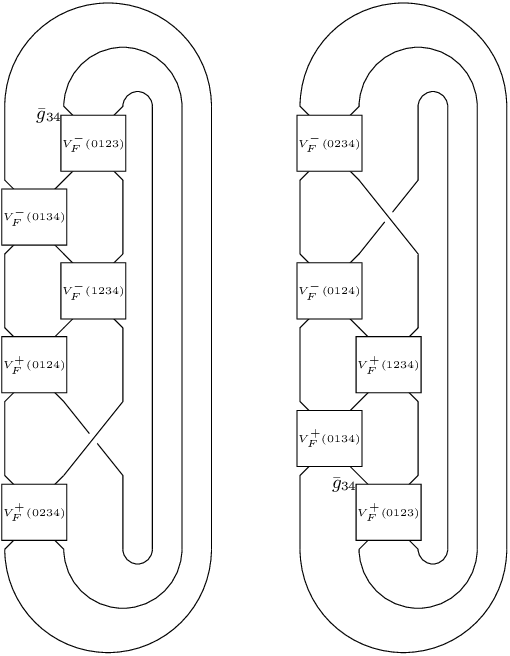}}
\caption{(Left): $\tilde{Z}^{+}_{F}(01234)$; (Right): $\tilde{Z}^{-}_{F}(01234)$}
\label{fig:25jbox}
\end{figure}

Let $\sigma = (01234) \in \T^4$ be any $4$-simplex with the induced ordering on its vertices, consider the boundary map,
$$\partial(01234) = (1234) - (0234) + (0134) - (0124) + (0123).$$
We define $Z^{\epsilon(\sigma)}_{F}(\sigma)$ as follows. If $\epsilon(\sigma) = +$, we first define a linear  functional
\begin{align*}
\begin{array}{rrcl}
\tilde{Z}^{+}_{F}(01234)\ \colon &V^{+}_{F}(0234) \otimes V^{+}_{F}(0124) \otimes V^{-}_{F}(1234) \otimes V^{-}_{F}(0134) \otimes V^{-}_{F}(0123)& \longrightarrow & \C
\end{array}
\end{align*}
by the graph diagram shown in Figure \ref{fig:25jbox} (Left), the meaning of which requires some explanations. This is a graph diagram as defined in Section \ref{subsec:stratification} where the segments are implicitly colored and directed, but the rectangles are not. The three long vertical segments on the right part of the diagram (which correspond to taking the trace) are directed upwards, and all remaining segments are directed downwards. The colors of the segment are assigned in such a way that any morphism from the morphism space indicated in a rectangle can be used to color that rectangle. For instance, for the rectangle containing $V_F^{+}(0234) = \Hom(034 \otimes \lsup{\overline{34}}{023}, 024 \otimes 234)$, the four segments incident to it, (lower left, lower right, upper left, upper right), are colored by ($034, \ \lsup{\overline{34}}{023}, \ 024, \   234$), respectively. Note that here $(034)$ really denotes $f_{034}$. For the top rectangle containing $V_F^{-}(0123)$, the $\overline{g}_{34}$ on the upper left corner means the action of $\overline{g}_{34}$ on $V_F^{-}(0123)$. Since $V_F^{-}(0123) = \Hom(013 \otimes 123, 023 \otimes \lsup{\overline{23}}{012})$, the four segments incident to the top rectangle, (lower left, lower right, upper left, upper right) are colored by $(\lsup{\overline{34}}{013}, \ \lsup{\overline{34}}{123}, \lsup{\overline{34}}{023}, \lsup{\overline{34}\cdot \overline{23}}{013})$, respectively. In the case of a strict \CatNameShort{$G$}, the only crossing in the diagram is colored by the identity map. In the nonstrict case, it is colored by some natural isomorphism involving the $\eta_a(g,h)\,'$s defined in Section \ref{subsec:skeletonization}. One can check that the above rules for the coloring are compatible. Given an element $\phi_0 \otimes \phi_1 \otimes\phi_2 \otimes\phi_3 \otimes\phi_4$ in the domain of $\tilde{Z}^{+}_{F}(01234)$, we now color each rectangle by  some $\phi_i$ which is from the space the rectangle contains. The resulting colored graph diagram can be interpreted as a morphism in $\Hom(\unit, \unit)\simeq \C$ according to Section \ref{subsec:stratification}. Since any such morphism is a scalar times $Id_{\unit}$, we define $\tilde{Z}^{+}_{F}(01234)(\phi_0 \otimes \phi_1 \otimes\phi_2 \otimes\phi_3 \otimes\phi_4)$ to be that scalar, or less rigorously to be (the evaluation of) the colored graph diagram.

\begin{remark}
We give some intuition on how the graph diagrams in Figure \ref{fig:25jbox} are derived. Roughly, if one views each rectangle as a vertex, then the two graph diagrams are both the $1$-skeleton of the dual triangulation of $\Sphere^3 = \partial(01234)$. Of course, one has to be careful in projecting the $1$-skeleton to the plane, so that the projected diagram is consistent with the colorings and the interpretation of the colored diagrams as the trace of certain morphisms in the category makes sense. See also \cite{crane1997state, mackaay1999spherical} for a similar diagram. 
\end{remark}

Using the nondegenerate pairing $\langle\;,\;\rangle$, we can define a linear map $Z^{+}_{F}(01234)$,
\begin{align*}
\begin{array}{rrcl}
Z^{+}_{F}(01234)\ \colon & V^{+}_{F}(0234) \otimes V^{+}_{F}(0124)& \longrightarrow & V^{+}_{F}(1234) \otimes V^{+}_{F}(0134) \otimes V^{+}_{F}(0123),
\end{array}
\end{align*}
such that
$$\langle Z^{+}_{F}(01234)(\phi_0 \otimes \phi_1), \phi_2 \otimes \phi_3 \otimes \phi_4 \rangle \ = \ \tilde{Z}^{+}_{F}(01234)(\phi_0 \otimes \phi_1 \otimes\phi_2 \otimes\phi_3 \otimes\phi_4) $$.


Similarly, if $\epsilon(\sigma) = -$, consider instead the functional
\begin{align*}
\begin{array}{rrcl}
\tilde{Z}^{-}_{F}(01234)\ \colon V^{-}_{F}(0234) \otimes V^{-}_{F}(0124) \otimes V^{+}_{F}(1234) \otimes V^{+}_{F}(0134) \otimes V^{+}_{F}(0123) &\longrightarrow & \C
\end{array}
\end{align*}
defined by the graph diagram shown in Figure \ref{fig:25jbox} (Right), which is obtained by reflecting the one in Figure \ref{fig:25jbox} (Left) along a horizon line. By the same way of interpreting the diagram as above, we get a linear map,
\begin{align*}
\begin{array}{rrcl}
Z^{-}_{F}(01234)\ \colon  &V^{+}_{F}(1234) \otimes V^{+}_{F}(0134) \otimes V^{+}_{F}(0123) &\longrightarrow & V^{+}_{F}(0234) \otimes V^{+}_{F}(0124),
\end{array}
\end{align*}
such that
$$\langle \tilde{\phi}_0 \otimes \tilde{\phi}_1, Z^{-}_{F}(01234)(\tilde{\phi}_2 \otimes \tilde{\phi}_3 \otimes \tilde{\phi}_4)  \rangle \ = \ \tilde{Z}^{-}_{F}(01234)(\tilde{\phi}_0 \otimes \tilde{\phi}_1 \otimes\tilde{\phi}_2 \otimes\tilde{\phi}_3 \otimes\tilde{\phi}_4) .$$

From the boundary equations,
\begin{align*}
\partial(+01234) &= (1234) - (0234) + (0134) - (0124) + (0123), \\
\partial(-01234) &= -(1234) + (0234) - (0134) + (0124) - (0123),
\end{align*}
we see that the domain of $Z^{\pm}_{F}(01234)$ always corresponds to the negative boundaries of $\pm(01234)$ and the codomain to positive boundaries of $\pm(01234)$. Each $3$-simplex $\tau$ is the intersection of exactly two $4$-simplices $\sigma_1$ and $\sigma_2$. If $\tau$ is a negative boundary in $\epsilon(\sigma_1)\sigma_1$, it must be a positive boundary in $\epsilon(\sigma_2)\sigma_2$, and vice versa. Thus, $V^{+}_{F}(\tau)$ appears exactly once in the domain for some $4$-simplex and codomain for some other $4$-simplex. Then we have
\begin{align*}
\begin{array}{rrcl}
\bigotimes\limits_{\sigma \in \T^4}Z_{F}^{\epsilon(\sigma)}(\sigma)\ \colon \bigotimes\limits_{\tau \in \T^3} V^{+}_{F}(\tau) \longrightarrow \bigotimes\limits_{\tau \in \T^3} V^{+}_{F}(\tau).
\end{array}
\end{align*}

\begin{definition}
\label{def:partition1}
Given a \CatNameShort{$G$} $\cross{\Cat}{G}$ and an ordered triangulation $\T$ of a $4$-manifold $M$, the partition function $Z_{\cross{\Cat}{G}}(M;\T)$ of the pair $(M,\T)$ is defined by,
\begin{align}
Z_{\cross{\Cat}{G}}(M;\T) &= \sum\limits_{F = (g,f)} \frac{(D^2/|G|)^{|\T^0|} \left(\prod\limits_{\beta \in \T^2}d_{f(\beta)}\right) \Tr\left(\bigotimes\limits_{\sigma \in \T^4}Z^{\epsilon(\sigma)}_{F}(\sigma)\right)}{(D^2)^{|\T^1|}},
\end{align}
where $F$ runs through all $\cross{\Cat}{G}$-colorings of $\T$, and $D^2:= \sum\limits_{a \in \L{\cross{\Cat}{G}}}d_a^2$ with $d_a$ the quantum dimension of the object $a$.
\end{definition}

In the above definition, $D^2$ is called the total dimension square. For $g \in G$, we define 
\begin{align*}
D^2_g:= \sum\limits_{a \in \L{\cross{\Cat}{G}} \cap \Cat_g} d_a^2.
\end{align*} 
Let $\text{Gr}(\cross{\Cat}{G})$ be the subset of $G$ that contains all elements $g$ such that $\Cat_g$ is not zero. Then $\text{Gr}(\cross{\Cat}{G})$ is a normal subgroup of $G$. By \cite{Turaev2010Homotopy2}, we have $D^2_{g} = D^2_{e}$ for any $g \in \text{Gr}(\cross{\Cat}{G})$. Of course, $D^2_g = 0$ for $g \notin \text{Gr}(\cross{\Cat}{G})$. Hence, $D^2 = D^2_e\, |\text{Gr}(\cross{\Cat}{G})|$.

\begin{theorem}[Main Theorem]
\label{thm:main}
Let $\cross{\Cat}{G}, M, \T$ be as above, then the partition function $Z_{\cross{\Cat}{G}}(M;\T)$ is an invariant of smooth closed oriented $4$-manifolds.
\end{theorem}

To prove the Main Theorem, one needs to show that $Z_{\cross{\Cat}{G}}(M;\T)$ is independent of,
\begin{enumerate}
 \item the choice of a representative for the color of each $2$-simplex,
 \item the ordering of vertices of a triangulation,
 \item the choice of a triangulation.
\end{enumerate}

These will be proved in Section \ref{sec:invariance representative}, Section \ref{sec:invariance ordering}, and Section \ref{sec:invariance pachner}, respectively.

If $\cross{\Cat}{G}$ and $\cross{\D}{G}$ are two equivalent \CatNameShort{$G$}s, let $F$ be an equivalence between such two categories. Note that $F$ maps each $\Cat_g$ to $\D_g$. It is not hard to see that the action of $F$ on morphism spaces preserve the pairing in Equation \ref{equ:pairingdef} and the functionals $\tilde{Z}^{\pm}_{F}(01234)$, hence we have $Z_{\cross{\Cat}{G}}(M;\T) = Z_{\cross{\D}{G}}(M;\T)$. That is, $Z_{\cross{\Cat}{G}}(M;\T)$ only depends on the equivalence class of $\cross{\Cat}{G}$.

Below we derive two other equivalent formulas for the partition function, which will be used in the proof of invariance. They each have different flavors, and depending on the situation it is more convenient to use one over another.

\subsection{Definition of Partition Function II}
\label{subsec:partition2}
The nondegenerate pairing $\lrangle$ defined in Equation \ref{equ:pairingdef} for two objects $X, Y$ in a spherical fusion category induces canonical isomorphisms $\Hom(X,Y) \simeq \Hom(Y,X)^{*}$ and $\Hom(Y,X) \simeq \Hom(X,Y)^{*}$. Let $\lrangle_{X,Y}^{*}$ be the dual map of $\lrangle_{X,Y}$,
\begin{equation*}
\begin{array}{ccccc}
\lrangle_{X,Y}^{*}\ \colon &\C& \longrightarrow &\Hom(X,Y)^{*} \otimes \Hom(Y,X)^{*} \simeq & \Hom(Y,X) \otimes \Hom(X,Y),
\end{array}
\end{equation*}
and let $\phi_{X,Y} = \lrangle_{X,Y}^{*}(1)$. If $\{v_i\}_{i \in I}$ and $\{w_j\}_{j \in I}$ are a basis of $\Hom(X,Y)$ and $\Hom(Y,X)$, respectively, such that $\langle v_i, w_j \rangle = c_i\delta_{i,j}$ for some nonzero numbers $c_i$, then $\phi_{X,Y} = \sum\limits_{i \in I} c_i^{-1} w_i \otimes v_i$.

Now assume a coloring $F = (g,f)$ and a representative for the color of each $2$-simplex have been chosen. For each $3$-simplex $\tau = (0123) \in \T^3$, there is thus the element
$$\phi_{F;\tau}:= \phi_{f_{023} \otimes \lsuprsub{\overline{g}_{23}}{f}{012}, f_{013} \otimes f_{123}} \ \in \ V^{-}_{F}(0123) \otimes V^{+}_{F}(0123) .$$
Let $V_F := \bigotimes\limits_{\tau \in \T^3} V^{-}_{F}(\tau) \otimes V^{+}_{F}(\tau)$, $\phi_F:= \bigotimes\limits_{\tau \in \T^3} \phi_{F;\tau}$, then $\phi_F \in V_F$. From the definition of $\tilde{Z}^{\epsilon(\sigma)}_{F}(\sigma)$, we see that each $V_F^{\pm}(\tau)$ appears exactly once as a tensor component of the domain for some $\sigma \in \T^4$. Therefore, the map $\bigotimes\limits_{\sigma \in \T^4} \tilde{Z}^{\epsilon(\sigma)}_{F}(\sigma)$ is a functional on $V_F$.

\begin{proposition}[Definition of Partition Function II]
\label{prop:partition2}
Let $\phi_F = \bigotimes\limits_{\tau \in \T^3} \phi_{F;\tau} \in V_F$ be as above, then $ Z_{\cross{\Cat}{G}}(M;\T)$ is given by the following formula,
\begin{align}
\label{equ:partition2}
Z_{\cross{\Cat}{G}}(M;\T) &= \sum\limits_{F = (g,f)} \frac{(D^2/|G|)^{|\T^0|} \left(\prod\limits_{\beta \in \T^2}d_{f(\beta)}\right) }{(D^2)^{|\T^1|}} \left(\bigotimes\limits_{\sigma \in \T^4}\tilde{Z}^{\epsilon(\sigma)}_{F}(\sigma)\right)(\phi_{F}).
\end{align}
\begin{proof}
It suffices to prove, for a fixed coloring $F = (g,f)$ and a chosen representative for the color of each $2$-simplex, we have,
\begin{align*}
\Tr\left(\bigotimes\limits_{\sigma \in \T^4}Z^{\epsilon(\sigma)}_{F}(\sigma)\right) &= \left(\bigotimes\limits_{\sigma \in \T^4}\tilde{Z}^{\epsilon(\sigma)}_{F}(\sigma)\right)(\phi_{F}).
\end{align*}
For any $\tau \in \T^3$, choose a basis $\{v_i^{+}(\tau)\colon i \in I_{\tau}\}$ of $ V_{F}^{+}(\tau)$ and a basis $\{v_j^{-}(\tau)\colon j \in I_{\tau}\}$ of $ V_{F}^{-}(\tau)$, such that $\langle v_i^{+}(\tau), v_j^{-}(\tau) \rangle = \delta_{i,j}$. Then
\begin{align*}
\phi_{F, \tau} = \sum\limits_{i \in I_{\tau}} v_i^{-}(\tau) \otimes v_i^{+}(\tau)
\end{align*}

Thus,
\begin{align*}
\Tr\left(\bigotimes\limits_{\sigma \in \T^4}Z^{\epsilon(\sigma)}_{F}(\sigma)\right) &= \sum\limits_{\tau \in \T^3, i_{\tau} \in I_{\tau}} \Bigg\langle \bigotimes\limits_{\tau \in \T^3}v_{i_{\tau}}^{-}(\tau) , \left(\bigotimes\limits_{\sigma \in \T^4}Z^{\epsilon(\sigma)}_{F}(\sigma)\right)\left(\bigotimes\limits_{\tau \in \T^3}v_{i_{\tau}}^{+}(\tau)\right)\Bigg\rangle \\
 &= \sum\limits_{\tau \in \T^3, i_{\tau} \in I_{\tau}} \left(\bigotimes\limits_{\sigma \in \T^4}\tilde{Z}^{\epsilon(\sigma)}_{F}(\sigma)\right) \left(\bigotimes\limits_{\tau \in \T^3}v_{i_{\tau}}^{+}(\tau) \bigotimes\limits_{\tau \in \T^3}v_{i_{\tau}}^{-}(\tau)\right) \\
 &= \left(\bigotimes\limits_{\sigma \in \T^4}\tilde{Z}^{\epsilon(\sigma)}_{F}(\sigma)\right)(\phi_{F}).
\end{align*}

\end{proof}
\end{proposition}

\subsection{Definition of Partition Function III}
\label{subsec:partition3}

We give a third formulation of the partition function as a state sum model, which is convenient in terms of calculations. For simplicity, let us assume the category $\cross{\Cat}{G}$ is multiplicity-free. The more general case can be treated in a similar way except that the formula involved will be more complicated.

Let $\tilde{L}(\cross{\Cat}{G})$ be an arbitrary complete set of representatives, one for each isomorphism class of simple objects. Recall from Section \ref{subsec:skeletonization} that, for each triple $(a,b,c)$ of simple objects such that $N_{ab}^c = 1$, namely, $(a,b,c)$ being admissible, we can choose a basis element $B_{c}^{ab} \in \Hom(c, a \otimes b)$ and $B_{ab}^c \in \Hom(a\otimes b, c)$ such that for any $c,c' \in \tilde{L}(\cross{\Cat}{G})$,
\begin{equation*}
B_{ab}^{c'} \circ B_{c}^{ab} = \delta_{c,c'} Id_c \quad \text{and} \quad \sum\limits_{c} B_{c}^{ab}\circ B_{ab}^c = Id_{a \otimes b}.
\end{equation*}
The graphical representations of $B_{ab}^c, B_{c}^{ab}$ and their relations are illustrated in Figure \ref{fig:identities}.

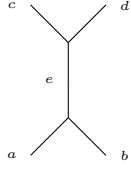
\begin{figure}
\centering
 \begin{tikzpicture}[scale = 0.5]
\Ishape{0}{0}{a}{b}{c}{d}{e}
 \end{tikzpicture}
 \caption{Graphical representation of $B_{ab,cd}^e$}\label{fig:B_abcde}
\end{figure}

For simple objects $(a,b,c,d)$, let $B_{ab,cd}^e = B_{e}^{cd} \circ B_{ab}^e$, which is represented by the graph diagram in Figure \ref{fig:B_abcde}. Then $$\{B_{ab,cd}^e\colon e \in \tilde{L}(\cross{\Cat}{G}), (a,b,e), (c,d,e) \text{ admissible}\}$$ forms a basis of $\Hom(a\otimes b, c \otimes d)$, and moreover, $\langle B_{ab,cd}^e,B_{cd,ab}^{e'}\rangle = \delta_{e,e'} d_{e}$.

\begin{definition}
\label{def:extendedcoloring}
An extended $\cross{\Cat}{G}$-coloring of the pair $(M, \T)$ is a triple $\hat{F} = (g,f,t)$, $g: \T^1 \rightarrow G$, $f: \T^2 \rightarrow \L{\cross{\Cat}{G}}$, $t: \T^3 \rightarrow \L{\cross{\Cat}{G}}$, such that,
 \begin{enumerate}
   \item for any $\beta = (012) \in \T^2$, $f(012) \in \Cat_{\overline{g}_{02}g_{01}g_{12}}$;
   \item for any $\tau = (0123) \in \T^3$, $t(0123) \in \Cat_{\overline{g}_{03}g_{01}g_{12}g_{23}}$.
 \end{enumerate}
\end{definition}

\begin{figure}
\centering
 \begin{tikzpicture}[scale = 0.5]
 \begin{scope}
  \ShortIshape{0}{0}{f_{023}}{}{f_{013}}{f_{123}}{}
  \draw (2,-1) node{{\tiny $\lsuprsub{\overline{g}_{23}}{f}{012}$}};
  \draw (0,-2) node{($a$)};
  \draw (-0.8,0.5) node{{\tiny $t_{0123}$}};
\end{scope}
\draw (3, 0.5) node{$=$};
\begin{scope}[xshift = 6cm]
  \IshapeL{0}{0}{t_{0123}}{-}{+}
  \draw (0,-2) node{($b$)};
\end{scope}
\begin{scope}[xshift = 12cm]
  \ShortIshape{0}{0}{f_{013}}{f_{123}}{f_{023}}{}{}
  \draw (2,2) node{{\tiny $\lsuprsub{\overline{g}_{23}}{f}{012}$}};
  \draw (0,-2) node{($c$)};
  \draw (3, 0.5) node{$=$};
  \draw (-0.8,0.5) node{{\tiny $t_{0123}$}};
\end{scope}
\begin{scope}[xshift = 18cm]
  \IshapeL{0}{0}{t_{0123}}{+}{-}
  \draw (0,-2) node{($d$)};
\end{scope}
 \end{tikzpicture}
 \caption{$(a)\colon B_{f_{023}\lsuprsub{\overline{g}_{23}}{f}{012}, f_{013} f_{123}}^{t_{0123}} $; $(b)\colon$ short notation; \\ $(c)\colon B_{f_{013}f_{123}, f_{023}\lsuprsub{\overline{g}_{23}}{f}{012}}^{t_{0123}} $; $(d)\colon$ short notation.}\label{fig:3-simplexbasis}
\end{figure}
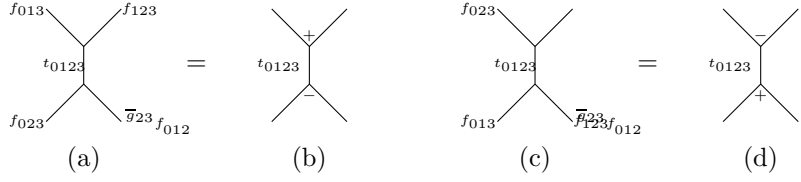

As before, we choose arbitrarily a representative for each $f(012)$ and $t(0123)$, and denote these representatives as $f_{012} $ and $ t_{0123}$ or simply as $012$ and $0123$. Given an extended coloring $\hat{F} = (g,f,t)$, then $F = (g,f)$ is a coloring according to Definition \ref{def:coloring}. As noted in Section \ref{subsec:partition1}, for each $3$-simplex $\tau = (0123)$, $f_{023} \otimes \lsuprsub{\overline{g}_{23}}{f}{012}$ and $f_{013} \otimes f_{123}$ are both in the sector $\Cat_{\overline{g}_{03}g_{01}g_{12}g_{23}}$, which is why in Definition \ref{def:extendedcoloring} we also require $t(0123)$ to be in the sector $\Cat_{\overline{g}_{03}g_{01}g_{12}g_{23}}$. More explicitly, as $t$ varies, the set
\begin{align*}
\Big\{B_{f_{023}\lsuprsub{\overline{g}_{23}}{f}{012}, f_{013} f_{123}}^{t_{0123}}\colon \ t_{0123} \in \Cat_{\overline{g}_{03}g_{01}g_{12}g_{23}}, \ t_{0123} \in \Ltil{\cross{\Cat}{G}}\Big\}
\end{align*}
forms a basis of $V_F^{+}(0123)$. Similarly, the set
\begin{align*}
\Big\{B_{f_{013}f_{123}, f_{023}\lsuprsub{\overline{g}_{23}}{f}{012}}^{t_{0123}}\colon \ t_{0123} \in \Cat_{\overline{g}_{03}g_{01}g_{12}g_{23}},\ t_{0123} \in \Ltil{\cross{\Cat}{G}}\Big\}
\end{align*}
forms a basis of $V_F^{-}(0123)$. Graphically, these two sets of bases are represented in Figure \ref{fig:3-simplexbasis}. Note that
\begin{align*}
\Big\langle B_{f_{023}\lsuprsub{\overline{g}_{23}}{f}{012}, f_{013} f_{123}}^{t_{0123}}\ ,\ B_{f_{013}f_{123}, f_{023}\lsuprsub{\overline{g}_{23}}{f}{012}}^{t'_{0123}} \Big\rangle &= d_{t_{0123}} \delta_{t_{0123}, t'_{0123}}.
\end{align*}

\begin{figure}
\centerline{\includegraphics{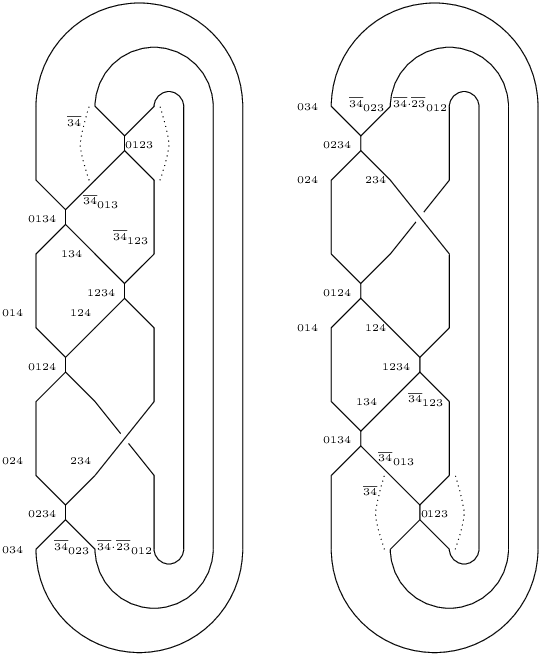}}
\caption{(Left): $\hat{Z}^{+}_{\hat{F}}(01234)$; (Right): $\hat{Z}^{-}_{\hat{F}}(01234)$}
\label{fig:25j}
\end{figure}

\begin{figure}
     \centering
     \includegraphics[width = 0.15\textwidth]{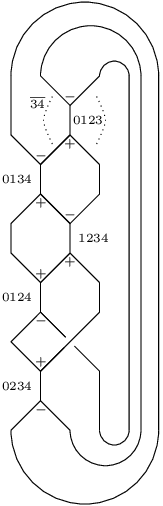}
     \qquad \qquad
     \includegraphics[width = 0.15\textwidth]{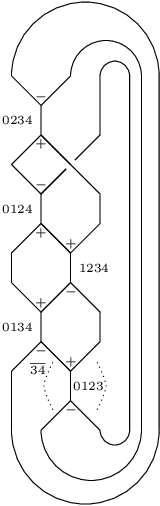}
\caption{Compact form: (Left): $\hat{Z}^{+}_{\hat{F}}(01234)$; (Right): $\hat{Z}^{-}_{\hat{F}}(01234)$}\label{fig:25jcompact}
\end{figure}

Rewriting the formula in Proposition \ref{prop:partition2} under this $B$-basis, we obtain a state-sum model. Explicitly, for each $\sigma \in \T^4$, the $25j$-symbol $\hat{Z}^{\epsilon(\sigma)}_{\hat{F}}(\sigma)$ is defined to be the evaluation of the graph diagram in Figure \ref{fig:25j} (Left) if $\epsilon(\sigma) = +$ and Figure \ref{fig:25j} (Right) if $\epsilon(\sigma) = -$. As in Section \ref{subsec:partition1}, the three long vertical segments representing the trace in each diagram are implicitly directed upwards and all other segments are directed downwards. And the $\overline{34}$ symbol in each diagram means the action of $\overline{34}$ on the morphism enclosed by the parenthesis. In the diagram we also have dropped the letters $g, f, t$, and thus the labels, e.g., $024, 0234$, etc. denote the colors on the corresponding simplices. In a more compact form, $\hat{Z}^{\epsilon(\sigma)}_{\hat{F}}(\sigma)$  is defined by the diagram in Figure \ref{fig:25jcompact} (left) in the case $\epsilon(\sigma) = +$ and by Figure \ref{fig:25jcompact} (right) otherwise.

\begin{proposition}[Definition of Partition Function III]
\label{prop:partition3}
The partition function has the state sum model,
\begin{align}
\label{equ:statesumdef}
Z_{\cross{\Cat}{G}}(M;\T) &= \sum\limits_{\hat{F} = (g,f,t)} \frac{(D^2/|G|)^{|\T^0|} \left(\prod\limits_{\beta \in \T^2}d_{f(\beta)}\right) \left(\prod\limits_{\sigma \in \T^4}\hat{Z}^{\epsilon(\sigma)}_{\hat{F}}(\sigma) \right) }{(D^2)^{|\T^1|} \left(\prod\limits_{\tau \in \T^3}d_{t(\tau)}\right)}
\end{align}
\begin{proof}
The right hand side of Equation \ref{equ:statesumdef} can be written as
\begin{align*}
\sum\limits_{F = (g,f)} \frac{(D^2/|G|)^{|\T^0|} \prod\limits_{\beta \in \T^2}d_{f(\beta)}  }{(D^2)^{|\T^1|} } \sum\limits_{\hat{F} = (F,t)}\frac{\prod\limits_{\sigma \in \T^4}\hat{Z}^{\epsilon(\sigma)}_{\hat{F}}(\sigma)}{\prod\limits_{\tau \in \T^3}d_{t(\tau)}}
\end{align*}
For each $3$-simplex $\tau = (0123)$, we have
\begin{align*}
\phi_{F, \tau} &= \sum_{t_{0123}} \, \frac{1}{d_{t_{0123}}}\, B_{f_{013}f_{123}, f_{023}\lsuprsub{\overline{g}_{23}}{f}{012}}^{t_{0123}} \ \otimes \ B_{f_{023}\lsuprsub{\overline{g}_{23}}{f}{012}, f_{013} f_{123}}^{t_{0123}}.
\end{align*}
Thus, for a fixed coloring $F = (g,f)$,
\begin{align*}
\left(\bigotimes\limits_{\sigma \in \T^4}\tilde{Z}^{\epsilon(\sigma)}_{F}(\sigma)\right)(\phi_{F}) &= \sum\limits_{\hat{F} = (F,t)}\frac{\prod\limits_{\sigma \in \T^4}\hat{Z}^{\epsilon(\sigma)}_{\hat{F}}(\sigma) }{\prod\limits_{\tau \in \T^3}d_{t(\tau)}}.
\end{align*}
\end{proof}
\end{proposition}

\begin{figure}
\centering
\begin{tikzpicture}[scale = 0.5]
\draw (0,2.5) node{$\sum\limits_{a,d}$};
\begin{scope}[xshift = 2cm]
\draw (0,0) -- (2,2) -- (2,3) -- (0,5);
\draw (2,0) -- (1,1);
\draw (4,0) -- (2,2) -- (2,3) -- (4,5);
\draw (1,4) -- (2,5);
\draw (0,6) node{{\tiny $034$}};
\draw (2,6) node{{\tiny $\lsup{\overline{34}}{023}$}};
\draw (4,6) node{{\tiny $\lsup{\overline{34}\cdot\overline{23}}{012}$}};
\draw (0,-1) node{{\tiny $034$}};
\draw (2,-1) node{{\tiny $\lsup{\overline{34}}{023}$}};
\draw (4,-1) node{{\tiny $\lsup{\overline{34}\cdot\overline{23}}{012}$}};
\draw (2.5,2.5) node{{\tiny $d$}};
\draw (1.7,1.3) node{{\tiny $a$}};
\draw (1.7,3.7) node{{\tiny $a$}};
\end{scope}

\draw (8,2.5) node{$=$};

\begin{scope}[xshift = 10cm]
\draw (0,0) -- (0,5);
\draw (2,0) -- (2,5);
\draw (4,0) -- (4,5);
\draw (0,6) node{{\tiny $034$}};
\draw (2,6) node{{\tiny $\lsup{\overline{34}}{023}$}};
\draw (4,6) node{{\tiny $\lsup{\overline{34}\cdot\overline{23}}{012}$}};
\draw (0,-1) node{{\tiny $034$}};
\draw (2,-1) node{{\tiny $\lsup{\overline{34}}{023}$}};
\draw (4,-1) node{{\tiny $\lsup{\overline{34}\cdot\overline{23}}{012}$}};
\end{scope}

\end{tikzpicture}
\caption{Identity morphism}\label{fig:Fbasisidentity}
\end{figure}
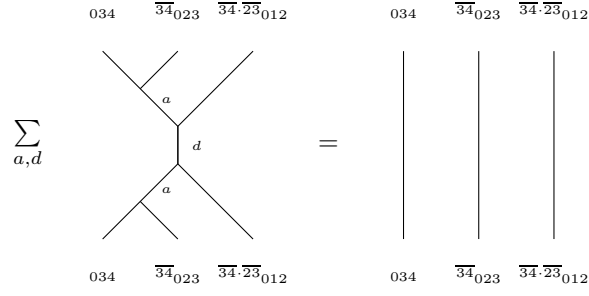

The $25j$-symbol $\hat{Z}^{\epsilon(\sigma)}_{\hat{F}}(\sigma)$ in Figure \ref{fig:25j} can be expressed as a concrete formula in terms of the data $(N_{ab}^c,\ F^{abc}_{d;nm},\  U_g(a,b;c),\  \eta_a(g,h),\  R_c^{ab})$ given in Section \ref{subsec:skeletonization}. Explicitly, let $\sigma = (01234)$ and assume an extended coloring $\hat{F}$ has been assigned to the triangulation, then $\hat{Z}^{+}_{\hat{F}}(\sigma)$ and $\hat{Z}^{-}_{\hat{F}}(\sigma)$ are given by Equations \ref{equ:Zpformua} and \ref{equ:Zmformula}, respectively. The way to obtain these formulas is by inserting the identity morphism shown in Figure \ref{fig:Fbasisidentity} to the bottom of the diagrams in Figure \ref{fig:25j} and using the diagrammatic equations in Section \ref{subsec:skeletonization} to simplify the diagrams.

 \begin{align}\label{equ:Zpformua}
 \hat{Z}^{+}_{\hat{F}}(\sigma) \ = \ \sum\limits_{d,a} &{}\ F_{d;a,0234}^{024,234,\lsup{\overline{34}\cdot\overline{23}}{012}}\, \eta^{-1}_{012}(\overline{34},\overline{23})\, \eta^{-1}_{012}(\overline{24}\cdot 23\cdot 34, \overline{34}\cdot\overline{23}) \, R_{a}^{234, \lsup{\overline{34}\cdot \overline{23}}{012}} \,   \nonumber\\
   &{}\ (F^{-1})_{d;0124,a}^{024,\lsup{\overline{24}}{012}, 234}\, F_{d;1234,0124}^{014,124,234} \, (F^{-1})_{d;0134,1234}^{014,134,\lsup{\overline{34}}{123}} \, F_{d;\lsup{\overline{34}}{0123},0134}^{034, \lsup{\overline{34}}{013}, \lsup{\overline{34}}{123}}\nonumber \\
   &{}\ U_{\overline{34}}(023,\lsup{\overline{23}}{012};0123) \, U^{-1}_{\overline{34}}(013,123;0123)\, (F^{-1})_{d;0234,\lsup{\overline{34}}{0123}}^{034,\lsup{\overline{34}}{023},\lsup{\overline{34}\cdot\overline{23}}{012}} \, d_d
 \end{align}

%

 \begin{align}\label{equ:Zmformula}
 \hat{Z}^{-}_{\hat{F}}(\sigma) \ = \ \sum\limits_{d,a} &{}\  U^{-1}_{\overline{34}}(023,\lsup{\overline{23}}{012};0123) \, U_{\overline{34}}(013,123;0123)\, F_{d;\lsup{\overline{34}}{0123},0234}^{034,\lsup{\overline{34}}{023},\lsup{\overline{34}\cdot\overline{23}}{012}}\,   \nonumber\\
   &{}\  (F^{-1})_{d;0134,\lsup{\overline{34}}{0123}}^{034, \lsup{\overline{34}}{013}, \lsup{\overline{34}}{123}} \, F_{d;1234,0134}^{014,134,\lsup{\overline{34}}{123}} \, (F^{-1})_{d;0124,1234}^{014,124,234} \, F_{d;a,0124}^{024,\lsup{\overline{24}}{012}, 234} \nonumber \\
   &{}\ \eta_{012}(\overline{24}\cdot 23\cdot 34, \overline{34}\cdot\overline{23}) (R_{a}^{234, \lsup{\overline{34}\cdot \overline{23}}{012}})^{-1}\eta_{012}(\overline{34},\overline{23})(F^{-1})_{d;0234,a}^{024,234,\lsup{\overline{34}\cdot\overline{23}}{012}}\, d_d
 \end{align}

\section{Examples and Variations}
\label{sec:example}
In this section, we give several examples of the partition function defined in Section \ref{sec:partition}. Some of the known invariants in literature emerge as special cases of our construction. At the end of this section, a variation of the partition function is introduced when the \CatNameShort{$G$} has a trivial grading. In this case, a cocycle in $H^4(G, U(1))$ can be included to produce a different invariant. A particular case reduces to the twisted Dijkgraaf-Witten invariant.
\subsection{Crane-Yetter Invariant}
The Crane-Yetter invariant was introduced in \cite{crane1993categorical}, where the authors gave a state sum construction of $4$-manifold invariants with the modular tensor category $Rep(U_q(sl_2))$, where $Rep(U_q(sl_2))$ is the category of representations of $U_q(sl_2)$ with $q$ some principal $4r$-th root of unity. Later on the construction was generalized to any semisimple ribbon category $\Cat$\cite{crane1997state} \footnote{In \cite{crane1997state}, such a category was called semisimple tortile category.}. This generalized state sum invariant is still called the Crane-Yetter invariant denoted by $\CY_{\Cat}(\cdot)$.

For a \CatNameShort{$G$} $\cross{\Cat}{G} = \bigoplus\limits_{g \in G}\Cat_g$, if $G = \{e\}$ is the trivial group, then $\cross{\Cat}{G}$ has only one sector $\Cat = \Cat_e$ and it is a ribbon fusion  category as noted in Section \ref{subsec:GBSFC_Def}. In this case, the only color on a $1$-simplex is the unit $e$, so we can just assume there is no color at all on $1$-simplices. The colors on each $2$- and $3$-simplex run through a complete set of representatives in $\Cat$. For each colored $4$-simplex, the partition function is given by Figure \ref{fig:25j} where all the relevant group elements and group actions are trivial, and this is equal to the $15j$-symbol defined in \cite{crane1997state}. The state-sum formula is given by
\begin{align}
\label{equ:statesumGtrivial}
Z_{\cross{\Cat}{\{e\}}}(M;\T) &= \sum\limits_{\hat{F} = (f,t)} \frac{(D^2)^{|\T^0|} \left(\prod\limits_{\beta \in \T^2}d_{f(\beta)}\right) \left(\prod\limits_{\sigma \in \T^4}\hat{Z}^{\epsilon(\sigma)}_{\hat{F}}(\sigma) \right) }{(D^2)^{|\T^1|} \left(\prod\limits_{\tau \in \T^3}d_{t(\tau)}\right)}.
\end{align}

It is then direct to see that the resulting partition function $Z_{\cross{\Cat}{\{e\}}}(M)$ is exactly the Crane-Yetter invariant $\CY_{\Cat}(M)$ (up to the translation of some conventions).

\subsection{Yetter's Invariant from Homotopy $2$-types}
\label{subsec:yetter}
A (strict) categorical group is a rigid tensor category $\G$ such that $(\otimes, \unit, (\cdot)^{*})$ satisfies the axioms of a group strictly and that every morphism is invertible. Note that here the morphism spaces are not required to be a vector space. Namely, in a categorical group, all the structural isomorphisms $a, l, r$ are identity maps and $A^* \otimes A = A \otimes A^* = \unit, \  f \otimes (f^{-1})^* = Id, \ (f^{-1})^* \otimes f  = Id$ for any object $A$ and morphism $f$.
There is a one-to-one correspondence between categorical groups and crossed modules which are defined below.

Let $G,H$ be two finite groups. A crossed module is a quadruple $(H, G, \rho, \phi)$ where $\rho\colon H \longrightarrow G$ is a group morphism and $\phi\colon G \times H \longrightarrow H$ is a group action of $G$ on $H$, such that the following two conditions are satisfied.
\begin{enumerate}
 \item $\rho$ commutes with the $G$-action, i.e., for any $g \in G, h \in H,$ $\rho(\phi(g,h)) = \lsup{g}{\rho(h)},$ where the right hand side is the conjugation action.
 \item The action $\phi$ extends the conjugation action of $H$ on itself, i.e., for any $h,h' \in H$, $\phi(\rho(h'),h) = \lsup{h'}{h}$.
\end{enumerate}
 We will usually write the action as $\phi(g,h) = \lsup{g}{h}$. The inverse of a group element $g$ is also denoted by $\overline{g}$.

Given a categorical group $\G$, let $G = \G^0$ and $H = \sqcup_{g \in G} \Hom(\unit,g)$. It is clear that $G$ is a group with the product $\otimes$, the unit $e = \unit$, and the inverse $\overline{(\cdot)} = (\cdot)^{*}$. Define $\rho\colon H \longrightarrow G$ to be the target map, namely $\rho(h)$ is the target of $g$, and define the $G$-action by $\phi(g,h) = Id_g \otimes h \otimes Id_{\overline{g}}$. We leave it as an exercise to check that $H$ is also a group, and $(H,G, \rho,\phi)$ is a crossed module, which we denote by $\Mod(\G)$.

For the converse direction, given a crossed module $(H,G, \rho,\phi)$, define the categorical group $\G(H,G,\rho,\phi)$ as follows.
\begin{enumerate}
  \item Objects are elements of $G$. For $g_1, g_2 \in G$, $\Hom(g_1,g_2):= \rho^{-1}(\overline{g_1}g_2) \subset H$. The composition of morphisms are the product in $H$.
  \item For $h \in \Hom(g_1,g_2), h' \in \Hom(g_1',g_2') $, $g_1 \otimes g_1' := g_1g_1'$, and $h \otimes h' := \lsup{\overline{g_1'}}{h}h'$. The unit $\unit$ is the unit in $G$. The structure isomorphisms $a, l,r $ are identity maps.
  \item For $g \in G$, the dual $g^*:= \overline{g}$, and the $b_g, d_g$ are identities.
\end{enumerate}

\begin{proposition}
If $(H,G, \rho,\phi)$ is a crossed module and $\G$ is a categorical group, then $\Mod(\G(H,G,\rho,\phi)) = (H,G,\rho,\phi)$ and $\G(\Mod(\G)) = \G$.
\end{proposition}
\begin{proof}
Direct verification.
\end{proof}

On the other hand, from a crossed module $(H, G,\rho, \phi)$ we can also construct a \CatNameShort{$G$} denoted by $\D = \D(H,G,\rho,\phi) = \bigoplus\limits_{g \in G}\Cat_g$. As a category, $\D(H,G,\rho,\phi)$ is the same as $\Vect_H$, the category of $H$-graded finite dimensional vector spaces. We identify the simple objects in $\D$ with elements of $H$. For each $g \in G$, the $g$-sector $\Cat_g$ is spanned by all simple objects in $\rho^{-1}(g)$. This defines a $G$-grading on $\D$ due to the fact that $\rho$ is a group morphism. The $G$-action on simple objects of $\D$ is defined to be the action $\phi$ on $H$. It can be checked that this is a well-defined action, and the first condition in the definition of crossed modules guarantees that the $G$-action $g$ sends the $g'$-sector to the $\lsup{g}{g'}$-sector.

 Let $h \in \Cat_g, h' \in \Cat_{g'}$, namely, $\rho(h) = g, \rho(h') = g'$, then $h \otimes h' = hh' \in \Cat_{gg'},$ and $\lsup{g}{h'} \otimes h = \lsup{\rho(h)}{h'} \otimes h = \lsup{h}{h'} \otimes h = hh' \in \Cat_{gg'}$, where the second equality is due to the second condition in the definition of crossed modules. We define the $G$-crossed braiding by the identity map, namely,
 \begin{align*}
 \begin{array}{cccc}
 c_{h,h'}\ := Id\ \colon & h \otimes h'& \longrightarrow &\lsup{g}{h'} \otimes h 
 \end{array}
 \end{align*}

 We now describe the partition function associated with $\D(H,G,\rho,\phi)$. Let $\hat{F} = (g,f,t)$ be an extended coloring. Then for each $2$-simplex $(012)$, $f_{012} \in \Cat_{\overline{g}_{02}g_{01}g_{12}}$. For any 3-simplex $(0123)$, if the space $V^{\pm}(0123)$ is to be nonzero, we need to have
 \begin{equation}
 \label{equ:crossmodule}
 f_{023} \lsup{\overline{g}_{23}}{f_{012}} = f_{013}f_{123} = t_{0123}
 \end{equation}
Thus the color on a $3$-simplex is uniquely determined by those on its boundary faces. Given an extended coloring such that the condition from Equation \ref{equ:crossmodule} is satisfied for every $3$-simplex, it is direct to see that for each $4$-simplex $\sigma$, $\hat{Z}^{\pm}(\sigma) = 1$ by Figure \ref{fig:25j}.

 \begin{definition}
 \label{def:crossmodulecoloring}
 Given a crossed module $(H,G,\rho,\phi)$ and an ordered triangulation $\T$ of $M$, an {\it admissible} coloring is a map $F= (g,f)$, $g\colon \T^1 \longrightarrow G$, $f\colon \T^2 \longrightarrow H$, such that,
 \begin{enumerate}
  \item for any $2$-simplex $(012)$, $\rho(f_{012}) = \overline{g}_{02}g_{01}g_{12}$;
  \item for any $3$-simplex $(0123)$, $f_{023} \lsup{\overline{g}_{23}}{f_{012}} = f_{013}f_{123}$.
 \end{enumerate}
 \end{definition}

 An admissible coloring is a $\G(H,G,\rho,\phi)$-color in the sense of \cite{yetter1993tqft}. The following property shows actually the partition function associated with $\D(H,G,\rho,\phi)$ is exactly the Yetter's invariant $\Y_{\G(H,G,\rho,\phi)}$ associated with $\G(H,G,\rho,\phi)$ in \cite{yetter1993tqft}.
 \begin{proposition}
 \label{prop:yetter}
 Let $\D = \D(H,G,\rho,\phi)$ be the \CatNameShort{$G$} obtained from a crossed module, then
 \begin{alignat}{3}
 Z_{\D}(M;\T)\quad &&= \quad\frac{|H|^{|\T^0|-|\T^1|}}{|G|^{|\T^0|}} \# (H,G,\rho,\phi) \quad&&=\quad \Y_{\G(H,G,\rho,\phi)}(M),
 \end{alignat}
 where $\# (H,G,\rho,\phi)$ is the number of {\it admissible} colorings.

 \end{proposition}
 \begin{proof}
 In $\D$, it is clear that the quantum dimension of each simple object is $1$, and the total dimension square $D^2 = |H|$. Then the first equality follows from Equation \ref{equ:statesumdef}. The second equality follows from \cite{yetter1993tqft}.
 \end{proof}

 If $H = \{e\}$ is the trivial group, there is a unique group morphism $\rho_0$ from $H$ to $G$ and a unique action (the trivial action) $\phi_0$ of $G$ on $H$. Then according to Definition \ref{def:crossmodulecoloring}, an {\it admissible} coloring is simply a map $g\colon \T^1 \longrightarrow G$, such that for each $2$-simplex $(012)$, $\overline{g}_{02}g_{01}g_{12} = e$. In this case the partition function is reduced to the untwisted Dijkgraaf-Witten invariant $\DW_G(M)$ \cite{dijkgraaf1990topological}.
 \begin{proposition}
 Let $\D_0 = \D(\{e\},G, \rho_0,\phi_0)$, then
 \begin{alignat}{4}
 Z_{\D_0}(M;\T)\quad &&= \quad\frac{1}{|G|^{|\T^0|}} \# (\{e\},G,\rho_0,\phi_0) \quad &&= \quad \frac{|\Hom(\pi_1(M),G)|}{|\Hom(\pi_0(M),G)|} \quad &&= \quad\DW_G(M),
 \end{alignat}
 where $\pi_0(M)$ is the set of connected components of $M$ and $\Hom(\pi_0(M),G)$ is the set of maps from $\pi_0(M)$ to $G$.

 \end{proposition}
 \begin{proof}
 The first equality is by Proposition \ref{prop:yetter}.

 Note that the two sides of the second equality are both multiplicative with respect to disjoint union of connected components. Thus it suffices to prove the equality for a connected manifold $M$, namely, $|\pi_0(M)| = 1$.

Choose a maximal spanning tree $K$, which is a sub complex of $\T^1$ with $|\T^0|-1$ edges. Then it is easy to see that there is a $|G|^{|\T^0|-1}$ to one correspondence between the set of {\it admissible} colorings and $\Hom(\pi_1(M),G)$.
 \end{proof}

 \subsection{Trivial $G$-grading with a Trivial $G$-action}
 \label{subsec:trivialgradingtrivialaction}
 Given a ribbon fusion category $\Cat$ and a finite group $G$, we can construct a \CatNameShort{$G$} $\cross{\Cat}{G} = \bigoplus\limits_{g \in G} \Cat_g$, where $\Cat_g = \Cat$ if $g = e$, and $\Cat_g = 0$ otherwise, and $G$ acts on $\cross{\Cat}{G}$ by the identity functor. We consider the partition function associated to $\cross{\Cat}{G}$.

 Since the nontrivial part of $\cross{\Cat}{G}$ is constrained in the trivial sector, the coloring $g$ on $1$-simplices needs to satisfy the Dijkgraaf-Witten coloring rule, namely, $\overline{g}_{02}g_{01}g_{12} = e$ for each $2$-simplex $(012)$. The colorings on $2$- and $3$-simplices are independent of those on $1$-simplices, and they are the same as the colorings of Crane-Yetter model. Moreover, the partition function corresponding to each $4$-simplex does not depend on the colors on $1$-simplices and is the same as that of Crane-Yetter since the group action here is trivial. We therefore have,
 \begin{proposition}
 If $\cross{\Cat}{G}$ has a trivial $G$-grading and a trivial $G$-action where the trivial sector of $\cross{\Cat}{G}$ is $\Cat$, then $Z_{\cross{\Cat}{G}}(M;\T) = \CY_{\Cat}(M) \DW_{G}(M)$.
 \end{proposition}
 \begin{proof}
 By Equation \ref{equ:statesumdef} and the argument above,
 \begin{align*}
Z_{\cross{\Cat}{G}}(M;\T) &= \sum\limits_{\hat{F} = (g,f,t)} \frac{(D^2/|G|)^{|\T^0|} \left(\prod\limits_{\beta \in \T^2}d_{f(\beta)}\right) \left(\prod\limits_{\sigma \in \T^4}\hat{Z}^{\epsilon(\sigma)}_{\hat{F}}(\sigma) \right) }{(D^2)^{|\T^1|} \left(\prod\limits_{\tau \in \T^3}d_{t(\tau)}\right)} \\
             &= \sum\limits_{g} \frac{1}{|G|^{|\T^0|}}\sum\limits_{f,t} \frac{(D^2)^{|\T^0|} \left(\prod\limits_{\beta \in \T^2}d_{f(\beta)}\right) \left(\prod\limits_{\sigma \in \T^4}\hat{Z}^{\epsilon(\sigma)}(\sigma) \right) }{(D^2)^{|\T^1|} \left(\prod\limits_{\tau \in \T^3}d_{t(\tau)}\right)}\\
             &= \sum\limits_{g} \frac{1}{|G|^{|\T^0|}} \CY_{\Cat}(M) \\
             &= \DW_{G}(M)\CY_{\Cat}(M),
\end{align*}
where $\hat{Z}^{\epsilon(\sigma)}(\sigma)$ is the $25j$-symbol of $\sigma$ for which we have hidden the dependence on the coloring.

 \end{proof}

 In Section \ref{subsec:variation}, we will give a variation of the construction of the partition function, in which the twisted Dijkgraaf-Witten invariant appears as a special case.

\subsection{Variation: Trivial $G$-grading}
\label{subsec:variation}
In this section, we study a variation of the construction of the partition function, which produces a different invariant of $4$-manifolds.

Let $\Cat$ be a ribbon fusion category and let $\cross{\Cat}{G} = \bigoplus\limits_{g \in G} \Cat_g$ be a \CatNameShort{$G$} with a trivial $G$-grading, namely, $\Cat_e = \Cat$ and $\Cat_g = 0$ for $g \neq e$. Thus $\cross{\Cat}{G}$ is a ribbon fusion category endowed with a $G$-action. In this case, we show that a $4$-cocycle $\omega \in H^4(G,U(1))$ can be introduced in the construction of the partition function.

Since the $G$-grading is trivial, the colors $g$ on $1$-simplices must satisfy the condition $\bar{g}_{02}g_{01}g_{12} = e$ for each $2$-simplex $(012)$. The colors on $2$- and $3$-simplices run through a complete set of representatives in $\Cat$. For each fixed coloring, we will introduce an $\omega$ factor to the partition function. Explicitly, using the notations in Section \ref{sec:partition}, the new invariant can be defined as follows.

\begin{definition}
\label{def:newpartition}
Given a \CatNameShort{$G$} $\cross{\Cat}{G}$ with a trivial $G$-grading, a $4$-cocycle $\omega \in H^4(G,U(1))$, and an ordered triangulation $\T$ of a $4$-manifold $M$, the partition function $Z_{\cross{\Cat}{G},\omega}(M;\T)$ of the pair $(M,\T)$ is defined by,
\begin{equation}
Z_{\cross{\Cat}{G},\omega}(M;\T) = \sum\limits_{F = (g,f)} \frac{(D^2/|G|)^{|\T^0|} \left(\prod\limits_{\beta \in \T^2}d_{f(\beta)}\right) \Tr\left(\bigotimes\limits_{\sigma \in \T^4}Z^{\epsilon(\sigma)}_{F}(\sigma)\omega(g_{01},g_{12},g_{23},g_{34})^{\epsilon(\sigma)}\right)}{(D^2)^{|\T^1|}} ,
\end{equation}
or as a state sum model,

\begin{equation}
\label{equ:newstatesumdef}
Z_{\cross{\Cat}{G},\omega}(M;\T) = \sum\limits_{\hat{F} = (g,f,t)} \frac{(D^2/|G|)^{|\T^0|} \left(\prod\limits_{\beta \in \T^2}d_{f(\beta)}\right) \left(\prod\limits_{\sigma \in \T^4}\hat{Z}^{\epsilon(\sigma)}_{\hat{F}}(\sigma)\omega(g_{01},g_{12},g_{23},g_{34})^{\epsilon(\sigma)} \right) }{(D^2)^{|\T^1|} \left(\prod\limits_{\tau \in \T^3}d_{t(\tau)}\right)},
\end{equation}
where we write each $\sigma$ as $(01234)$ with the induced ordering, and  $\omega(\cdots)^{\epsilon(\sigma)}$ is $\omega(\cdots)$ if $\epsilon(\sigma) = +$ and $\omega(\cdots)^{-1} $ otherwise.
\end{definition}

\begin{theorem}
$Z_{\cross{\Cat}{G}, \omega}(M,\T)$ is independent of the choice of the triangulation $\T$, and is thus an invariant of closed smooth oriented $4$-manifolds.
\end{theorem}

The proof of invariance of $Z_{\cross{\Cat}{G}, \omega}(M;\T)$ can be processed in a similar way as that of $Z_{\cross{\Cat}{G}}(M;\T)$ with slight modifications, so we will omit the details. For instance, to show that it is invariant under the $3$-$3$ Pachner move, Equation \ref{equ:3-3} will be replaced by Equation \ref{equ:3-3 new},

\begin{align}
\label{equ:3-3 new}
    & \sum\limits_{\substack{I_2, \\ I_3}} \Big( \frac{d_{024}}{d_{0124}d_{0234}d_{0245}} \hat{Z}^{+}(01234)\hat{Z}^{+}(01245)\hat{Z}^{+}(02345) \nonumber \\
    &{} \qquad \cdot \omega(01,12,23,34)\omega(01,12,24,45)\omega(02,23,34,45) \Big)\nonumber \\
    =&
 \sum\limits_{\substack{I'_2, \\ I'_3}} \Big( \frac{d_{135}}{d_{0135}d_{1235}d_{1345}} \hat{Z}^{+}(01235)\hat{Z}^{+}(01345)\hat{Z}^{+}(12345) \nonumber \\
   &{} \qquad \cdot \omega(01,12,23,35)\omega(01,13,34,45)\omega(12,23,34,45)\Big)
\end{align}

In Equation \ref{equ:3-3 new}, the product of the three $\omega$ factors on the left hand side is equal to the product of the three factors on the right hand side, and this is precisely because $\omega$ is a $4$-cocycle and the colors on $1$-simplices satisfy $\bar{g}_{02}g_{01}g_{12}=e$ for each $2$-simplex $(012)$. Cancelling the $\omega$ factors on both sides, we get back to Equation \ref{equ:3-3}.

\begin{proposition}
If $\cross{\Cat}{G} = \Cat$ is a \CatNameShort{$G$} with a trivial $G$-grading and a trivial $G$-action, and $\omega \in H^4(G,U(1))$, then
\begin{equation}
Z_{\cross{\Cat}{G},\omega}(M;\T) = \DW_{G}^{\omega}(M) \CY_{\Cat}(M),
\end{equation}
where $\DW_{G}^{\omega}(M)$ is the twisted Dijkgraaf-Witten invariant \cite{dijkgraaf1990topological}. In particular, If $\cross{\Cat}{G} = \Vect$, then $Z_{\cross{\Cat}{G},\omega} = \DW_{G}^{\omega}(M)$.
\end{proposition}
\begin{proof}
A similar argument as that in Section \ref{subsec:trivialgradingtrivialaction} shows that $Z_{\cross{\Cat}{G},\omega} $ can be written as follows:
 \begin{align*}
Z_{\cross{\Cat}{G},\omega}(M;\T) &= \sum\limits_{\hat{F} = (g,f,t)} \frac{(D^2/|G|)^{|\T^0|} \left(\prod\limits_{\beta \in \T^2}d_{f(\beta)}\right) \left(\prod\limits_{\sigma \in \T^4}\hat{Z}^{\epsilon(\sigma)}_{\hat{F}}(\sigma)\omega(g_{01},g_{12},g_{23},g_{34})^{\epsilon(\sigma)} \right) }{(D^2)^{|\T^1|} \left(\prod\limits_{\tau \in \T^3}d_{t(\tau)}\right)}, \\
             &= \sum\limits_{g} \frac{\prod\limits_{\sigma \in \T^4} \omega(g_{01},g_{12},g_{23},g_{34})^{\epsilon(\sigma)}}{|G|^{|\T^0|}}\sum\limits_{f,t} \frac{(D^2)^{|\T^0|} \left(\prod\limits_{\beta \in \T^2}d_{f(\beta)}\right) \left(\prod\limits_{\sigma \in \T^4}\hat{Z}^{\epsilon(\sigma)}(\sigma) \right) }{(D^2)^{|\T^1|} \left(\prod\limits_{\tau \in \T^3}d_{t(\tau)}\right)}\\
             &= \sum\limits_{g} \frac{\prod\limits_{\sigma \in \T^4} \omega(g_{01},g_{12},g_{23},g_{34})^{\epsilon(\sigma)}}{|G|^{|\T^0|}} \CY_{\Cat}(M) \\
             &= \DW_{G}^{\omega}(M)\CY_{\Cat}(M).
\end{align*}
\end{proof}

\section{Proof of Theorem \ref{thm:main}}
\label{sec:proofmain}
In this section we will prove in turn that $Z_{\cross{\Cat}{G}}(M;\T)$ is independent of the choice of representatives (Section \ref{sec:invariance representative}), the ordering of vertices of a triangulation (Section \ref{sec:invariance ordering}), and the choice of a triangulation (Section \ref{sec:invariance pachner}).

\subsection{Invariance under Choice of Representatives}
\label{sec:invariance representative}
 In defining the partition function $Z_{\cross{\Cat}{G}}(M;\T)$, for each coloring $F = (g,f)$, we chose arbitrarily a representative $f_{\beta}$ for the triangle color $f(\beta)$, $\beta \in \T^2$. In the following, we show that $\Tr\left(\bigotimes\limits_{\sigma\in \T^4} Z_{F}^{\epsilon(\sigma)}(\sigma)\right)$ is independent of the choice of representatives, and thus $Z_{\cross{\Cat}{G}}(M;\T)$ is also independent of the choice of representatives by Definition \ref{def:partition1}.

Let $\phi_A\colon A \overset{\sim}{\longrightarrow} A',\ \phi_B\colon B \overset{\sim}{\longrightarrow} B'$ be isomorphisms in $\cross{\Cat}{G}$. Denote by $T_{\phi_A}^{\phi_B}$ the following linear isomorphism,
\begin{align*}
\begin{array}{rrcl}
T_{\phi_A}^{\phi_B}\ \colon & \Hom(A,B)     &\overset{\sim}{\longrightarrow} &\Hom(A',B') \\
                            & \phi_{A}^{B}  &\longmapsto       & \phi_B\circ \phi_{A}^{B} \circ \phi_A^{-1}.\\
\end{array}
\end{align*}
Thus, Diagram \ref{diag:T_phi^psi} commutes.

\begin{equation}
\label{diag:T_phi^psi}
\insertimage{1}{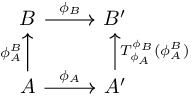}
\end{equation}
The isomorphism $T_{\phi_A}^{\phi_B}$ satisfies functorial properties, namely,
\begin{enumerate}
  \item If $\phi_C\colon C \overset{\sim}{\longrightarrow} C'$, then $T_{\phi_B}^{\phi_C}(\phi_B^C) \circ T_{\phi_A}^{\phi_B}(\phi_A^B) \ =\  T_{\phi_A}^{\phi_C}( \phi_B^C \circ \phi_A^B)$.
  \item $T_{\phi_A}^{\phi_A}(Id_A)\ =\ Id_{A'}$.
\end{enumerate}

\begin{lemma}
\label{lem:TpreservesPairing}
The isomorphism $T_{\phi_A}^{\phi_B}$ preserves the pairing, i.e., for $\phi_{A}^B \in \Hom(A, B), \phi_{B}^A \in \Hom(B, A)$, $$\langle \phi_{B}^{A} , \phi_{A}^{B}\rangle = \langle T_{\phi_B}^{\phi_A}(\phi_B^A), T_{\phi_A}^{\phi_B}(\phi_A^B)\rangle.$$
\begin{proof}
Direct verification.
\end{proof}
\end{lemma}

For the rest of the section, we fix a coloring $F = (g,f)$. For each $2$-simplex $(ijk) \in \T^2$, assume we have arbitrarily chosen two representatives of $f(ijk)$ denoted by $ijk$ and $ijk'$, respectively. To distinguish these two choices, we attach an apostrophe to all quantities related to the second the choice, e.g., $V_F^{+}(0123)', Z_{F}^{-}(01234)',$ etc.

For each $\beta = (ijk) \in \T^2,$ choose any isomorphism $\phi_{ijk}\colon ijk \longrightarrow ijk'$. If $\tau = (0123) \in \T^3$, let $\phi_{-\tau}$, $\phi_{+\tau}$ be the isomorphisms defined below.
\begin{align*}
\begin{array}{rccc}
\phi_{-\tau}:= \phi_{023} \otimes \lsup{\overline{23}}{(\phi_{012})}\ \colon &023 \otimes \lsup{\overline{23}}{012}& \overset{\sim}{\longrightarrow} & 023' \otimes \lsup{\overline{23}}{012'},\\
\phi_{+\tau}:= \phi_{013} \otimes \phi_{123}\ \colon & 013 \otimes 123 &\overset{\sim}{\longrightarrow} & 013' \otimes 123'. 
\end{array}
\end{align*}

Then we have the isomorphisms $T_{\phi_{-\tau}}^{\phi_{+\tau}}\colon V_F^{+}(\tau) \longrightarrow V_F^{+}(\tau)'$, $T_{\phi_{+\tau}}^{\phi_{-\tau}}\colon V_F^{-}(\tau) \longrightarrow V_F^{-}(\tau)'$. See Diagram \ref{diag:Ttau}. When it is clear from the context, we will drop the subscript/superscript and write $T_{\phi_{-\tau}}^{\phi_{+\tau}}$ as $T$.

\begin{equation}
\label{diag:Ttau}
\insertimage{1}{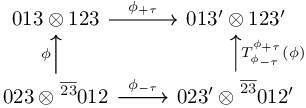}
\end{equation}

\begin{lemma}
\label{lem:TtildeZcommute}
For any $4$-simplex $\sigma = (01234) \in \T^4$, Diagrams \ref{diag:TtildeZplus} and \ref{diag:TtildeZminus} both commute.
\begin{equation}
\label{diag:TtildeZplus}
\insertimage{1}{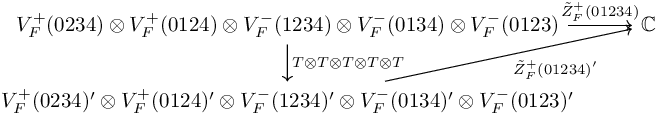}
\end{equation}
\begin{equation}
\label{diag:TtildeZminus}
\insertimage{1}{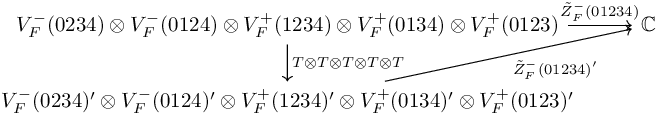}
\end{equation}
\end{lemma}
\begin{proof}
We only prove the case for Diagram \ref{diag:TtildeZplus}. The other case can be proved in the same way. Let $\phi_0 \otimes \phi_1 \otimes \phi_2 \otimes \phi_3 \otimes \phi_4 $ be in the domain of $\tilde{Z}^{+}_{F}(01234).$ Consider the following diagram, where it is not hard to see that the the two maps on the two sides of each horizontal arrow represent the same map.
\begin{equation}
\insertimage{1}{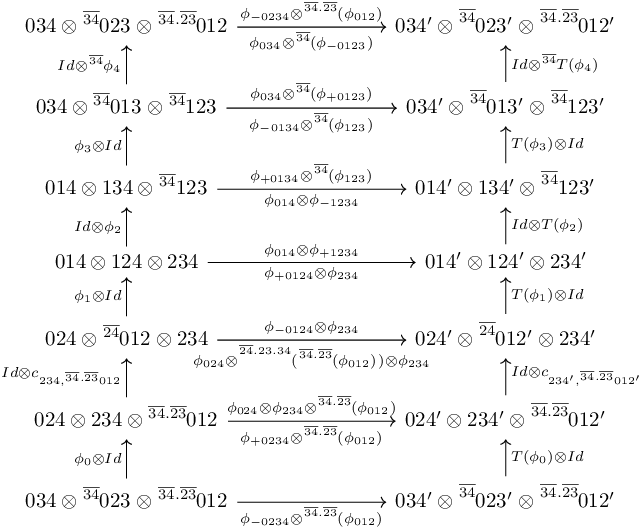}
\end{equation}
The second square diagram from the bottom commutes by the naturality of the $G$-crossed braiding. All other square diagrams commute by Diagram \ref{diag:Ttau}. Thus, the whole diagram commutes.

Denote the composition of the vertical maps on the left (resp. right) side by $L$ (resp. $R$), then $\tilde{Z}_F^{+}(01234)(\phi_0 \otimes \phi_1 \otimes \phi_2 \otimes \phi_3 \otimes \phi_4) = \Tr(L), \ \tilde{Z}_F^{+}(01234)'\left(T(\phi_0) \otimes T(\phi_1) \otimes T(\phi_2) \otimes T(\phi_3) \otimes T(\phi_4)\right) = \Tr(R)$. Since $L, R$ are conjugate by the above diagram, thus $\Tr(L) = \Tr(R)$ and we have $\tilde{Z}_F^{+}(01234) =  \tilde{Z}_F^{+}(01234)'(T^{\otimes 5}).$
\end{proof}

\begin{proposition}
\label{prop:TZ}
For any $\sigma = (01234) \in \T^4$, $Z_{F}^{\epsilon(\sigma)}(\sigma)$ commutes with $T$. More precisely, the following diagrams commute.
\begin{equation}
\label{diag:TZplus}
\insertimage{1}{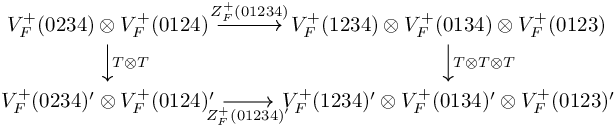}
\end{equation}

\begin{equation}
\label{diag:TZminus}
\insertimage{1}{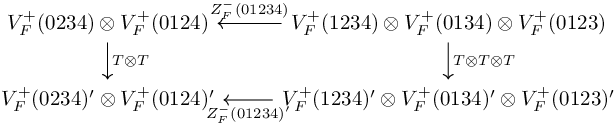}
\end{equation}
\begin{proof}
Again we only prove the case of Diagram \ref{diag:TZplus}.

Let $\phi_0 \otimes \phi_1 \in V^{+}_{F}(0234) \otimes V^{+}_{F}(0124),\ \phi_2 \otimes \phi_3 \otimes \phi_4 \in  V^{+}_{F}(1234) \otimes V^{+}_{F}(0134) \otimes V^{+}_{F}(0123),$ then
\begin{align}\label{equ:TTTZ}
  & \Big\langle (T \otimes T \otimes T) \circ Z_F^{+}(01234) (\phi_0 \otimes \phi_1), (T \otimes T \otimes T)(\phi_2 \otimes \phi_3 \otimes \phi_4) \Big\rangle  \nonumber\\
=& \Big\langle Z_F^{+}(01234) (\phi_0 \otimes \phi_1), \phi_2 \otimes \phi_3 \otimes \phi_4 \Big\rangle \nonumber \\
=& \tilde{Z}_F^{+}(01234)(\phi_0 \otimes \phi_1 \otimes \phi_2 \otimes \phi_3 \otimes \phi_4)
\end{align}
\begin{align}\label{equ:ZTT}
& \Big\langle  Z_F^{+}(01234)' \circ (T \otimes T) (\phi_0 \otimes \phi_1), (T \otimes T \otimes T)(\phi_2 \otimes \phi_3 \otimes \phi_4) \Big\rangle \nonumber \\
=& \tilde{Z}_F^{+}(01234)'\Big(T(\phi_0) \otimes T(\phi_1)\otimes T(\phi_2) \otimes T(\phi_3) \otimes T(\phi_4)\Big) \nonumber\\
=& \tilde{Z}_F^{+}(01234)' \circ (T \otimes T \otimes T \otimes T \otimes T)(\phi_0 \otimes \phi_1 \otimes \phi_2 \otimes \phi_3 \otimes \phi_4) \nonumber \\
=& \tilde{Z}_F^{+}(01234)(\phi_0 \otimes \phi_1 \otimes \phi_2 \otimes \phi_3 \otimes \phi_4)
\end{align}
The first \lq\lq = '' in Equation \ref{equ:TTTZ} is by Lemma \ref{lem:TpreservesPairing}, and the third \lq\lq = '' in Equation \ref{equ:ZTT} is by Lemma \ref{lem:TtildeZcommute}.

Since the equalities in Equation \ref{equ:TTTZ} and \ref{equ:ZTT} hold for any $\phi_i, i = 0,1,2,3,4$. this implies Diagram \ref{diag:TZplus} commutes.
\end{proof}
\end{proposition}

\begin{theorem}
Given a coloring $F = (g,f)$, we have $$\Tr\left(\bigotimes\limits_{\sigma \in \T^4} Z_F^{\epsilon(\sigma)}(\sigma)\right) = \Tr\left(\bigotimes\limits_{\sigma \in \T^4} Z_F^{\epsilon(\sigma)}(\sigma)'\right).$$ As a consequence, the partition function $Z_{\cross{\Cat}{G}}(M;\T)$ is independent on the choice of representatives for each triangle color.
\end{theorem}
\begin{proof}
It suffices to prove Diagram \ref{diag:TZtotal} commutes, which follows from Proposition \ref{prop:TZ} since each $V_F^{+}(\tau)$ is acted on by exactly one $Z_F^{\epsilon(\sigma)}(\sigma)$.
\begin{equation}
\label{diag:TZtotal}
\insertimage{1}{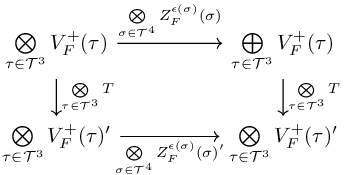}
\end{equation}
\end{proof}

\subsection{Invariance under Change of Ordering}
\label{sec:invariance ordering}
When defining the partition function, we require that the triangulation be equipped with an ordering on its vertices. In this subsection, it is shown that the partition function is actually independent of the ordering. We take Equation \ref{equ:partition2} as the definition of the partition function and notations from Section \ref{subsec:partition2} will be used.

Let $\Sym_n$ be the permutation group acting on the set $\{0,1, \cdots, n-1\}$ with the standard generators $s_i \in \Sym_n, i = 0,1, \cdots, n-2,$ swapping $i$ and $i+1$. For each $s \in \Sym_n$ and any nonempty subset $T \subset \{0,1, \cdots, n-1\}$, if we order the elements in $T$ and in $s(T)$, respectively, by $0,1,\cdots, |T|-1$ according to their relative order, then we get a permutation $s_T$, called the restriction of $s$ on $T$, on $\{0,1,\cdots, |T|-1\}$, namely, $s_T(i) := j$ if $s$ maps the $i$-th greatest element in $T$ to the $j$-th greatest element in $s(T)$. For example, if $s = s_0 \in \Sym_3$, then $s_{\{0,1\}} = s_0 \in \Sym_2, \ s_{\{0,2\}} = Id \in \Sym_2$.

\begin{figure}
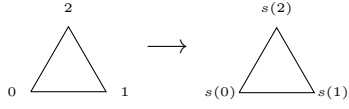

\centering
$\Triangle{0}{0}{0}{1}{2}{}{}{}{} \longrightarrow  \Triangle{0}{0}{s(0)}{s(1)}{s(2)}{}{}{}{} $
\caption{Reordering by a permutation $s$}\label{fig:reodering triangle}
\end{figure}

Let $\T$ be a fixed ordered triangulation with vertices ordered by $0,1,\cdots, N-1$, where $N = |\T^0|$. Let $\T'$ be any ordered triangulation obtained from $\T$ by reordering its vertices. Apparently, each reordering of vertices corresponds to an element of $\Sym_N$. To make it clear, if $s \in \Sym_N$, then we obtain $\T'$ from $\T$ by replacing the label $i$ of a vertex by $s(i)$. See Figure \ref{fig:reodering triangle} for the reordering of a $2$-simplex. {\it That is, we are thinking that the \lq physical' unordered triangulation is fixed all the time; only the labels on vertices are permuted.} For a $k$-simplex $\sigma = (a_0a_1\cdots a_{k-1}) \in \T^k$, denote by $\sigma_s$ the corresponding $k$-simplex $(\sigma(a_0)\sigma(a_1)\cdots\sigma(a_{k-1})) \in (\T')^k$. Note that our convention for expressing a $k$-simplex is to list its labels in increasing order, thus here $\sigma_s$ should really be $(\sigma(a_0)\sigma(a_1)\cdots\sigma(a_{k-1}))$ rearranged in increasing order. For example, if $s(0)=3, s(1) = 2, s(2) = 4$, then the $2$-simplex $(012) \in \T^2$ corresponds to $(012)_s = (234) \in (\T'^2)$. Also note that $\sigma$ and $\sigma_s$ are the same as an unordered simplex. Since every permutation is a composition of $s_i\,'$s, it suffices to prove the partition function is invariant under every $s_i$, the swap of two labels $i$ and $i+1$. If the vertices labeled by $i$ and $i+1$ are not within any $4$-simplex, then after swapping $i$ and $i+1$, the relative order of vertices of any $4$-simplex remains the same, and thus the partition function also remains the same. Hence, we only need to consider the case where the vertices labeled by $i$ and $i+1$ are within some $4$-simplex. (They could be contained in more than one $4$-simplex.)

Here is the basic idea of the proof. Let $s = s_i \in \Sym_{N}$ be some swap, and let $\T'$ be the ordered triangulation obtained from $\T$ by $s$. If $\sigma$ is any $k$-simplex, $s_{\sigma}$ means the restriction of $s$ on the set of vertices of $\sigma$. Note that here we have defined both $s_{\sigma}$ and $\sigma_{s}$. We will show there is a one-to-one correspondence between colorings $F = (g,f)$ of $\T$ and colorings $F' = (g',f')$ of $\T'$ such that,
 \begin{enumerate}
   \item for each $2$-simplex $\beta \in \T^2$, $d_{f(\beta)} = d_{f'(\beta_{s})}$.
   \item for each $3$-simplex $\tau \in \T^3$, there is an isomorphism,
   \begin{align*}
   \begin{array}{cccc}
   \Phi_{s, \tau}\ \colon &V_{F}^{\pm}(\tau) &\overset{\simeq}{\longrightarrow} & V_{F'}^{\pm \epsilon(s_{\tau})}(\tau_s),
   \end{array}
   \end{align*}
   where $\epsilon(s_{\tau})$ is the sign of $\epsilon(s_{\tau})$ taking values in $\{+,-\}$ and we use the convention that $++ = -- = +,\  +- = -+ = -$. The isomorphism
   \begin{align*}
   \begin{array}{cccc}
   \Phi_{s,\tau} \otimes \Phi_{s,\tau}\ \colon &V_{F}^{-}(\tau) \otimes V_{F}^{+}(\tau) &\longrightarrow & V_{F'}^{-}(\tau_{s}) \otimes V_{Fs}^{+}(\tau_s)\ \left(\text{or } V_{F'}^{+}(\tau_{s}) \otimes V_{Fs}^{-}(\tau_s)\right)
   \end{array}
   \end{align*}
   maps $\phi_{F, \tau}$ to $\phi_{F',\tau_{s}}$ (up to a permutation of the tensor components). Let $\Phi_{s} = \bigotimes\limits_{\tau \in \T^3} (\Phi_{s,\tau} \otimes \Phi_{s,\tau})$. Then $\Phi_{s}$ maps $\phi_{F}$ to $\phi_{F'}$.

   \item for each $4$-simplex $\sigma = (01234)$, let
   \begin{align*}
    V_F^{+}(\partial \sigma) &= V^{+}_{F}(0234) \otimes V^{+}_{F}(0124) \otimes V^{-}_{F}(1234) \otimes V^{-}_{F}(0134) \otimes V^{-}_{F}(0123) \\
    V_F^{-}(\partial \sigma) &= V^{-}_{F}(0234) \otimes V^{-}_{F}(0124) \otimes V^{+}_{F}(1234) \otimes V^{+}_{F}(0134) \otimes V^{+}_{F}(0123).
    \end{align*}
    Then $\tilde{Z}_F^{\epsilon(\sigma)}(\sigma)$ is a functional on $V_F^{\epsilon(\sigma)}(\partial \sigma)$. The map
    $$\Phi_{s,\sigma}:= \Phi_{s, 0234} \otimes \Phi_{s, 0124} \otimes \Phi_{s, 1234}\otimes \Phi_{s, 0134}\otimes \Phi_{s, 0123}$$ is an isomorphism from $V_{F}^{\epsilon(\sigma)}(\partial \sigma)$ to $V_{F'}^{\epsilon(\sigma_s)}(\partial \sigma_s),$
    and moreover $\tilde{Z}_F^{\epsilon(\sigma)}(\sigma) = \tilde{Z}_{F'}^{\epsilon(\sigma_s)}(\sigma_s) \circ \Phi_{s, \sigma}$.
 \end{enumerate}

 \noindent The above conditions are sufficient to show $Z_{\cross{\Cat}{G}}(M;\T) = Z_{\cross{\Cat}{G}}(M;\T')$, as illustrated in Theorem \ref{thm:reordering}. Now we build the correspondence between $F$ and $F'$. Given a coloring $F = (g,f)$, define $F' = (g',f')$ as follows.

For each $1$-simplex $\alpha \in \T^1$, define $g'(\alpha_{s}) = g(\alpha)$ if $s_{\alpha}$ is the identity permutation and $g'(\alpha_{s}) = g(\alpha)^{-1}$ otherwise. (In the latter case, $s_{\alpha}$ is the swap permutation in $\Sym_2$.)

For each $2$-simplex $\beta = (012) \in \T^2$, there are three possibilities for $s_{\beta}$, namely, $s_{\beta} = Id, \, s_0, \, s_1$. Define
\begin{equation}
f'(\beta_s) =
\begin{cases}
f(\beta) & s_{\beta} = Id \\
f(\beta)^* & s_{\beta} = s_0 \\
\lsup{g(12)}{f(\beta)^*} & s_{\beta} = s_1.
\end{cases}
\end{equation}

\begin{figure}
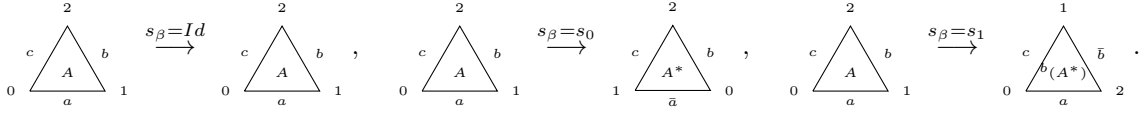

\centering
$
\Triangle{0}{0}{0}{1}{2}{a}{b}{c}{A} \overset{s_{\beta} = Id}{\longrightarrow} \Triangle{0}{0}{0}{1}{2}{a}{b}{c}{A}, \quad \Triangle{0}{0}{0}{1}{2}{a}{b}{c}{A} \overset{s_{\beta} = s_0}{\longrightarrow} \Triangle{0}{0}{1}{0}{2}{\bar{a}}{b}{c}{A^{*}} , \quad \Triangle{0}{0}{0}{1}{2}{a}{b}{c}{A}  \overset{s_{\beta} = s_1}{\longrightarrow} \Triangle{0}{0}{0}{2}{1}{a}{\bar{b}}{c}{\lsup{b}{(A^{*})}}.
$
\caption{The change of coloring on a $2$-simplex}\label{fig:2simplexcolorchange}
\end{figure}

Graphically, this is illustrated in Figure \ref{fig:2simplexcolorchange}. One can check that $F' = (g',f')$ satisfies the constraint of a coloring. The correspondence is apparently one-to-one. Also, it is clear that $d_{f(\beta)} = d_{f'(\beta_s)}$.

For each $3$-simplex $\tau = (0123)$, $s_{\tau}$ could be $Id, s_0, s_1$ or $s_2$. If $s_{\tau} = Id$, for every boundary $2$-simplex $(ijk)$ of $\tau$, $s_{(ijk)} = Id$, thus $(ijk)$ and $(ijk)_{s}$ have the same color and $V_{F}^{\pm}(\tau)  = V_{F'}^{\pm}(\tau_s)$. In this case, $\epsilon(s_{\tau}) = +$, and define
\begin{align*}
\begin{array}{cccc}
   \Phi_{Id, \tau}\ = \ Id\ \colon &V_{F}^{\pm}(\tau) &\overset{\simeq}{\longrightarrow} &V_{F'}^{\pm \epsilon(s_{\tau})}(\tau_{Id}).
   \end{array}
\end{align*}

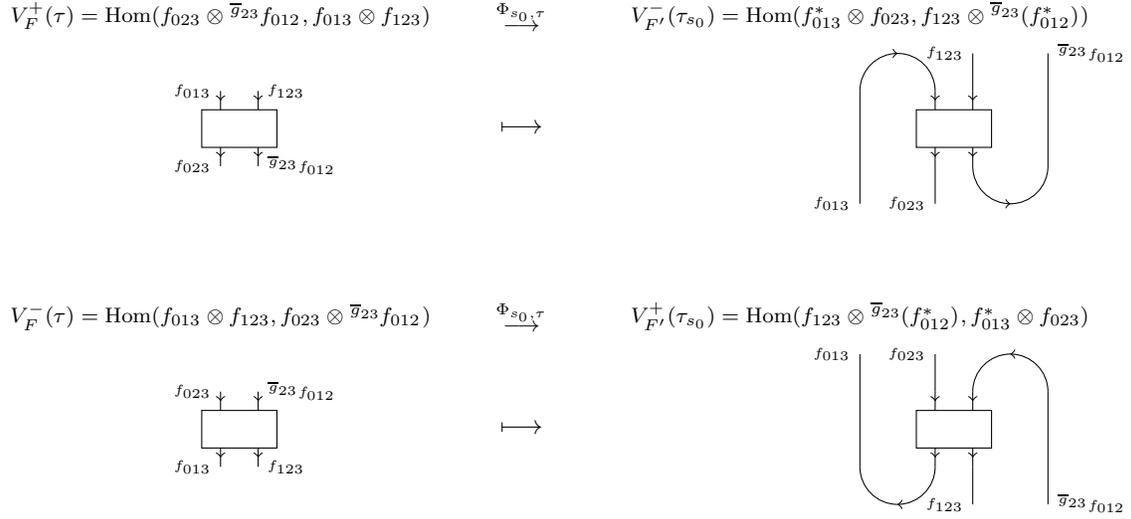
\begin{figure}
 \centering
\begin{tikzpicture}[scale = 0.5]
\begin{scope}[yshift = 8cm]
 \begin{scope}[decoration={
    markings,
    mark=at position 0.5 with {\arrow{>}}}]
   \RecMor{0}{0}{f_{023}}{\lsup{\overline{g}_{23}}{f_{012}}}{f_{013}}{f_{123}}
  \end{scope}
  \begin{scope}[xshift = 8cm]
  \draw (0,1) node{$\longmapsto$};
  \end{scope}
  \begin{scope}[xshift = 17cm, decoration={
    markings,
    mark=at position 0.5 with {\arrow{>}}}]
   \RecMor{2}{0}{}{}{}{}
   \draw [postaction={decorate}] (3,0) arc(-180:0:1cm);
   \draw [postaction={decorate}] (0,2) arc(180:0:1cm);
   \draw (0,-1)node[left]{{\tiny $f_{013}$}} -- (0,2);
   \draw (5,0) -- (5,3)node[right]{{\tiny $\lsup{\overline{g}_{23}}{f_{012}}$}};
   \draw (2,-1)node[left]{{\tiny $f_{023}$}} -- (2,0);
   \draw (3,2) --(3,3)node[left]{{\tiny $f_{123}$}};
  \end{scope}
  \begin{scope}[yshift = 4cm]
  \draw (0,0) node{{\footnotesize $V_{F}^{+}(\tau) = \Hom(f_{023} \otimes \lsup{\overline{g}_{23}}{f_{012}}, f_{013} \otimes f_{123})$}};
  \draw (8,0) node{{\footnotesize $\overset{\Phi_{s_0,\tau}}{\longrightarrow}$}};
  \draw (17,0) node{{\footnotesize $V_{F'}^{-}(\tau_{s_0}) = \Hom( f_{013}^* \otimes f_{023}, f_{123} \otimes \lsup{\overline{g}_{23}}{(f_{012}^*)})$}};
  \end{scope}
\end{scope}

 \begin{scope}[decoration={
    markings,
    mark=at position 0.5 with {\arrow{>}}}]
   \RecMor{0}{0}{f_{013}}{f_{123}}{f_{023}}{\lsup{\overline{g}_{23}}{f_{012}}}
  \end{scope}
  \begin{scope}[xshift = 8cm]
  \draw (0,1) node{$\longmapsto$};
  \end{scope}
  \begin{scope}[xshift = 17cm, decoration={
    markings,
    mark=at position 0.5 with {\arrow{>}}}]
   \RecMor{2}{0}{}{}{}{}
   \draw [postaction={decorate}] (2,0) arc(0:-180:1cm);
   \draw [postaction={decorate}] (5,2) arc(0:180:1cm);
   \draw (0,0) -- (0,3)node[left]{{\tiny $f_{013}$}};
   \draw (5,-1)node[right]{{\tiny $\lsup{\overline{g}_{23}}{f_{012}}$}} -- (5,2);
   \draw (2,2) -- (2,3)node[left]{{\tiny $f_{023}$}};
   \draw (3,-1)node[left]{{\tiny $f_{123}$}} --(3,0);
  \end{scope}
  \begin{scope}[yshift = 4cm]
  \draw (0,0) node{{\footnotesize $V_{F}^{-}(\tau) = \Hom(f_{013} \otimes f_{123}, f_{023} \otimes \lsup{\overline{g}_{23}}{f_{012}})$}};
  \draw (8,0) node{{\footnotesize $\overset{\Phi_{s_0,\tau}}{\longrightarrow}$}};
  \draw (17,0) node{{\footnotesize $V_{F'}^{+}(\tau_{s_0}) = \Hom(f_{123} \otimes \lsup{\overline{g}_{23}}{(f_{012}^*)},  f_{013}^* \otimes f_{023})$}};
  \end{scope}
\end{tikzpicture}
\caption{Definition of $\Phi_{s_0,\tau}$} \label{fig:V(s0)}
\end{figure}

If $s_{\tau} = s_0$, then $s_{(023)} = s_{(123)} = Id, \, s_{(012)} = s_{(013)} = s_0, \, s_{(23)} = Id.$ Note that by definition, the correspondence between $2$-simplices of $\tau$ and those of $\tau_s$ is given by: $(023)_s = (123), (123)_s = (023), (012)_s = (012), (013)_s = (013)$. Hence in $\tau_s$, the coloring of each $2$-simplex is given as follows: $f'(123) = f(023), f'(023) = f(123), f'(012) = f(012)^*, f'(013) = f(013)^*$. Also we have $g'(23) = g(23)$. Hence,
\begin{alignat*}{3}
V_{F'}^{+}(\tau_{s})\quad &&=\quad\Hom(f'_{023} \otimes \lsuprsub{\overline{g'}_{23}}{f'}{012}, f'_{013} \otimes f'_{123}) \quad &&=\quad \Hom(f_{123} \otimes \lsup{\overline{g}_{23}}{(f_{012}^*)}, f_{013}^* \otimes f_{023}) \\
V_{F'}^{-}(\tau_{s}) \quad &&=\quad\Hom(f'_{013} \otimes f'_{123}, f'_{023} \otimes \lsuprsub{\overline{g'}_{23}}{f'}{012}) \quad &&=\quad \Hom( f_{013}^* \otimes f_{023}, f_{123} \otimes \lsup{\overline{g}_{23}}{(f_{012}^*)})
\end{alignat*}
Define the isomorphism $\Phi_{s_0, \tau}\colon V_{F}^{\pm}(\tau) \longrightarrow V_{F'}^{\mp}(\tau_{s_0})$ by graph diagrams as shown in Figure \ref{fig:V(s0)}.

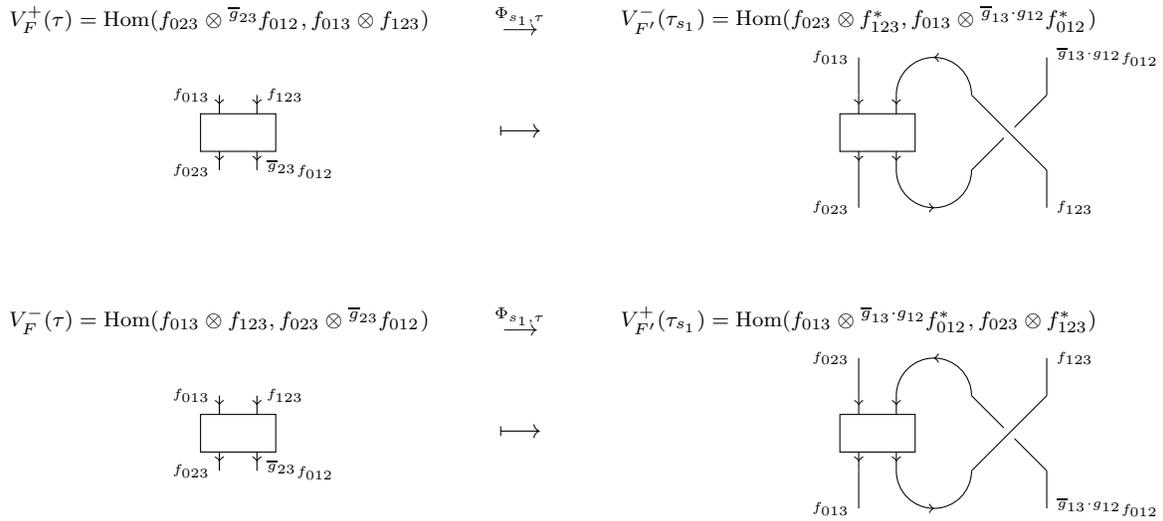
\begin{figure}
 \centering
\begin{tikzpicture}[scale = 0.5]
\begin{scope}[yshift = 8cm]
 \begin{scope}[decoration={
    markings,
    mark=at position 0.5 with {\arrow{>}}}]
   \RecMor{0}{0}{f_{023}}{\lsup{\overline{g}_{23}}{f_{012}}}{f_{013}}{f_{123}}
  \end{scope}
  \begin{scope}[xshift = 8cm]
  \draw (0,1) node{$\longmapsto$};
  \end{scope}
  \begin{scope}[xshift = 15cm, decoration={
    markings,
    mark=at position 0.5 with {\arrow{>}}}]
   \RecMor{2}{0}{}{}{}{}
   \draw [postaction={decorate}] (3,0) arc(-180:0:1cm);
   \draw [postaction={decorate}] (5,2) arc(0:180:1cm);
   \InvBraid{5}{0}
   \draw (2,-1)node[left]{{\tiny $f_{023}$}} -- (2,0);
   \draw (7,2) -- (7,3)node[right]{{\tiny $\lsup{\overline{g}_{13}\cdot g_{12}}{f_{012}}$}};
   \draw (2,3)node[left]{{\tiny $f_{013}$}} -- (2,2);
   \draw (7,0) --(7,-1)node[right]{{\tiny $f_{123}$}};
  \end{scope}
  \begin{scope}[yshift = 4cm]
  \draw (0,0) node{{\footnotesize $V_{F}^{+}(\tau) = \Hom(f_{023} \otimes \lsup{\overline{g}_{23}}{f_{012}}, f_{013} \otimes f_{123})$}};
  \draw (8,0) node{{\footnotesize $\overset{\Phi_{s_1,\tau}}{\longrightarrow}$}};
  \draw (17,0) node{{\footnotesize $V_{F'}^{-}(\tau_{s_1}) = \Hom(f_{023} \otimes f_{123}^*, f_{013} \otimes \lsup{\overline{g}_{13}\cdot g_{12}}{f_{012}^*})$}};
  \end{scope}
\end{scope}

 \begin{scope}[decoration={
    markings,
    mark=at position 0.5 with {\arrow{>}}}]
   \RecMor{0}{0}{f_{023}}{\lsup{\overline{g}_{23}}{f_{012}}}{f_{013}}{f_{123}}
  \end{scope}
  \begin{scope}[xshift = 8cm]
  \draw (0,1) node{$\longmapsto$};
  \end{scope}
  \begin{scope}[xshift = 15cm, decoration={
    markings,
    mark=at position 0.5 with {\arrow{>}}}]
   \RecMor{2}{0}{}{}{}{}
   \draw [postaction={decorate}] (3,0) arc(-180:0:1cm);
   \draw [postaction={decorate}] (5,2) arc(0:180:1cm);
   \Braid{5}{0}
   \draw (2,-1)node[left]{{\tiny $f_{013}$}} -- (2,0);
   \draw (7,0) -- (7,-1)node[right]{{\tiny $\lsup{\overline{g}_{13}\cdot g_{12}}{f_{012}}$}};
   \draw (2,3)node[left]{{\tiny $f_{023}$}} -- (2,2);
   \draw (7,2) --(7,3)node[right]{{\tiny $f_{123}$}};
  \end{scope}
  \begin{scope}[yshift = 4cm]
  \draw (0,0) node{{\footnotesize $V_{F}^{-}(\tau) = \Hom(f_{013} \otimes f_{123}, f_{023} \otimes \lsup{\overline{g}_{23}}{f_{012}})$}};
  \draw (8,0) node{{\footnotesize $\overset{\Phi_{s_1,\tau}}{\longrightarrow}$}};
  \draw (17,0) node{{\footnotesize $V_{F'}^{+}(\tau_{s_1}) = \Hom(f_{013} \otimes \lsup{\overline{g}_{13}\cdot g_{12}}{f_{012}^*},f_{023} \otimes f_{123}^*)$}};
  \end{scope}
\end{tikzpicture}
\caption{Definition of $\Phi_{s_1,\tau}$} \label{fig:V(s1)}
\end{figure}

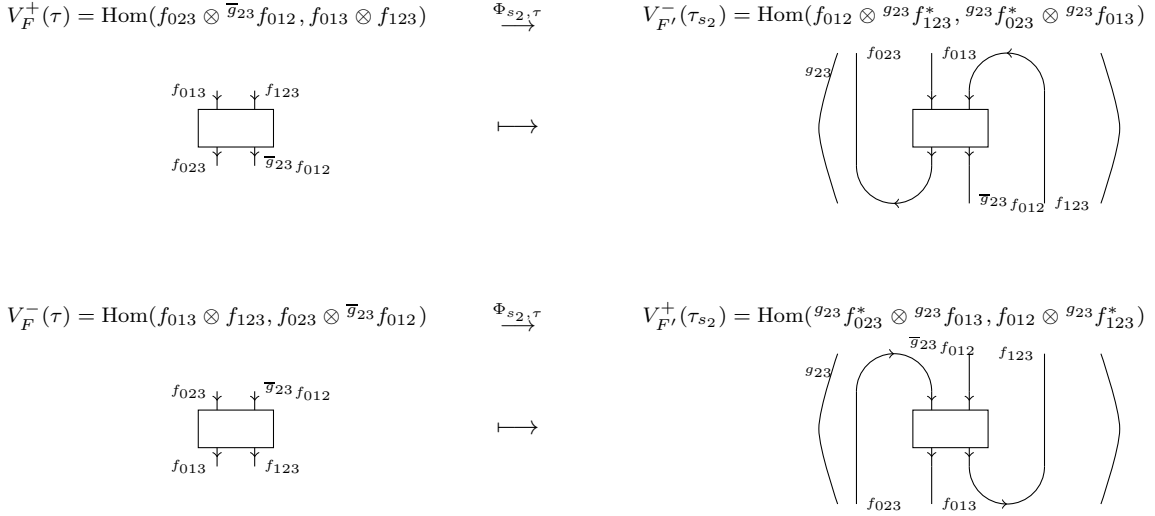
\begin{figure}
 \centering
\begin{tikzpicture}[scale = 0.5]
\begin{scope}[yshift = 8cm]
 \begin{scope}[decoration={
    markings,
    mark=at position 0.5 with {\arrow{>}}}]
   \RecMor{0}{0}{f_{023}}{\lsup{\overline{g}_{23}}{f_{012}}}{f_{013}}{f_{123}}
  \end{scope}
  \begin{scope}[xshift = 8cm]
  \draw (0,1) node{$\longmapsto$};
  \end{scope}
 \begin{scope}[xshift = 17cm, decoration={
    markings,
    mark=at position 0.5 with {\arrow{>}}}]
   \RecMor{2}{0}{}{}{}{}
   \draw [postaction={decorate}] (2,0) arc(0:-180:1cm);
   \draw [postaction={decorate}] (5,2) arc(0:180:1cm);
   \draw (0,0) -- (0,3)node[right]{{\tiny $f_{023}$}};
   \draw (5,-1)node[right]{{\tiny $f_{123}$}} -- (5,2);
   \draw (2,2) -- (2,3)node[right]{{\tiny $f_{013}$}};
   \draw (3,-1)node[below]{{\tiny $\lsup{\overline{g}_{23}}{f_{012}}$}} --(3,0);

   \draw plot[smooth] coordinates{(-0.5,-1)(-1,1)(-0.5,3)};
   \draw plot[smooth] coordinates{(6.5,-1)(7,1)(6.5,3)};
   \draw (-1.2,2.5) node{{\tiny $g_{23}$}};
  \end{scope}
  \begin{scope}[yshift = 4cm]
  \draw (0,0) node{{\footnotesize $V_{F}^{+}(\tau) = \Hom(f_{023} \otimes \lsup{\overline{g}_{23}}{f_{012}}, f_{013} \otimes f_{123})$}};
  \draw (8,0) node{{\footnotesize $\overset{\Phi_{s_2,\tau}}{\longrightarrow}$}};
  \draw (18,0) node{{\footnotesize $V_{F'}^{-}(\tau_{s_2}) = \Hom(f_{012} \otimes \lsup{g_{23}}{f_{123}^*}, \lsup{g_{23}}{f_{023}^*} \otimes\lsup{g_{23}}{f_{013}})$}};
  \end{scope}
\end{scope}

 \begin{scope}[decoration={
    markings,
    mark=at position 0.5 with {\arrow{>}}}]
   \RecMor{0}{0}{f_{013}}{f_{123}}{f_{023}}{\lsup{\overline{g}_{23}}{f_{012}}}
  \end{scope}
  \begin{scope}[xshift = 8cm]
  \draw (0,1) node{$\longmapsto$};
  \end{scope}
  \begin{scope}[xshift = 17cm, decoration={
    markings,
    mark=at position 0.5 with {\arrow{>}}}]
   \RecMor{2}{0}{}{}{}{}
   \draw [postaction={decorate}] (3,0) arc(-180:0:1cm);
   \draw [postaction={decorate}] (0,2) arc(180:0:1cm);
   \draw (0,-1)node[right]{{\tiny $f_{023}$}} -- (0,2);
   \draw (5,0) -- (5,3)node[left]{{\tiny $f_{123}$}};
   \draw (2,-1)node[right]{{\tiny $f_{013}$}} -- (2,0);
   \draw (3,2) --(3,3);
   \draw (2.5,3.3)node{{\tiny $\lsup{\overline{g}_{23}}{f_{012}}$}};

    \draw plot[smooth] coordinates{(-0.5,-1)(-1,1)(-0.5,3)};
   \draw plot[smooth] coordinates{(6.5,-1)(7,1)(6.5,3)};
   \draw (-1.2,2.5) node{{\tiny $g_{23}$}};
  \end{scope}
  \begin{scope}[yshift = 4cm]
  \draw (0,0) node{{\footnotesize $V_{F}^{-}(\tau) = \Hom(f_{013} \otimes f_{123}, f_{023} \otimes \lsup{\overline{g}_{23}}{f_{012}})$}};
  \draw (8,0) node{{\footnotesize $\overset{\Phi_{s_2,\tau}}{\longrightarrow}$}};
  \draw (18,0) node{{\footnotesize $V_{F'}^{+}(\tau_{s_2}) = \Hom(\lsup{g_{23}}{f_{023}^*} \otimes\lsup{g_{23}}{f_{013}}, f_{012} \otimes \lsup{g_{23}}{f_{123}^*})$}};
  \end{scope}
\end{tikzpicture}
\caption{Definition of $\Phi_{s_2,\tau}$} \label{fig:V(s2)}
\end{figure}

Similarly, if $s = s_1$ or $s_2$, we leave it as an exercise to work out the details of the colorings of $\tau_s$. The isomorphisms $\Phi_{s_1, \tau}$ and $\Phi_{s_2, \tau}$ are defined as shown in Figure \ref{fig:V(s1)} and \ref{fig:V(s2)}, respectively. The following lemma is straightforward.

\begin{lemma}
\label{lem:preservepairing}
The isomorphism $\Phi_{s, \tau} \otimes \Phi_{s, \tau}$ preserves the pairing on $V_{F}^{-}(\tau) \otimes V_{F}^{+}(\tau)$. Consequently, $\Phi_{s, \tau} \otimes \Phi_{s, \tau}$ maps $\phi_{F, \tau}$ to $\phi_{F', \tau_s}$.
\end{lemma}

Finally we consider the reordering of $4$-simplices.  Let $\sigma = (01234)$ be any $4$-simplex of $\T$, then $s_{\sigma} = Id,\, s_0, \, s_1, \,s_2$, or $s_3$. Let
\begin{align*}
\begin{array}{cccc}
\Phi_{s,\sigma} \ =\ \Phi_{s, 0234} \otimes \Phi_{s, 0124} \otimes \Phi_{s, 1234}\otimes \Phi_{s, 0134}\otimes \Phi_{s, 0123}\ \colon & V_{F}^{\epsilon(\sigma)}(\sigma)& \longrightarrow & V_{F'}^{\epsilon(\sigma_s)}(\sigma_{s}).
\end{array}
\end{align*}
If $s_{\sigma} = Id$, it is clear that for any $3$-simplex $\tau$ of $\sigma$ we have $V_{F}^{\pm}(\tau) = V_{F'}^{\pm}(\tau_{s})$, $s_{\tau} = Id$ and thus $\Phi_{s, \tau} = Id$ implying $\Phi_{s, \sigma} = Id$. In this case, the equality $\tilde{Z}_F^{\epsilon(\sigma)}(\sigma) = \tilde{Z}_{F'}^{\epsilon(\sigma_s)}(\sigma_s) \circ \Phi_{s, \sigma}$ automatically holds. The following lemma asserts that in other cases it also holds.

\begin{lemma}
\label{lem:ZtildePhi}
Let $\sigma = (01234)$ be any $4$-simplex of $\T$, then $\tilde{Z}_F^{\epsilon(\sigma)}(\sigma) = \tilde{Z}_{F'}^{\epsilon(\sigma_s)}(\sigma_s) \circ \Phi_{s, \sigma}$.
\begin{proof}
If $s_{\sigma} = Id$, the equality is shown above. If $s_{\sigma} = s_0, \, s_1, \,s_2$, or $s_3$, then $\epsilon(\sigma) = -\epsilon(\sigma_s)$. We need to show $\tilde{Z}_F^{\pm}(\sigma) = \tilde{Z}_{F'}^{\mp}(\sigma_s) \circ \Phi_{s, \sigma}. $

We show the proof for the case $s_{\sigma} = s_0$. Other cases can be done similarly. If $s= s_0$, then the correspondence between $3$-simplices of $\sigma$ and those of $\sigma_s$ are given by
\begin{align*}
(0234)_s = (1234), \, (0124)_s = (0124), \, (1234)_s = (0234), \, (0134)_s = (0134), \, (0123)_s = (0123).
\end{align*}
And the restriction of $s$ on each $3$-simplex is given by
\begin{align*}
s_{(0234)} = Id, \, s_{(0124)} = s_0, \, s_{(1234)} = Id, \, s_{(0134)} = s_0, \, s_{(0123)} = s_0.
\end{align*}
Thus, $\Phi_{s, \sigma}$ can be expressed as
\begin{multline*}
V_F^{+}(0234) \otimes V_F^{+}(0124)\otimes V_F^{-}(1234)\otimes V_F^{-}(0134)\otimes V_F^{-}(0123) \\ \longrightarrow V_{F'}^{+}(1234) \otimes V_{F'}^{-}(0124)\otimes V_{F'}^{-}(0234)\otimes V_{F'}^{+}(0134)\otimes V_{F'}^{+}(0123),
\end{multline*}
\noindent mapping $\psi_{0} \otimes\psi_{1} \otimes\psi_{2} \otimes\psi_{3} \otimes\psi_{4} $ to $ \psi_{0} \otimes \Phi_{s_0, (0124)}(\psi_1) \otimes \psi_{2} \otimes \Phi_{s_0, (0134)}(\psi_3) \otimes \Phi_{s_0, (0123)}(\psi_4).$

Then identity $\tilde{Z}_F^{+}(\sigma) = \tilde{Z}_{F'}^{-}(\sigma_s) \circ \Phi_{s, \sigma}$ is obtained by inserting the $\psi_i$ into the corresponding rectangle in Figure \ref{fig:25jbox} (Left), inserting the image of $\psi_i$ under the above map into the corresponding rectangle in Figure \ref{fig:25jbox} (Right), and showing the resulting two graph diagrams are isotopic.
\end{proof}
\end{lemma}

Now we are ready to prove the main theorem of the subsection on the independence of the partition function on the ordering of vertices.
\begin{theorem}
\label{thm:reordering}
Let $\T, \T'$ be two ordered triangulations of $M$ which are the same when the ordering is ignored, then $Z_{\cross{\Cat}{G}}(M;\T) = Z_{\cross{\Cat}{G}}(M;\T')$.
\begin{proof}
As analyzed at the beginning of this subsection, it suffices to consider the case where $\T'$ is obtained from $\T$ by some $s = s_i \in \Sym_{N}$, $N = |\T^0|$. Given any coloring $F = (g,f)$ of $\T$, let $F' = (g',f')$ be the corresponding coloring of $\T'$. Then for any $2$-simplex $\beta \in \T$, we have $d_{f(\beta)} = d_{f'(\beta_s)}$. Let
\begin{align*}
\begin{array}{cccc}
\Phi_{s}\ =\ \bigotimes\limits_{\tau \in \T^3} \Phi_{s, \tau} \otimes \Phi_{s, \tau}\  \colon & V_{F} \left(= \bigotimes\limits_{\tau \in \T^3} V_{F}^{-}(\tau)\otimes V_{F}^{+}(\tau)\right) &\overset{\simeq}{\longrightarrow} &V_{F'}.
\end{array}
\end{align*}
By Lemma \ref{lem:preservepairing}, $(\Phi_{s, \tau} \otimes \Phi_{s, \tau})(\phi_{F, \tau}) = \phi_{F', \tau_s}$. Therefore, $\Phi_s(\phi_F) = \phi_{F'}$.

By Lemma \ref{lem:ZtildePhi}, for any $4$-simplex $\sigma$, $\tilde{Z}_F^{\epsilon(\sigma)}(\sigma) = \tilde{Z}_{F'}^{\epsilon(\sigma_s)}(\sigma_s) \circ \Phi_{s, \sigma}$. Thus,
\begin{align*}
\bigotimes\limits_{\sigma \in \T^4} \tilde{Z}_F^{\epsilon(\sigma)}(\sigma) &= \bigotimes\limits_{\sigma_s \in (\T')^4} \tilde{Z}_{F'}^{\epsilon(\sigma_s)}(\sigma_s) \circ \bigotimes\limits_{\sigma\in \T^4} \Phi_{s, \sigma}.
\end{align*}

By virtual of the definition, $\Phi_s = \bigotimes\limits_{\sigma\in \T^4} \Phi_{s, \sigma}$. Therefore,
\begin{align*}
\prod\limits_{\beta \in \T^2}d_{f(\beta)}\,\left(\bigotimes\limits_{\sigma \in \T^4}\tilde{Z}^{\epsilon(\sigma)}_{F}(\sigma)\right)(\phi_{F}) &= \prod\limits_{\beta_s \in (\T')^2}d_{f'(\beta_s)}\,\left(\bigotimes\limits_{\sigma_s \in (\T')^4} \tilde{Z}_{F'}^{\epsilon(\sigma_s)}(\sigma_s) \circ \Phi_s\right)(\phi_{F}) \\
  &= \prod\limits_{\beta_s \in (\T')^2}d_{f'(\beta_s)}\,\left(\bigotimes\limits_{\sigma_s \in (\T')^4} \tilde{Z}_{F'}^{\epsilon(\sigma_s)}(\sigma_s) \right)(\phi_{F'})
\end{align*}

\end{proof}
\end{theorem}

\subsection{Pachner Moves}
\label{subsec:pachner}
In this subsection, we work in the category of piecewise linear (PL) manifolds. A PL manifold is a topological manifold endowed with an equivalence class of PL structures and a triangulation of a PL manifold is a simplicial complex that corresponds to a choice of a representative from the equivalence class of PL structures. Any two triangulations of the same manifold are related by a sequence of Pachner moves. See \cite{rourke2012introduction} for detailed discussions. 

Let $\sigma_{n+1}$ be any $(n+1)$-simplex. Its boundary, denoted by $\partial(\sigma_{n+1})$, consists of $n+2$ $n$-simplices. Let $\partial(\sigma_{n+1}) = I \sqcup J$ be a bi-partition of $\partial(\sigma_{n+1})$ such that $|I| = k$ and $|J| = n+2-k$. Two triangulations $\T $ and $ \T'$ are related by a Pachner move of type $(k,n+2-k)$ if $I \subset \T$ and $\T' = (\T \setminus I) \cup J$ for some $I, J$. Apparently, Pachner moves of types $(k,n+2-k)$ and $(n+2-k,k)$ are inverse operations if the partitions are taken to be $(I,J)$ and $(J,I)$ respectively. 

\begin{theorem}\cite{PACHNER1991129}
Let $M = M^n$ be a closed PL manifold and $\T,\ \T'$ be two triangulations of $M$, then there is a sequence of triangulations of $M$, $\T = \T_0,\ \T_1,\ \cdots,\ \T_k = \T'$, such that any two neighboring triangulations in the sequence are related by a Pachner move of type $(k, n+2-k )$, for some $1 \leq k \leq \frac{n}{2}+1$, or its inverse.
\end{theorem}

We remark that for $n=4$, the category of PL manifolds is equivalent to the category of smooth manifolds. Namely, each PL manifold admits a unique smooth structure, and vice versa. 

Some examples of Pachner moves are listed here.

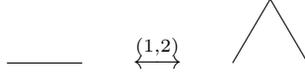
\begin{figure}
\centering
 \begin{tikzpicture}[scale = 0.5]
  \begin{scope}[xshift = -2cm]
  \draw (2,0) -- (4,0);
  \draw (6,0.3) node{$\overset{(1,2)}{\longleftrightarrow}$};
  \draw (8,0) -- (9,1.7) -- (10,0);
  \end{scope}
 \end{tikzpicture}
 \caption{Pachner move $(1,2)$; $n=1$}\label{fig:pachner dim=1}
\end{figure}

If $dim(M) = n = 1,$ we have Pachner moves of type $(1,2)$ and $(2,1)$ which are inverse to each other. See Figure \ref{fig:pachner dim=1}.

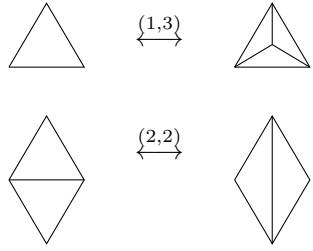
\begin{figure}
\centering
 \begin{tikzpicture}[scale = 0.5]
  \begin{scope}[xshift = 0cm]
   \draw (0,0) -- (1,1.7) -- (2,0) -- (0,0);
  \end{scope}
  \draw (4,1) node{$\overset{(1,3)}{\longleftrightarrow}$};
  \begin{scope}[xshift = 6cm]
   \draw (0,0) -- (1,1.7) -- (2,0) -- (0,0);
   \draw (0,0) -- (1,0.6) -- (1,1.7);
   \draw (2,0) -- (1,0.6);
  \end{scope}

  \begin{scope}[yshift = -3cm]
   \begin{scope}[xshift = 0cm]
    \draw (0,0) -- (1,1.7) -- (2,0) -- (0,0);
    \draw (0,0) -- (1,-1.7) -- (2,0);
   \end{scope}
    \draw (4,1) node{$\overset{(2,2)}{\longleftrightarrow}$};
    \begin{scope}[xshift = 6cm]
     \draw (0,0) -- (1,1.7) -- (2,0) -- (1,-1.7) -- (0,0);
     \draw (1,1.7) -- (1,-1.7);
    \end{scope}
  \end{scope}
 \end{tikzpicture}
 \caption{Pachner moves $(1,3)$, $(2,2)$; $n=2$}\label{fig:pachner dim=2}
\end{figure}

If $n = 2$, we have Pachner moves of type $(1,3)$, $(2,2)$, and $(3,1)$. See Figure \ref{fig:pachner dim=2}.

We are more interested in the case $n=4$, where we have Pachner moves of types $(3,3), (2,4), (1,5)$ and their inverses. Given a $5$-simplex $\sigma_5$, we order its vertices by $0,1,2,3,4$, and denote the face which does not contain the vertex $i$ by $(0 \cdots \hat{i}\cdots 4)$. For each type of Pachner moves, we pick a specific partition $I, J$, and call it the {\it typical} Pachner move. The {\it typical} Pachner moves for $n=4$ are as follows:
\begin{align*}
(02345) (01245) (01234) \qquad &\overset{(3,3)}{\longleftrightarrow} \qquad (12345) (01345) (01235) \\
(02345) (01245) (01234) (12345) \qquad &\overset{(2,4)}{\longleftrightarrow} \qquad  (01345) (01235) \\
(02345) (01245) (01234) (12345) (01345) \qquad &\overset{(1,5)}{\longleftrightarrow} \qquad   (01235)
\end{align*}

\subsection{Invariance under Pachner Moves}
\label{sec:invariance pachner}
We prove that the partition function $Z_{\cross{\Cat}{G}}(M;\T)$ defined in Section \ref{sec:partition} is invariant under Pachner moves. Since in Section \ref{sec:invariance ordering} it has been shown that $Z_{\cross{\Cat}{G}}(M;\T)$ does not depend on the ordering of the vertices, we only need to consider {\it typical} Pachner moves of types $3$-$3$, $2$-$4$, and $1$-$5$. Let $\T, \T'$ be two ordered triangulations whose vertices are ordered in such a way that they only differ by a {\it typical} Pachner move. In each of the cases, let $I_k = \T^k \setminus (\T')^k$ be the set of $k$-simplices which belong to $\T$ but not $\T'$, and similarly let $I'_k = (\T')^k \setminus \T^k$. Tables \ref{table:comparison 3-3}, \ref{table:comparison 2-4}, and \ref{table:comparison 1-5} list the differences between $\T$ and $\T'$ corresponding to each {\it typical} Pachner move $3$-$3$, $2$-$4$, and $1$-$5$, respectively. Without loss of generality, assume the $4$-simplices $02345$, $01245$, $01234$ in $\T$ are positive, then $12345$, $01345$, $01235$ in $\T'$ are positive, and if any of $12345$, $01345$, $01235$ are transported to $\T$, they become negative.

\begin{table}
\centering
\begin{tabular}{c| ccc|ccc}
      & \multicolumn{3}{|c|}{$I_k$} & \multicolumn{3}{|c}{$I'_k$} \\
\hline
$k=4$ & $02345$ & $01245$ & $01234$  & $12345$ & $01345$ & $01235$ \\
$k=3$ & $0124$  & $0234$  & $0245$   & $0135$  & $1235$  & $1345$  \\
$k=2$ & $024$   &         &          & $135$   &         &         \\
\end{tabular}
\caption{Comparison of $\T$ and $\T'$ for $3$-$3$ move}\label{table:comparison 3-3}
\end{table}

\begin{table}
\centering
\begin{tabular}{c| ccc|cc}
                       & \multicolumn{3}{|c|}{$I_k$} & \multicolumn{2}{|c}{$I'_k$} \\
\hline
\multirow{2}{*}{$k=4$} & $02345$ & $01245$ & $01234$   & $01345$ & $01235$ \\
                       & $12345$ &         &           &         &          \\
\hline
\multirow{2}{*}{$k=3$} & $0124$  & $0234$  & $0245$   & $0135$   &    \\
                       & $2345$  & $1245$  & $1234$   &          &    \\
\hline
\multirow{2}{*}{$k=2$} & $024$   &         &          &          &                \\
                       & $245$   & $234$   & $124$    &          &       \\
\hline
$k=1$                  & $24$    &         &          &          & \\
\end{tabular}
\caption{Comparison of $\T$ and $\T'$ for $2$-$4$ move}\label{table:comparison 2-4}
\end{table}

\begin{table}
\centering
\begin{tabular}{c| cccccc|c}
                       & \multicolumn{6}{|c|}{$I_k$}                            & $I'_k$ \\
\hline
\multirow{2}{*}{$k=4$} & $02345$ & $01245$ & $01234$ & $12345$ &       &       & $01235$  \\
                       & $01345$ &         &         &         &       &       &           \\
\hline
\multirow{2}{*}{$k=3$} & $0124$  & $0234$  & $0245$  & $2345$  & $1245$& $1234$ &       \\
                       & $0134$  & $0145$  & $1345$  & $0345$  &       &        &             \\
\hline
\multirow{2}{*}{$k=2$} & $024$   & $245$   & $234$   & $124$   &       &        &                         \\
                       & $034$   & $014$   & $134$   & $045$   & $145$ & $345$  &            \\
\hline
\multirow{2}{*}{$k=1$} & $24$    &         &         &          &       &        &           \\
                       & $04$    & $14$    & $34$    & $45$     &       &        &            \\
\hline
$k = 0$                & $4$     &         &         &          &       &        &  \\
\end{tabular}
\caption{Comparison of $\T$ and $\T'$ for $1$-$5$ move}\label{table:comparison 1-5}
\end{table}

By a complete set of representatives we mean a set of objects which contains exactly one representative for each isomorphism class of simple objects. For the rest of the section, a summation variable concerning simple objects in the category is always assumed to run through an arbitrary complete set of representatives, unless otherwise stated. Recall that an extended coloring is a map $\hat{F} = (g,f,t)$, which assigns a group element to each $1$-simplex, an isomorphism class of simple objects to each $2$-, and $3$-simplex, such that these assignments satisfy some restrictions, see Definition \ref{def:extendedcoloring}. Then to define the partition function, a representative is arbitrarily chosen for each $2$- and $3$-simplex with an extended coloring. The partition function is a sum over all extended colorings, see Equation \ref{equ:statesumdef}. In light of the observations above, we can rephrase the definition of the partition function as follows.

We let each $1$-simplex run through the set of all elements in $G$, and let each $2$-, $3$-simplex run through an arbitrary complete set of representatives (i.e., each $2$-, $3$-simplex has its own complete set). Then an extended coloring can be viewed as a choice of value for $1$-, $2$-, and $3$-simplices, so that the resulting configuration satisfies the restriction in Definition \ref{def:extendedcoloring}. If a configuration is not an extended coloring, we set its contribution to zero. As before, we denote the color of a simplex by the simplex itself. Then Equation \ref{equ:statesumdef} can be rewritten as Equation \ref{equ:statesumdefrewritten}, where we have omitted the dependence of $\hat{Z}^{\epsilon(\sigma)}(\sigma)$ $\hat{F}$.

\begin{equation}
\label{equ:statesumdefrewritten}
Z_{\cross{\Cat}{G}}(M;\T) = \sum\limits_{\T^1,\T^2,\T^3}\frac{(D^2/|G|)^{|\T^0|}}{(D^2)^{|\T^1|}} \frac{\prod\limits_{\beta \in \T^2}d_{\beta}}{\prod\limits_{\tau \in \T^3}d_{\tau}}\prod\limits_{\sigma \in \T^4}\hat{Z}^{\epsilon(\sigma)}(\sigma)
\end{equation}

By Equation \ref{equ:statesumdefrewritten}, to prove $Z_{\cross{\Cat}{G}}(M;\T) = Z_{\cross{\Cat}{G}}(M;\T')$, it suffices to show that Equations \ref{equ:3-3}, \ref{equ:2-4}, and \ref{equ:1-5} hold, corresponding to the {\it typical} Pachner move $3$-$3$, $2$-$4$, $1$-$5$, respectively, where all simplices which do not occur in the summation are assumed to have a fixed coloring on them.

\begin{align}
\label{equ:3-3}
    & \sum\limits_{\substack{I_2, \\ I_3}} \frac{d_{024}}{d_{0124}d_{0234}d_{0245}} \hat{Z}^{+}(01234)\hat{Z}^{+}(01245)\hat{Z}^{+}(02345) \nonumber \\
    =&
\sum\limits_{\substack{I'_2, \\ I'_3}} \frac{d_{135}}{d_{0135}d_{1235}d_{1345}} \hat{Z}^{+}(01235)\hat{Z}^{+}(01345)\hat{Z}^{+}(12345)
\end{align}

\begin{align}
\label{equ:2-4}
    & \sum\limits_{\substack{I_1, \\ I_2, \\ I_3}}\frac{1}{D^2} \frac{\prod\limits_{\beta \in I_2}d_{\beta}}{\prod\limits_{\tau \in I_3}d_{\tau}} \hat{Z}^{+}(01234)\hat{Z}^{+}(01245)\hat{Z}^{+}(02345)\hat{Z}^{-}(12345) \nonumber \\
    =&
\sum\limits_{\substack{I_3'}} \frac{1}{d_{0135}} \hat{Z}^{+}(01235)\hat{Z}^{+}(01345)
\end{align}

\begin{align}
\label{equ:1-5}
    & \sum\limits_{\substack{I_0, I_1, I_2, I_3}}\frac{(D^2/|G|)^{|I_0|}}{(D^2)^{|I_1|}} \frac{\prod\limits_{\beta \in I_2}d_{\beta}}{\prod\limits_{\tau \in I_3}d_{\tau}} \hat{Z}^{+}(01234)\hat{Z}^{+}(01245)\hat{Z}^{+}(02345)\hat{Z}^{-}(12345)\hat{Z}^{-}(01345) \nonumber \\
    =& \hat{Z}^{+}(01235)
\end{align}

We give the proof of Equation \ref{equ:3-3} below and leave the proof of Equation \ref{equ:2-4} and Equation \ref{equ:1-5} as an exercise. In the proof, the following simple lemmas will be used heavily. We take one more convention that {\it all diagrams of morphisms in figures of this subsection are assumed to have their top and bottom identified}, namely, the diagrams in figures are actually the trace of the morphisms drawn.

\begin{lemma}[Merging Formula]
\label{lem:trace identity}
Let $A, B$ be two objects of $\cross{\Cat}{G}$, and $\{e_{A,B,i}:i \in I\}, \{e_{B,A,j}: j \in I\}$ be a basis of $\Hom(A,B)$ and $\Hom(B,A)$, respectively, such that $\langle e_{A,B,i}, e_{B,A,j} \rangle = \alpha_i \delta_{i,j}$, then for $F \in \Hom(A,B), G \in \Hom(B,A)$, we have $\langle G, F\rangle = \sum\limits_{i \in I} \frac{1}{\alpha_i}\langle G,e_{A,B,i}\rangle \langle F,e_{B,A,i}\rangle$. In particular,

\begin{enumerate}
\item Let $F \in \Hom(a \otimes b, c \otimes d), \ G  \in \Hom(c \otimes d, a \otimes b),$ then
$$\langle G , F\rangle = \sum\limits_{x} \frac{1}{d_{x}} \langle F , B_{cd,ab}^x\rangle \langle G , B_{ab,cd}^x\rangle.$$
\item Let $F \in \Hom(a \otimes b \otimes c, d \otimes e \otimes f), \ G  \in \Hom(d \otimes e \otimes f, a \otimes b \otimes c),$ then
$$\langle G , F\rangle = \sum\limits_{x,y,z} \frac{d_{z}}{d_{x}d_{y}} \langle F , (Id \otimes B_{zf,bc}^{y}) \circ (B_{de,az}^{x} \otimes Id)\rangle  \langle G , (B_{az,de}^{x} \otimes Id) \circ (Id \otimes B_{bc,zf}^{y})\rangle.$$
\end{enumerate}
The identities in Part $1$ and Part $2$ are also illustrated in Figure \ref{fig:trace identity 1} and Figure \ref{fig:trace identity 2}, respectively.
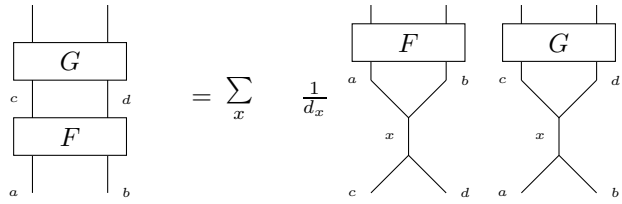
\begin{figure}
\centering
 \begin{tikzpicture}[scale = 0.5]
  \begin{scope}
   \draw (0,0) -- (0,5);
   \draw (2,0) -- (2,5);
   \draw [fill = white] (-0.5,1) rectangle (2.5,2) node[pos=.5] {$F$};
   \draw [fill = white] (-0.5,3) rectangle (2.5,4) node[pos=.5] {$G$};
   \draw (-0.5, 0) node{{\tiny $a$}};
   \draw (2.5, 0) node{{\tiny $b$}};
   \draw (-0.5, 2.5) node{{\tiny $c$}};
   \draw (2.5, 2.5) node{{\tiny $d$}};
  \end{scope}

  \begin{scope}[xshift = 4.5cm]
   \draw (0,2.5)node{$=$};
   \draw (1,2.5) node{$\sum\limits_{{\tiny x}}$};
   \draw (3,2.5) node {$\frac{1}{d_x}$};
  \end{scope}

  \begin{scope}[xshift = 9cm]
   \ShortIshape{1}{1}{c}{d}{a}{b}{x}
   \draw (0,3) -- (0,5);
   \draw (2,3) -- (2,5);
   \draw [fill = white](-0.5,3.5) rectangle (2.5, 4.5) node[pos = .5] {$F$};
  \end{scope}

   \begin{scope}[xshift = 13cm]
   \ShortIshape{1}{1}{a}{b}{c}{d}{x}
   \draw (0,3) -- (0,5);
   \draw (2,3) -- (2,5);
   \draw [fill = white](-0.5,3.5) rectangle (2.5, 4.5) node[pos = .5] {$G$};
  \end{scope}
 \end{tikzpicture}
 \caption{Trace Identity (1)} \label{fig:trace identity 1}
\end{figure}

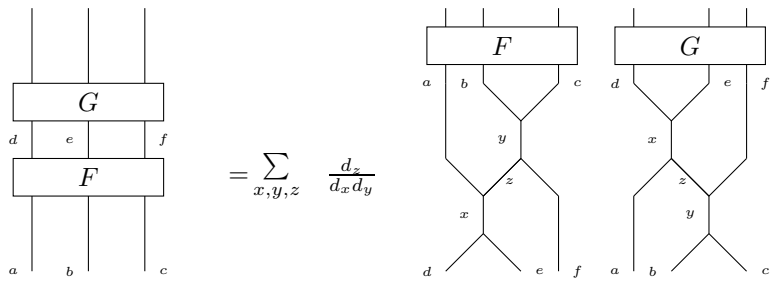
\begin{figure}
\centering
  \begin{tikzpicture} [scale = 0.5]
    \begin{scope}
      \VlineDot{0}{0}{7}
      \VlineDot{1.5}{0}{7}
      \VlineDot{3}{0}{7}
      \draw [fill = white] (-0.5,2) rectangle (3.5,3) node[pos = .5]{$F$};
      \draw [fill = white] (-0.5,4) rectangle (3.5,5) node[pos = .5]{$G$};
      \draw (-0.5, 0) node{{\tiny $a$}};
      \draw (1, 0) node{{\tiny $b$}};
      \draw (3.5, 0) node{{\tiny $c$}};
      \draw (-0.5, 3.5) node{{\tiny $d$}};
      \draw (1, 3.5) node{{\tiny $e$}};
      \draw (3.5, 3.5) node{{\tiny $f$}};
    \end{scope}

    \begin{scope}[xshift = 5.5cm]
     \draw (0,2.5)node{$=$};
     \draw (1,2.5) node{$\sum\limits_{{\tiny x,y,z}}$};
     \draw (3,2.5) node {$\frac{d_z}{d_xd_y}$};
   \end{scope}

   \begin{scope}[xshift = 11cm]
    \ShortIshape{1}{1}{d}{e}{}{}{x}
    \draw (3,0) -- (3,2);
    \ShortIshape{2}{3}{}{}{b}{c}{y}
    \draw (0,3) -- (0,5);
    \draw (0,5) -- (0,7);
    \draw (1,5) -- (1,7);
    \draw (3,5) -- (3,7);
    \draw [fill = white] (-0.5,5.5) rectangle (3.5, 6.5) node[pos = .5]{$F$};
    \draw (-0.5,5)node{{\tiny $a$}};
    \draw (3.5,0)node{{\tiny $f$}};
    \draw (1.7,2.4)node{{\tiny $z$}};


   \end{scope}

   \begin{scope}[xshift = 16cm]
    \ShortIshape{2}{1}{b}{c}{}{}{y}
    \draw (0,0) -- (0,2);
    \ShortIshape{1}{3}{}{}{d}{e}{x}
    \draw (3,3) -- (3,5);
    \draw (0,5) -- (0,7);
    \draw (2,5) -- (2,7);
    \draw (3,5) -- (3,7);
    \draw [fill = white] (-0.5,5.5) rectangle (3.5, 6.5) node[pos = .5]{$G$};
    \draw (-0.5,0)node{{\tiny $a$}};
    \draw (3.5,5)node{{\tiny $f$}};
    \draw (1.3,2.4)node{{\tiny $z$}};


   \end{scope}
  \end{tikzpicture}
  \caption{Trace Identity (2)}\label{fig:trace identity 2}
\end{figure}

%
\end{lemma}

\begin{lemma}
\label{lem:dimension identity}
Let $g,g' \in \text{Gr}(\cross{\Cat}{G}) \subset G$ be fixed and $a, \ b$ be simple objects of $\cross{\Cat}{G}$ such that $a \otimes b \in \Cat_{gg'} $, then the equality in Figure \ref{fig:dimension identity} holds, where $c,\ d, \ e$ are simple objects in the $g$-, $g'$-, $(gg')$-sector, respectively.
\begin{figure}
\centering
 \begin{tikzpicture}[scale = 0.5]
  \draw (0,3) node{$\sum\limits_{\substack{c \in \Cat_g, d \in \Cat_{g'}, \\ e \in \Cat_{gg'}}} \frac{d_c d_d}{d_e}$};
  \begin{scope}[xshift = 2cm]
  \ShortIshape{3}{1}{a}{b}{c}{d}{e}
  \ShortIshape{3}{4}{c}{d}{a}{b}{e}
  \draw (7,3) node{$= \ \frac{D^2}{|\text{Gr}(\cross{\Cat}{G})|}$};
  \draw (9,0) -- (9,6);
  \draw (11,0) -- (11,6);
  \draw (8.5,0) node{\tiny{$a$}};
  \draw (8.5,6) node{\tiny{$a$}};
  \draw (11.5,0) node{\tiny{$b$}};
  \draw (11.5,6) node{\tiny{$b$}};
  \end{scope}
 \end{tikzpicture}
 \caption{Dimension Identity}\label{fig:dimension identity}
\end{figure}
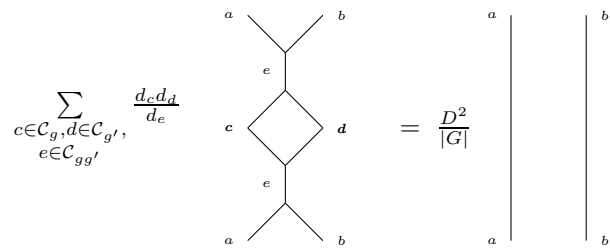
\end{lemma}
\begin{proof}
\begin{align*}
\begin{array}{cclcl}
\text{LHS } &=& \sum\limits_{c,d,e} \frac{d_c d_d}{d_e}N_{cd}^e B_{ab,ab}^e  &=& \sum\limits_{e \in \Cat_{gg'}} \sum\limits_{c\in \Cat_g} \frac{d_c}{d_e}B_{ab,ab}^e \sum\limits_{d \in \Cat_{g'}} d_{d^*}N_{e^* c}^{d^*} \\
            &=& \sum\limits_{e \in \Cat_{gg'}} \sum\limits_{c\in \Cat_g} \frac{d_c^2d_{e^*}}{d_e}B_{ab,ab}^e &=& \sum\limits_{e \in \Cat_{gg'}} D^2_g\, B_{ab,ab}^e \\
            &=& \frac{D^2}{|\text{Gr}(\cross{\Cat}{G})|}Id_{a\otimes b}                              & &
            \end{array}
\end{align*}
\end{proof}

As an application of the above two lemmas and also a warm-up for the proof of invariance under Pachner moves, we compute the invariant for the $4$-sphere $\Sphere^4$. 
\begin{proposition}
\label{prop:S^4_eval}
For any $\cross{\Cat}{G}$, we have $Z_{\cross{\Cat}{G}}(\Sphere^4) = \frac{D^2}{|G|}$.
\begin{proof}
Take the triangulation $\T$ of $\Sphere^4$ which is obtained by identifying a positive $4$-simplex $\sigma_1 = +(01234)$ with a negative $4$-simplex $\sigma_2 = -(01234)$. We use Equation \ref{equ:statesumdefrewritten} to compute $Z_{\cross{\Cat}{G}}(\Sphere^4;\T)$. First of all, fix a coloring on $\T^1$ so that it is extendable to $\T^2$. The expression
\begin{align}
\label{equ:S^4_eval}
\sum\limits_{\T^2,\T^3}\frac{\prod\limits_{\beta \in \T^2}d_{\beta}}{\prod\limits_{\tau \in \T^3}d_{\tau}} \hat{Z}^{+}(\sigma_1)\hat{Z}^{-}(\sigma_2)
\end{align} 
can be computed as shown in Figure~\ref{fig:S^4}. For each graph diagram in the figure, the top and the bottom are identified, which means we are taking the trace of the morphisms represented by the diagrams. The coefficients of the diagrams are placed below them. The summation, such as the first one over $\T^2, \T^3$, means summation over all colorings of the simplices contained in the relevant set. Hence, the first term in the figure corresponds to the expression in Equation \ref{equ:S^4_eval}. Also, $\text{Gr}$ denotes the subgroup $\text{Gr}(\cross{\Cat}{G})$. We explain the four steps of calculations shown in the figure.
\begin{enumerate}[label=\protect\circled{\arabic*}]
\item Apply Lemma \ref{lem:trace identity}.
\item Apply Lemma \ref{lem:dimension identity} twice.
\item Apply Lemma \ref{lem:dimension identity} a third time and isotope the diagram.
\item Use the fact that each $D^2_g$ equals $\frac{D^2}{|\text{Gr}|}$ for $g \in \text{Gr}$.
\end{enumerate}
Then we have
\begin{align*}
\begin{array}{ccccccc}
Z_{\cross{\Cat}{G}}(\Sphere^4;\T) &=& \sum\limits_{\T^1} \frac{(D^2/|G|)^5}{(D^2)^{10}} \left(\frac{D^2}{|\text{Gr}|}\right)^6 &=& |G/\text{Gr}|^4\, |\text{Gr}|^{10}\, \frac{(D^2/|G|)^5}{(D^2)^{10}} \left(\frac{D^2}{|\text{Gr}|}\right)^6  &=& \frac{D^2}{|G|},
\end{array}
\end{align*}
where in the middle equality is due to the fact that the number of extendable colorings on $1$-simplices is $|G/\text{Gr}|^4\, |\text{Gr}|^{10}$.
\begin{figure}[!htbp]
\centering
\begin{tikzpicture}[scale = 0.5] 
   \begin{scope}[xshift = 0cm]
     \IshapeL{1}{1}{0234}{-}{+}
     \IshapeL{1}{4}{0124}{-}{+}
     \Braid{1}{2}
     \draw (3,0) -- (3,2);
     \draw (3,4)--(3,5);
     \IshapeR{2}{6}{1234}{+}{-}
     \draw (0,6) -- (0,7);
     \IshapeL{1}{8}{0134}{+}{-}
     \draw (3,8) -- (3,9);
     \IshapeRP{2}{10}{0123}{+}{-}{\overline{34}}
     \draw (0,10) -- (0,12);
     
     \draw (3.5,-2.5) node{$\sum\limits_{\T^2,\T^3} \frac{\prod\limits_{\beta \in \T^2}d_{\beta}}{\prod\limits_{\tau \in \T^3}d_{\tau}}$};
   \end{scope}

   \begin{scope}[xshift = 4cm]
     \IshapeRP{2}{1}{0123}{-}{+}{\overline{34}}
     \draw (0,0) -- (0,2);
     \IshapeL{1}{3}{0134}{-}{+}
     \draw (3,3) -- (3,4);
     \IshapeR{2}{5}{1234}{-}{+}
     \draw (0,5) -- (0,6);
     \IshapeL{1}{7}{0124}{+}{-}
     \draw (3,7)--(3,8);
     \InvBraid{1}{8}
     \IshapeL{1}{10}{0234}{+}{-}
     \draw (3,10) -- (3,12);
    \end{scope}
    
    \draw (10,6) node{$\overset{\circled{1}}{=}$};
    
    \begin{scope}[xshift = 12cm]
     \IshapeL{1}{1}{0234}{-}{+}
     \IshapeL{1}{4}{0124}{-}{+}
     \Braid{1}{2}
     \draw (3,0) -- (3,2);
     \draw (3,4)--(3,5);
     \IshapeR{2}{6}{1234}{+}{-}
     \draw (0,6) -- (0,9);
     \begin{scope}[yshift = 4cm]
     \IshapeR{2}{5}{1234}{-}{+}
     \draw (0,5) -- (0,6);
     \IshapeL{1}{7}{0124}{+}{-}
     \draw (3,7)--(3,8);
     \InvBraid{1}{8}
     \IshapeL{1}{10}{0234}{+}{-}
     \draw (3,10) -- (3,12);
     \end{scope}
     \draw (1.5,-2.5) node{$\sum\limits_{\substack{\T^2_{1}:= \T^2 - \{(013)\},\\ \T^3_{1}:= \T^3 - \{(0134),(0123)\} }} \frac{\prod\limits_{\beta \in \T^2_{1}}d_{\beta}}{\prod\limits_{\tau \in \T^3_{1}}d_{\tau}}$};
   \end{scope}
   
   \draw (19,6) node{$\overset{\circled{2}}{=}$};
   
   \begin{scope}[xshift = 22cm]
     \InvBraid{1}{0}
     \IshapeL{1}{4}{0124}{-}{+}
     \IshapeL{1}{7}{0124}{+}{-}
      \Braid{1}{2}
      \draw (0,0) -- (0,3);
      \draw (3,4) -- (3,9)node[right]{{\tiny $234$}};
     \draw (1.5,-2.5) node{$\sum\limits_{\substack{\T^2_{2}:= \{(012),(014), \\(024),  (124),(234)\},\\ \T^3_{2}:= \{(0124)\} }} \frac{\prod\limits_{\beta \in \T^2_{2}}d_{\beta}}{\prod\limits_{\tau \in \T^3_{2}}d_{\tau}} (\frac{D^2}{|\text{Gr}|})^2$ };
   \end{scope}
   
   \begin{scope}[xshift = 0cm, yshift = -10cm]
   \draw (0,2) node{$\overset{\circled{3}}{=}$};
   
   \begin{scope}[xshift = 2cm]
   \draw (0,0) -- (0,4)node[above]{{\tiny $024$}};
   \draw (1.5,0) -- (1.5,4)node[above]{{\tiny $\lsup{\overline{24}}{012}$}};
   \draw (3,0) -- (3,4)node[above]{{\tiny $234$}};
   \draw (1.5,-1.5) node{$\sum\limits_{\T^2_{3}:= \{(012),(024),(234)\}} (\prod\limits_{\beta \in \T^2_{3}}d_{\beta}) (\frac{D^2}{|\text{Gr}|})^3$ };
   \end{scope}

   \draw (10,2) node{$\overset{\circled{4}}{=}$};
   
   \draw (14, 2) node{$(\frac{D^2}{|\text{Gr}|})^6$};
   
   \end{scope}

 \end{tikzpicture} 
 \caption{Computing $Z_{\cross{\Cat}{G}}(\Sphere^4;\T)$.}\label{fig:S^4} 
\end{figure}
\end{proof}
\end{proposition}

%
\begin{theorem}
Let $\T, \T'$ be two ordered triangulations which differ by a {\it typical} $3$-$3$ Pachner move, see Table \ref{table:comparison 3-3} for their comparison. Then the identity in Equation \ref{equ:3-3} holds, assuming all simplices which do not present among the summation variables have been assigned a fixed coloring.
\end{theorem}

\begin{proof}
The proof is best illustrated using picture calculus and Lemma \ref{lem:trace identity}.
Note that the top and bottom of each diagram (see below) are identified, so within each diagram morphisms are cyclically ordered, namely, one can move a morphism on the bottom to the top and vice versa. Also, the summation terms are written below the diagrams for convenience. It is direct to check that the final diagrams in both LHS and RHS represent the same morphism, which justifies the equation. Below we give a brief explanation of each \lq\lq = '' sign in the diagrams.

LHS:
\begin{enumerate}[label=\protect\circled{\arabic*}]
 \item By definition.
 \item Apply $\overline{45}$ to the first diagram. Move the top two morphisms involving $0245, 0234$ in the third diagram to the bottom.
 \item Apply Lemma \ref{lem:trace identity} to the second and third diagram.
 \item Isotope the second diagram.
 \item Apply Lemma \ref{lem:trace identity}.
\end{enumerate}

RHS:
\begin{enumerate}[label=\protect\circled{\arabic*}]
 \item By definition.
 \item Apply Lemma \ref{lem:trace identity} to the first and third diagram.
 \item Apply Lemma \ref{lem:trace identity} to the two diagrams.
 \item Isotope the diagram. Note that the group element acting on $0123$ changes from $\overline{35}$ to $\overline{45} \, \overline{34}$ due to braiding.
\end{enumerate}


\begin{figure}[!htbp]
 \begin{tikzpicture}[scale = 0.5]
     \draw (-2,6) node{LHS $\overset{\circled{1}}{=}$ };
     \draw (6,-1.5) node{$\sum\limits_{I_2,I_3} \frac{d_{024}}{d_{0124}d_{0234}d_{0245}}$};
     \draw (20,-1.5) node{$\sum\limits_{I_2,I_3} \frac{d_{024}}{d_{0124}d_{0234}d_{0245}}$};
   \begin{scope}
     \IshapeL{1}{1}{0234}{-}{+}
     \IshapeL{1}{4}{0124}{-}{+}
     \Braid{1}{2}
     \draw (3,0) -- (3,2);
     \draw (3,4)--(3,5);
     \IshapeR{2}{6}{1234}{+}{-}
     \draw (0,6) -- (0,7);
     \IshapeL{1}{8}{0134}{+}{-}
     \draw (3,8) -- (3,9);
     \IshapeRP{2}{10}{0123}{+}{-}{\overline{34}}
     \draw (0,10) -- (0,12);
   \end{scope}

   \begin{scope}[xshift = 4cm]
     \IshapeL{1}{1}{0245}{-}{+}
     \IshapeL{1}{4}{0125}{-}{+}
     \Braid{1}{2}
     \draw (3,0) -- (3,2);
     \draw (3,4)--(3,5);
     \IshapeR{2}{6}{1245}{+}{-}
     \draw (0,6) -- (0,7);
     \IshapeL{1}{8}{0145}{+}{-}
     \draw (3,8) -- (3,9);
     \IshapeRP{2}{10}{0124}{+}{-}{\overline{45}}
     \draw (0,10) -- (0,12);
   \end{scope}

   \begin{scope}[xshift = 8cm]
     \IshapeL{1}{1}{0345}{-}{+}
     \IshapeL{1}{4}{0235}{-}{+}
     \Braid{1}{2}
     \draw (3,0) -- (3,2);
     \draw (3,4)--(3,5);
     \IshapeR{2}{6}{2345}{+}{-}
     \draw (0,6) -- (0,7);
     \IshapeL{1}{8}{0245}{+}{-}
     \draw (3,8) -- (3,9);
     \IshapeRP{2}{10}{0234}{+}{-}{\overline{45}}
     \draw (0,10) -- (0,12);
   \end{scope}

   \begin{scope}[xshift = 12cm]
     \draw (0,6) node{$\overset{\circled{2}}{=}$};
   \end{scope}

   \begin{scope}[xshift = 14cm]
     \IshapeRP{1}{1}{0234}{-}{+}{\overline{45}}
     \IshapeRP{1}{4}{0124}{-}{+}{\overline{45}}
     \Braid{1}{2}
     \draw (3,0) -- (3,2);
     \draw (3,4)--(3,5);
     \IshapeRP{2}{6}{1234}{+}{-}{\overline{45}}
     \draw (0,6) -- (0,7);
     \IshapeRP{1}{8}{0134}{+}{-}{\overline{45}}
     \draw (3,8) -- (3,9);
     \IshapeRP{2}{10}{0123}{+}{-}{}
     \draw (0,10) -- (0,12);
     \draw (0.7,11.1) node{{\tiny $\overline{45}\,\overline{34}$}};
   \end{scope}

    \begin{scope}[xshift = 18cm]
     \IshapeL{1}{1}{0245}{-}{+}
     \IshapeL{1}{4}{0125}{-}{+}
     \Braid{1}{2}
     \draw (3,0) -- (3,2);
     \draw (3,4)--(3,5);
     \IshapeR{2}{6}{1245}{+}{-}
     \draw (0,6) -- (0,7);
     \IshapeL{1}{8}{0145}{+}{-}
     \draw (3,8) -- (3,9);
     \IshapeRP{2}{10}{0124}{+}{-}{\overline{45}}
     \draw (0,10) -- (0,12);
   \end{scope}

   \begin{scope}[xshift = 22cm]
     \IshapeL{1}{1}{0245}{+}{-}
     \draw (3,0) -- (3,2);
     \IshapeRP{2}{3}{0234}{+}{-}{\overline{45}}
     \draw (0,3) -- (0,4);
     \begin{scope}[yshift=4cm]
     \IshapeL{1}{1}{0345}{-}{+}
     \IshapeL{1}{4}{0235}{-}{+}
     \Braid{1}{2}
     \draw (3,1) -- (3,2);
     \draw (3,4)--(3,5);
     \IshapeR{2}{6}{2345}{+}{-}
     \draw (0,6) -- (0,8);
     \end{scope}
   \end{scope}
   \begin{scope}[xshift = 26cm]
     \draw (0,6) node{$\overset{\circled{3}}{=}$};
   \end{scope}

 \begin{scope}[yshift = -27cm]
 \draw (4,-1.5) node{$\sum\limits_{024,0124,0234} \frac{d_{024}}{d_{0124}d_{0234}}$};
 \draw (14,-1.5) node{$\sum\limits_{024,0124,0234} \frac{d_{024}}{d_{0124}d_{0234}}$};
 \begin{scope}[xshift = 0cm]
     \IshapeRP{1}{1}{0234}{-}{+}{\overline{45}}
     \IshapeRP{1}{4}{0124}{-}{+}{\overline{45}}
     \Braid{1}{2}
     \draw (3,0) -- (3,2);
     \draw (3,4)--(3,5);
     \IshapeRP{2}{6}{1234}{+}{-}{\overline{45}}
     \draw (0,6) -- (0,7);
     \IshapeRP{1}{8}{0134}{+}{-}{\overline{45}}
     \draw (3,8) -- (3,9);
     \IshapeRP{2}{10}{0123}{+}{-}{}
     \draw (0,10) -- (0,12);
     \draw (0.7,11.1) node{{\tiny $\overline{45}\,\overline{34}$}};
 \end{scope}

  \begin{scope}[xshift = 4cm]
      \begin{scope}[yshift = -2cm]
     \IshapeL{1}{4}{0125}{-}{+}
     \Braid{1}{2}
     \draw (0,2)--(0,3);
     \draw (3,4)--(3,5);
     \IshapeR{2}{6}{1245}{+}{-}
     \draw (0,6) -- (0,7);
     \IshapeL{1}{8}{0145}{+}{-}
     \draw (3,8) -- (3,9);
     \IshapeRP{2}{10}{0124}{+}{-}{\overline{45}}
     \draw (0,10) -- (0,12);
     \draw (0,12) --(0,14);
     \draw (1,12) -- (1,14);

     \draw plot[smooth] coordinates{(3,2)(3.4,2)(3.4,12)(3,12)};
     \end{scope}

     \begin{scope}[yshift = 10cm]
     \draw (0,2) -- (0,3);
     \IshapeRP{2}{3}{0234}{+}{-}{\overline{45}}
     \draw (0,3) -- (0,4);
     \begin{scope}[yshift=4cm]
     \IshapeL{1}{1}{0345}{-}{+}
     \IshapeL{1}{4}{0235}{-}{+}
     \Braid{1}{2}
     \draw (3,1) -- (3,2);
     \draw (3,4)--(3,5);
     \IshapeR{2}{6}{2345}{+}{-}
     \draw (0,6) -- (0,8);
     \end{scope}
      \draw plot[smooth] coordinates{(3,2)(3.4,2)(3.4,12)(3,12)};
     \end{scope}

   \end{scope}

   \begin{scope}[xshift = 8.5cm]
    \draw (0,6) node {$\overset{\circled{4}}{=}$};
   \end{scope}

   \begin{scope}[xshift = 10cm]
     \IshapeRP{1}{1}{0234}{-}{+}{\overline{45}}
     \IshapeRP{1}{4}{0124}{-}{+}{\overline{45}}
     \Braid{1}{2}
     \draw (3,0) -- (3,2);
     \draw (3,4)--(3,5);
     \IshapeRP{2}{6}{1234}{+}{-}{\overline{45}}
     \draw (0,6) -- (0,7);
     \IshapeRP{1}{8}{0134}{+}{-}{\overline{45}}
     \draw (3,8) -- (3,9);
     \IshapeRP{2}{10}{0123}{+}{-}{}
     \draw (0,10) -- (0,12);
     \draw (0.7,11.1) node{{\tiny $\overline{45}\,\overline{34}$}};
 \end{scope}

 \begin{scope}[xshift = 14cm]
      \begin{scope}[yshift = -2cm]
     \IshapeL{1}{4}{0125}{-}{+}
     \Braid{1}{2}
     \draw (0,2)--(0,3);
     \draw (3,4)--(3,5);
     \IshapeR{2}{6}{1245}{+}{-}
     \draw (0,6) -- (0,7);
     \IshapeL{1}{8}{0145}{+}{-}
     \draw (3,8) -- (3,9);
     \IshapeRP{2}{10}{0124}{+}{-}{\overline{45}}
     \draw (0,10) -- (0,12);
     \draw (0,12) --(0,14);
     \draw (1,12) -- (1,14);
     \draw (5,2) -- (5,12);
     \InvBraid{3}{12}
     \end{scope}

     \begin{scope}[yshift = 10cm]
     \draw (0,2) -- (0,3);
     \IshapeRP{2}{3}{0234}{+}{-}{\overline{45}}
     \draw (0,3) -- (0,4);
     \begin{scope}[yshift=4cm]
     \IshapeL{1}{1}{0345}{-}{+}
     \IshapeL{1}{4}{0235}{-}{+}
     \Braid{1}{2}
     \draw (3,1) -- (3,2);
     \draw (3,4)--(3,5);
     \IshapeR{2}{6}{2345}{+}{-}
     \draw (0,6) -- (0,8);
     \end{scope}
     \draw (5,2) -- (5,12);
     \Braid{3}{12}
     \draw (0,12) -- (0,14);
     \draw (1,12) -- (1,14);
     \end{scope}
   \end{scope}

    \begin{scope}[xshift = 20cm]
    \draw (0,6) node {$\overset{\circled{5}}{=}$};
    \end{scope}

    \begin{scope}[xshift = 22cm]
      \begin{scope}[yshift = -2cm]
     \IshapeL{1}{4}{0125}{-}{+}
     \Braid{1}{2}
     \draw (0,2)--(0,3);
     \draw (3,4)--(3,5);
     \IshapeR{2}{6}{1245}{+}{-}
     \draw (0,6) -- (0,7);
     \IshapeL{1}{8}{0145}{+}{-}
     \draw (3,8) -- (3,9);
     \IshapeRP{4}{10}{1234}{+}{-}{\overline{45}}
     \draw (2,10) -- (2,11);
     \IshapeRP{3}{12}{0134}{+}{-}{\overline{45}}
     \draw (5,12) -- (5,13);
     \IshapeRP{4}{14}{0123}{+}{-}{}
     \draw (2,14) -- (2,16);
     \draw (2.7,15.1) node{{\tiny $\overline{45}\,\overline{34}$}};
     \draw (5,2) -- (5,9);
     \draw (0,10) -- (0,16);
     \draw (3,16) -- (3,17);
     \end{scope}

     \begin{scope}[yshift = 14cm]
     \IshapeL{1}{1}{0345}{-}{+}
     \IshapeL{1}{4}{0235}{-}{+}
     \Braid{1}{2}
     \draw (3,1) -- (3,2);
     \draw (3,4)--(3,5);
     \IshapeR{2}{6}{2345}{+}{-}
     \draw (0,6) -- (0,8);
     \draw (5,0) -- (5,8);
     \Braid{3}{8}
     \draw (0,8) -- (0,10);
     \draw (1,8) -- (1,10);
     \end{scope}

   \end{scope}

 \end{scope}
 \end{tikzpicture}
\end{figure}

\begin{figure}[!htbp]
 \begin{tikzpicture}[scale = 0.5]
  \draw (-2,6) node{RHS $\overset{\circled{1}}{=}$ };
  \draw (6,-1.5) node{$\sum\limits_{I'_2,I'_3} \frac{d_{135}}{d_{0135}d_{1235}d_{1345}}$};
   \begin{scope}
     \IshapeL{1}{1}{0235}{-}{+}
     \IshapeL{1}{4}{0125}{-}{+}
     \Braid{1}{2}
     \draw (3,0) -- (3,2);
     \draw (3,4)--(3,5);
     \IshapeR{2}{6}{1235}{+}{-}
     \draw (0,6) -- (0,7);
     \IshapeL{1}{8}{0135}{+}{-}
     \draw (3,8) -- (3,9);
     \IshapeRP{2}{10}{0123}{+}{-}{\overline{35}}
     \draw (0,10) -- (0,12);
   \end{scope}

   \begin{scope}[xshift = 5cm]
     \IshapeL{1}{1}{0345}{-}{+}
     \IshapeL{1}{4}{0135}{-}{+}
     \Braid{1}{2}
     \draw (3,0) -- (3,2);
     \draw (3,4)--(3,5);
     \IshapeR{2}{6}{1345}{+}{-}
     \draw (0,6) -- (0,7);
     \IshapeL{1}{8}{0145}{+}{-}
     \draw (3,8) -- (3,9);
     \IshapeRP{2}{10}{0134}{+}{-}{\overline{45}}
     \draw (0,10) -- (0,12);
   \end{scope}

   \begin{scope}[xshift = 10cm]
     \IshapeL{1}{1}{1345}{-}{+}
     \IshapeL{1}{4}{1235}{-}{+}
     \Braid{1}{2}
     \draw (3,0) -- (3,2);
     \draw (3,4)--(3,5);
     \IshapeR{2}{6}{2345}{+}{-}
     \draw (0,6) -- (0,7);
     \IshapeL{1}{8}{1245}{+}{-}
     \draw (3,8) -- (3,9);
     \IshapeRP{2}{10}{1234}{+}{-}{\overline{45}}
     \draw (0,10) -- (0,12);
   \end{scope}

   \begin{scope}[xshift = 15cm]
     \draw (0,6) node{$\overset{\circled{2}}{=}$};
   \end{scope}

 \begin{scope}[yshift = -25cm]
 \draw (4,-1.5) node{$\sum\limits_{135, 0135, 1345} \frac{d_{135}}{d_{0135}d_{1345}}$};
  \begin{scope}[xshift = 0cm]
     \IshapeL{1}{1}{0235}{-}{+}
     \IshapeL{1}{4}{0125}{-}{+}
     \Braid{1}{2}
     \draw (3,0) -- (3,2);
     \draw (3,4)--(3,5);
     \IshapeR{4}{6}{2345}{+}{-}
     \draw (2,6) -- (2,7);
     \IshapeL{3}{8}{1245}{+}{-}
     \draw (5,8) -- (5,9);
     \IshapeRP{4}{10}{1234}{+}{-}{\overline{45}}
     \draw (2,10) -- (2,11);
     \IshapeL{3}{12}{1345}{-}{+}
     \draw (5,12) -- (5,13);
     \Braid{3}{13}

     \IshapeL{1}{15}{0135}{+}{-}
     \draw (3,15) -- (3,16);
     \IshapeRP{2}{17}{0123}{+}{-}{\overline{35}}
     \draw (0,17) -- (0,19);

     \draw (5,0) -- (5,5);
     \draw (0,6) -- (0,14);
     \draw (5,15) -- (5,19);
   \end{scope}

   \begin{scope}[xshift = 6cm]
     \IshapeL{1}{1}{0345}{-}{+}
     \IshapeL{1}{4}{0135}{-}{+}
     \Braid{1}{2}
     \draw (3,0) -- (3,2);
     \draw (3,4)--(3,5);
     \IshapeR{2}{6}{1345}{+}{-}
     \draw (0,6) -- (0,7);
     \IshapeL{1}{8}{0145}{+}{-}
     \draw (3,8) -- (3,9);
     \IshapeRP{2}{10}{0134}{+}{-}{\overline{45}}
     \draw (0,10) -- (0,12);
   \end{scope}

    \begin{scope}[xshift = 10cm]
     \draw (0,6) node{$\overset{\circled{3}}{=}$};
   \end{scope}

    \begin{scope}[xshift = 12cm]
     \IshapeL{1}{1}{0235}{-}{+}
     \IshapeL{1}{4}{0125}{-}{+}
     \Braid{1}{2}
     \draw (3,0) -- (3,2);
     \draw (3,4)--(3,5);
     \IshapeR{4}{6}{2345}{+}{-}
     \draw (2,6) -- (2,7);
     \IshapeL{3}{8}{1245}{+}{-}
     \draw (5,8) -- (5,9);
     \IshapeRP{4}{10}{1234}{+}{-}{\overline{45}}

     \IshapeL{1}{11}{0145}{+}{-}
     \IshapeRP{2}{13}{0134}{+}{-}{\overline{45}}
     \draw (0,13) -- (0,14);
     \IshapeL{1}{15}{0345}{-}{+}
     \Braid{1}{16}
     \Braid{3}{18}
     \IshapeRP{2}{21}{0123}{+}{-}{\overline{35}}

     \draw (5,0) -- (5,5);
     \draw (0,6) -- (0,10);
     \draw (5,12) -- (5,18);
     \draw (3,15) -- (3,16);
     \draw (1,18) -- (1,20);
     \draw (0,17) -- (0,23);
     \draw (5,20) -- (5,23);

   \end{scope}

   \begin{scope}[xshift = 18cm]
     \draw (0,6) node{$\overset{\circled{4}}{=}$};
   \end{scope}

   \begin{scope}[xshift = 20cm]
     \IshapeL{1}{1}{0235}{-}{+}
     \IshapeL{1}{4}{0125}{-}{+}
     \Braid{1}{2}
     \draw (3,0) -- (3,2);
     \draw (3,4)--(3,5);
     \IshapeR{4}{6}{2345}{+}{-}
     \draw (2,6) -- (2,7);
     \IshapeL{3}{8}{1245}{+}{-}
     \draw (5,8) -- (5,9);
     \IshapeRP{4}{10}{1234}{+}{-}{\overline{45}}

     \IshapeL{1}{11}{0145}{+}{-}
     \IshapeRP{2}{13}{0134}{+}{-}{\overline{45}}
     \draw (0,13) -- (0,14);
     \IshapeL{1}{15}{0345}{-}{+}
     \IshapeRP{4}{16}{0123}{+}{-}{}
     \draw (2.7,17.1) node{{\tiny $\overline{45}\, \overline{34}$}};
     \Braid{2}{17}
     \Braid{3}{18}

     \draw (5,0) -- (5,5);
     \draw (0,6) -- (0,10);
     \draw (5,12) -- (5,15);
     \draw (0,17) -- (0,20);
     \draw (2,19) -- (2,20);

   \end{scope}

%
%
%
%
%
%
%
%
  \end{scope}

 \end{tikzpicture}
\end{figure}

\end{proof}


\section{Monoidal $2$-categories with Duals}
\label{sec:twocat}
By a $2$-category we always mean a {\it weak} $2$-category or a bicategory. In comparison, a strict $2$-category is a $2$-category where the composition of $1$-morphisms and the horizontal composition of $2$-morphisms are strictly associative. A monoidal $2$-category by definition is a tricategory with a single object. Monoidal $2$-categories with certain extra structure/properties such as duals, sphericity, and semisimplicity defined in a certain sense are supposed to be the input data for a state sum construction of $(3+1)$-$\TQFT$s. We call such categories spherical fusion $2$-categories. The main purpose of this section is to show that from a \CatNameShort{$G$} we can construct a monoidal $2$-category with duals and some additional structures. However, this monoidal $2$-category is not a spherical $2$-category in the sense of \cite{mackaay1999spherical}, while it does become a spherical $2$-category if we adapt the definition of sphericity introduced in \cite{barrett2012gray} for Gray categories to monoidal $2$-categories. 



\subsection{Constructing a Monoidal $2$-category from a \CatNameShort{$G$}}
For monoidal $2$-categories, we follow closely the conventions in \cite{schommer2011classification} (Section $2.3$ and Appendices $A$ and $C$), where a bicategory corresponds to a $2$-category presented here. See \cite{kapanov1994category, Baez2003705} for some relatively older references. Let $\cross{\Cat}{G} = \bigoplus\limits_{g \in G} \Cat_g$ be a \CatNameShort{$G$}. We construct a monoidal $2$-category $\D = \D(\cross{\Cat}{G})$ from $\cross{\Cat}{G}$.

The objects of $\D$ correspond to elements of $G$. That is, $\D^0 = G$. If $a,b \in \D^0$ are two objects, the category $\D(a,b)$ of $1$-morphisms from $a$ to $b$ is given by $\Cat_{\bar{a}b}$, where $\bar{a}$ means the inverse of $a$ as a group element. Thus, $1$-morphisms from $a$ to $b$ are objects, denoted by $A, B, C$, etc., in $\Cat_{\bar{a}b}$, and $2$-morphisms from $A$ to $B$ are morphisms, denoted by $\alpha, \beta$, etc., in $\Cat_{\bar{a}b}$ from $A$ to $B$. The horizontal composition functors are defined by the tensor products in $\cross{\Cat}{G}$, namely,   
\begin{align*}
\begin{array}{cccc}
c_{abc}\,\,\colon  &\D(b,c) \times \D(a,b)  &\longrightarrow  &\D(a,c),\\
                     & (B, A)     &\longmapsto     & A \otimes B,\\
                    & (\beta, \alpha) &\longmapsto      & \alpha \otimes \beta.
\end{array}
\end{align*}
The identity $1$-morphism $I_a \in \D(a,a)$ is the unit object $\unit$ of $\Cat_{e} \subset \cross{\Cat}{G}$. Clearly, the  associator and left/right unitors of the horizontal composition functors in $\D$ correspond precisely to those of the tensor products in $\cross{\Cat}{G}$. This makes $\D$ into a $2$-category. Note that so far we have only used the fact that $\cross{\Cat}{G}$ is a $G$-graded monoidal category where the tensor product respects to the group multiplication.

To make a distinction between the tensor product \lq$\otimes$' in $\cross{\Cat}{G}$ which corresponds to the horizontal composition in $\D$ and the tensor product in $\D$, we denote the latter by the symbol \lq$\boxtimes$'. The distinguished object (or the unit object) $\unit_{\D}$ is identity element $e$ of $G$. The tensor product $\boxtimes$ is a homomorphism of $2$-categories,
\begin{align*}
\begin{array}{cccc}
\boxtimes = (\boxtimes, \phi_{(A,A'),(B,B')}, \phi_{(a,a')})\,\, \colon & \D \times \D &\longrightarrow & \D,
\end{array}
\end{align*}
which is unpacked as follows. For $a,b \in \D^0$, we have $a \boxtimes b:= ab$. The functor 
\begin{align*}
\begin{array}{cccc}
\boxtimes\,\,\colon & \D(a,b) \times \D(a',b') &\longrightarrow & \D(aa',bb'),
\end{array}
\end{align*}
is defined by
\begin{align*}
A \boxtimes A' := \lsup{\overline{a'}}{A} \otimes A', \qquad \alpha \boxtimes \alpha' := \lsup{\overline{a'}}{\alpha} \otimes \alpha'. 
\end{align*}
For
\begin{equation*}
\begin{tikzcd}
a \arrow[r, "A"] & b \arrow[r, "B"] & c,
\end{tikzcd}
\qquad
\begin{tikzcd}
a' \arrow[r, "A'"] & b' \arrow[r, "B'"] & c',
\end{tikzcd}
\end{equation*}
the natural isomorphisms
\begin{equation*}
\begin{array}{rrcl}
\phi_{(A,A'),(B,B')} \,\,\colon & (\lsup{\overline{a'}}{A} \otimes A') \otimes (\lsup{\overline{b'}}{B} \otimes B') & \longrightarrow  & \lsup{\overline{a'}}{(A \otimes B)} \otimes (A'\otimes B'), \\
\phi_{(a,a')} \,\,\colon & \lsup{\overline{a'}}{\unit} \otimes \unit & \longrightarrow & \unit
\end{array}
\end{equation*}
are given, respectively, by the diagrams
\begin{equation*}
\begin{array}{c}
\begin{tikzcd}
(\lsup{\overline{a'}}{A} \otimes A') \otimes (\lsup{\overline{b'}}{B} \otimes B') \arrow[rr, "id \otimes c \otimes id"] && (\lsup{\overline{a'}}{A} \otimes \lsup{\overline{a'}}{B}) \otimes (A'\otimes B') \arrow[r, " "] & \lsup{\overline{a'}}{(A \otimes B)} \otimes (A'\otimes B'), 
\end{tikzcd} 
\\
\begin{tikzcd}
\lsup{\overline{a'}}{\unit} \otimes \unit \arrow[r, " "] &\unit \otimes \unit \arrow[r, " "] &\unit,
\end{tikzcd}
\end{array}
\end{equation*}
where $c$ is the $G$-crossed braiding and the arrows without labels denote natural isomorphisms arising from the action of $\overline{a'}$. Note that in the first diagram above, some associators need to inserted in order for the maps to make sense.
Clearly, $\phi_{(A,A'),(B,B')}$ is natural with respect to its arguments. The properties of the $G$-crossed braiding and the $G$-action ensure that $\phi_{(A,A'),(B,B')}$ and $\phi_{(a,a')}$ satisfy the required equations. The three transformations relating the associators and unitors of $\boxtimes$,
\begin{align*}
\begin{array}{rrcl}
(\alpha_{abc}, \alpha_{ABC})=\alpha\,\, \colon & (a \boxtimes b) \boxtimes c &\longrightarrow & a \boxtimes (b \boxtimes c), \\
(l_a,l_A) = l\,\,\colon & e \boxtimes a & \longrightarrow & a,\\
(r_a, r_A) = r\,\,\colon & a &\longrightarrow & a \boxtimes e,
\end{array}
\end{align*}
are defined as follows. The $1$-morphisms $\alpha_{abc}$, $l_a$, and $r_a$ are all equal to $\unit$ and the $2$-morphisms $\alpha_{ABC}$, $l_A$, and $r_A$ are given by the the natural isomorphisms,
\begin{align*}
\begin{array}{rrcl}
\alpha_{ABC}\,\,\colon & \unit \otimes \left(\lsup{\overline{bc}}{A} \otimes (\lsup{\overline{c}}{B} \otimes C)\right) &\longrightarrow & \left(\lsup{\bar{c}}{(\lsup{\bar{b}}{A} \otimes B)} \otimes C\right)\otimes \unit,\\
l_A\,\,\colon & \unit \otimes A &\longrightarrow & (\lsup{\bar{a}}{\unit} \otimes A) \otimes \unit,\\
r_A\,\,\colon & \unit \otimes (\lsup{e}{A} \otimes \unit) &\longrightarrow & A \otimes \unit.
\end{array}
\end{align*} 
Lastly, we consider the invertible modification $\pi$, the pentagonator, as shown below.
\begin{equation*}
\begin{tikzcd}
     & (a \boxtimes b) \boxtimes (c \boxtimes d) \arrow[dr, "\alpha"]&   \\
((a \boxtimes b) \boxtimes c) \boxtimes d \arrow[ur, "\alpha"]\arrow[d, "\alpha \boxtimes I"]&  \tikz{\node[rotate = 90] at (0,0) {$\Longrightarrow$};\node at (0.3,-0.1) {$\pi$};}    & a \boxtimes (b \boxtimes (c \boxtimes d)) \\
(a \boxtimes (b \boxtimes c)) \boxtimes d \arrow[rr, "\alpha"] &      & a \boxtimes ((b \boxtimes c) \boxtimes d)\arrow[u, "I \boxtimes \alpha"]
\end{tikzcd}
\end{equation*}
It is direct to see that the source and target of $\pi(a,b,c,d)$ are both naturally isomorphic to $\unit$, hence $\pi(a,b,c,d)$ is a nonzero scalar. The equality of the two Stasheff polytopes in Figures C$.1$ and C$.2$ of \cite{schommer2011classification} is equivalent to the cocycle condition,
\begin{align*}
\pi(b,c,d,e) \,\pi(a,bc,d,e)\, \pi(a,b,c,de)\,\, &= \,\,\pi(ab,c,d,e) \,\pi(a,b,cd,e)\, \pi(a,b,c,d).
\end{align*}
Note that here $e$ denotes a general object, but not the unit. This means $\pi$ is a $4$-cocycle representing a class in $H^4(G,\C)$. 

\begin{remark}
There are also three other invertible modifications $\mu$, $\lambda$, and $\rho$ resulting from weakening the Triangle Identity. We define all of them to be identically $1$. To satisfy the conditions relating these modifications and $\pi$, we need to choose the $4$-cocycle $\pi$ to be {\it normalized}, namely, $\pi(a,b,c,d) = 1$ whenever $a$, $b$, $c$, or $d$ is the identity element. But this choice is only for notational convenience. In general, one can always properly define $\mu$, $\lambda$, and $\rho$ by certain values of $\pi$ so that the relevant equations are satisfied.
\end{remark}

To summarize, by taking $\pi \in H^4(G,\C) $ trivial, a monoidal $2$-category $\D(\cross{\Cat}{G})$ is constructed from a \CatNameShort{$G$} $\cross{\Cat}{G}$ with a trivial pentagonator. One can also introduce a nontrivial $\pi$ to $\D(\cross{\Cat}{G})$ so that it becomes a monoidal 2-category with $\pi$ being the pentagonator. From the construction, it is direct to see that $\D$ is semistrict if and only if $\cross{\Cat}{G}$ is strict as a $G$-crossed braided category (see Section \ref{subsec:stratification}) and $\pi$ is identically $1$. An equivalence of $\cross{\Cat}{G}$ as a $G$-crossed braided category implies an equivalence of $\D(\cross{\Cat}{G})$ as a monoidal $2$-category. However, the converse does not seem to be true. 

\begin{remark}
\begin{enumerate}
\item It can be shown that monoidal functors between \CatNameShort{$G$}s and natural transformations between monoidal functors can be lifted to monoidal $2$-functors and $2$-transformations, respectively. However, for monoidal $2$-categories, the more general notions of pseudofunctors, tritransformations, and trimodifications also exist. See, for instance, \cite{gordon1995coherence, schommer2011classification}.
\item In Section \ref{subsec:variation}, the definition of the $4$-manifold invariant is generalized to include a $4$-cocycle in $H^4(G, \C)$ with the restriction that the \CatNameShort{$G$} is concentrated on the sector indexed by the identity element. However, we just showed that a $4$-cocycle can be combined with any  \CatNameShort{$G$} to produce a monoidal $2$-category in which the $4$-cocycle corresponds to the pentagonator. This leads us to speculate that the restriction on the \CatNameShort{$G$} in the definition of the invariant is not necessary. However, it is not clear that Equation \ref{def:newpartition} still gives a well-defined invariant without the restriction.
\end{enumerate}
\end{remark} 

\subsection{Constructing a \lq Spherical' $2$-category from a \CatNameShort{$G$}} 
We continue to show that $\D(\cross{\Cat}{G})$ has duals and certain other structures which makes it a spherical $2$-category in the sense of \cite{barrett2012gray}, but not in the sense of \cite{mackaay1999spherical}. Both of the two references define duals assuming the category is semistrict. This does not lose any generality since every monoidal $2$-category is equivalent to a semistrict one. That means, if one can define duals on a semistrict monoidal $2$-category, then the same definition also works for nonsemistrict ones by inserting appropriate natural isomorphisms whenever needed. Thus here we assume that $\cross{\Cat}{G}$ is strict and hence $\D$ is semistrict. For a direct comparison, we follow the conventions in \cite{mackaay1999spherical} to define duals. See \cite{Baez2003705} for the original definition.

Some notations are in order. Let $\D$ be any semistrict monoidal $2$-category. The distinguished object $\unit_{\D}$ is always denoted by $e$, while $\unit$ is reserved to mean the unit object in a \CatNameShort{$G$}. The identity in the group $G$ is also denoted by $e$. If $A\colon a \rightarrow b$, $B \colon b \rightarrow c$ are $1$-morphisms, their horizontal composition is written as $A \hc B := c_{abc}(B,A)$. Similarly, write $\alpha \hc \beta$ as the horizontal composition of two $2$-morphisms. If $A$ is a $1$-morphism and $\alpha$ is a $2$-morphism, then $A \hc \alpha$ means $I_A \hc \alpha$. Also, if $c$ is any object, then $c \boxtimes A$, $c \boxtimes \alpha$, and $A \boxtimes \alpha$ represent $I_c \boxtimes A$,   $I_{I_c} \boxtimes \alpha$, and $I_A \boxtimes \alpha$, respectively.  For instance, if $\D = \D(\cross{\Cat}{G})$, then $A \hc B = A \otimes B$, $\alpha \hc \beta = \alpha \otimes \beta$, $A \hc \alpha = I_A \otimes \alpha$. Also, if $A, A' \colon a \rightarrow b$, $\alpha\colon A \Rightarrow A'$ then $c \boxtimes A = A$, $A \boxtimes c = \lsup{\bar{c}}{A}$, $c \boxtimes \alpha = \alpha$, $\alpha \boxtimes c = \lsup{\bar{c}}{\alpha}$. If $B\colon c \rightarrow d$, then $B \boxtimes \alpha = \alpha$, $\alpha \boxtimes B = \lsup{\bar{c}}{\alpha}$.

The duals in a monoidal $2$-category  consist of the following structures.
\begin{enumerate}
\item For every $2$-morphism $\alpha\colon A \Rightarrow B$, there is a $2$-morphism $\alpha^{\dag}: B \Rightarrow A$.   
\item For every $1$-morphism $A\colon a \rightarrow b$, there is a $1$-morphism $A^*\colon b \rightarrow a$, a $2$-morphism $i_A\colon I_a \Rightarrow A \hc A^*$, and a $2$-morphism $e_A\colon A^* \hc A \Rightarrow I_{b}$ called the dual, the unit, and the counit of $A$, respectively.
\item For every object $a$, there is an object $a^*$, a $1$-morphism $i_a\colon e \rightarrow a \boxtimes a^*$, and a $1$-morphism $e_a\colon a^* \boxtimes a \rightarrow e$ called the dual, the unit, and the counit of $a$, respectively. Also there is a $2$-morphism $T_a\colon (i_a \boxtimes I_a) \hc (I_a \boxtimes e_a) \Rightarrow I_a$ called the triangulator of $a$.
\end{enumerate}  
A $2$-morphism $\alpha$ is called unitary if it is invertible and $\alpha^{\dag} = \alpha^{-1}$. For $\alpha\colon A \Rightarrow B$, we define a $2$-morphism $\alpha^* \colon B^* \Rightarrow A^*$ by
\begin{equation*}
\begin{tikzcd}
\alpha^*\,\,\colon\,\, B^* \arrow[rr, Rightarrow,"B^* \hc i_A"] & & B^*  \hc A \hc A^* \arrow[rr,Rightarrow, "B^* \hc \alpha \hc A^*"] && B^* \hc B \hc A^* \arrow[rr,Rightarrow, "e_B \hc A^*"] && A^*.
\end{tikzcd}
\end{equation*}
{\it Note that the roles of $\alpha^{\dag}$ and $\alpha^*$ are swapped in} \cite{mackaay1999spherical}. These structures defined above need to satisfy some consistency conditions \cite{mackaay1999spherical}.

To define duals in $\D(\cross{\Cat}{G})$,  we will assume $\cross{\Cat}{G}$ has a unitary structure, namely, for every morphism $\alpha\colon A \rightarrow B$, there is a morphism $\alpha^{\dag}: B \rightarrow A$. The map $(\cdot)^{\dag}$ is involutory and compatible other structures of the category. See Section $4.3$ of \cite{Wang2010topological} or Chapter $2$ of \cite{turaev1994quantum} for the precise definition.

For $\D = \D(\cross{\Cat}{G})$, the structures in the first two items above are  defined by the relevant structures in $\cross{\Cat}{G}$ in the obvious way. Namely, as a $2$-morphism, $\alpha^{\dag}$ corresponds to the unitary structure from $\cross{\Cat}{G}$. The dual $A^*$ of $A$ as a $1$-morphism is the dual of $A$ as an object in $\cross{\Cat}{G}$. For an object $a \in \D$, $a^*$ is defined to be $\bar{a}$, $i_a$ and $e_a$ are equal to $I_e = \unit$, and $T_a$ is the identity $2$-morphism on $I_a = \unit$. It is direct to see that $\alpha^*$ defined above is the same as the dual of $\alpha$ viewed as a morphism in $\cross{\Cat}{G}$. This is the reason why we swapped the notations $(\cdot)^{\dag}$ and $(\cdot)^{*}$. It is straightforward to check that this defines duals in $\D$.

We now proceed to define a pivotal structure on $\D = \D(\cross{\Cat}{G})$. Given a $1$-morphism $A\colon a \rightarrow b$, define the $1$-morphisms $\lsup{\#}{A}, \ A^{\#}\, \colon \bar{b} \rightarrow \bar{a}$ as follows.
\begin{equation*}
\begin{tikzcd}
\lsup{\#}{A}\,\,\colon\,\, \bar{b} \arrow[r, "i_{\bar{a}} \boxtimes \bar{b}"]   & \bar{a} \boxtimes a \boxtimes \bar{b}  \arrow[r, "\bar{a} \boxtimes A \boxtimes \bar{b}"]& \bar{a} \boxtimes b \boxtimes \bar{b} \arrow[r, "\bar{a} \boxtimes e_{\bar{b}}"] & \bar{a}, \\
A^{\#}\,\,\colon\,\, \bar{b} \arrow[r, "\bar{b} \boxtimes i_a"]   & \bar{b} \boxtimes a \boxtimes \bar{a}  \arrow[r, "\bar{b} \boxtimes A \boxtimes \bar{a}"]& \bar{b} \boxtimes b \boxtimes \bar{a} \arrow[r, "e_{b} \boxtimes \bar{a}"] & \bar{a}.
\end{tikzcd}
\end{equation*}
Hence, we have $\lsup{\#}{A} = \lsup{b}{A}$ and $A^{\#} = \lsup{a}{A}$. Note that the above definitions naturally extend to $2$-morphisms. Explicitly, for a $2$-morphism $\alpha\colon A \Rightarrow B$, where $A, B\colon a \rightarrow b$ are $1$-morphisms, define the $2$-morphisms $\lsup{\#}{\alpha}\colon \lsup{\#}{A}\Rightarrow \lsup{\#}{B}$ and $\alpha^{\#}\colon A^{\#} \Rightarrow B^{\#}$ by the following formulas,
\begin{alignat*}{3}
\lsup{\#}{\alpha}\quad &&:=\quad (i_{\bar{a}} \boxtimes \bar{b})\hc (\bar{a} \boxtimes \alpha \boxtimes \bar{b}) \hc (\bar{a} \boxtimes e_{\bar{b}})\quad &&=\quad \lsup{b}{\alpha},\\
\alpha^{\#} \quad      &&:=\quad (\bar{b} \boxtimes i_a) \hc (\bar{b} \boxtimes \alpha \boxtimes \bar{a}) \hc (e_{b} \boxtimes \bar{a})\quad &&=\quad \lsup{a}{\alpha}.
\end{alignat*}
Note that the twist in $\cross{\Cat}{G}$ provides a natural isomorphism,
\begin{align*}
\theta_{A^{\#}}^{-1} = (\lsup{a}{\theta_A})^{-1}\,\,\colon\,\, &\lsup{\#}{A} \Longrightarrow A^{\#}.
\end{align*}
We denote this isomorphism by $\phi_A$. It is direct to verify that $\phi_A$ satisfies Conditions $1$ through $4$ below, which makes $\D$ pivotal according \cite{mackaay1999spherical}.
\begin{enumerate}
\item For any $2$-morphism $\alpha\colon A \Rightarrow B$, the following diagram commutes,
\begin{equation*}
\begin{tikzcd}
\lsup{\#}{A} \arrow[r, Rightarrow, "\phi_A"]\arrow[d, Rightarrow, "\lsup{\#}{\alpha}"] & A^{\#} \arrow[d, Rightarrow, "\alpha^{\#}"] \\
\lsup{\#}{B} \arrow[r, Rightarrow, "\phi_B"] & B^{\#}.
\end{tikzcd}
\end{equation*}
\item $\phi_A^{*} = \phi_{A^*}$.
\item For $A\colon a \rightarrow b,$ $B\colon b \rightarrow c$, the following diagram commutes:
\begin{equation*}
\begin{tikzcd}
\lsup{\#}{(A \hc B)}\arrow[r, Rightarrow, "\phi_{(A \hc B)}"] \arrow[d, Rightarrow] & (A \hc B)^{\#}\arrow[d, Rightarrow] \\
\lsup{\#}{B} \hc \lsup{\#}{A} \arrow[r, Rightarrow, "\phi_B \hc \phi_A"] & B^{\#} \hc A^{\#},
\end{tikzcd}
\end{equation*}
where the vertical maps are given by the $G$-crossed braiding (or its inverse) up to natural isomorphisms arising from the group action.
\item For any $2$-morphism $\alpha$, we have 
\begin{align*}
\lsup{\#}{(\alpha^*)} = (\alpha^{\#})^{*} \quad \text{and} \quad (\lsup{\#}{A})^* = (A^*)^{\#}.
\end{align*}
\end{enumerate}
The conditions above can be derived from properties of the twist. See Proposition \ref{prop:twist_prop}. In particular, the third condition is equivalent to the relation between the twist and $G$-crossed braiding,
\begin{equation*}
\begin{tikzcd}
\lsup{c}{(A \otimes B)} \arrow[r, "\theta^{-1}"]\arrow[d, "c^{-1}_{(\lsup{c}{B}, \lsup{b}{A})}"] &\lsup{a}{(A \otimes B)} \arrow[d, "c_{(\lsup{a}{A},\lsup{a}{B})}"] \\
\lsup{c}{B} \otimes \lsup{b}{A} \arrow[r, "\theta^{-1} \otimes \theta^{-1}"] & \lsup{b}{B} \otimes \lsup{a}{A}. 
\end{tikzcd}
\end{equation*}

We introduce the left/right trace functors $\Tr_L, \Tr_R\colon \D(a,a) \rightarrow \D(e,e)$. For a $1$-morphism $A\colon a \rightarrow a$, $\Tr_L(A)$ and $\Tr_R(A)$ are defined respectively by,
\begin{equation*}
\begin{tikzcd}
\Tr_L(A)\,\,\colon\,\, e \arrow[r, "i_{\bar{a}}"] & \bar{a} \boxtimes a \arrow[r, "\bar{a} \boxtimes A"] & \bar{a} \boxtimes a \arrow[r, "e_a"] & e,\\
\Tr_R(A)\,\,\colon\,\, e \arrow[r, "i_{a}"] & a \boxtimes \bar{a} \arrow[r, "A \boxtimes \bar{a}"] & a \boxtimes \bar{a} \arrow[r, "e_{\bar{a}}"] & e
\end{tikzcd}
\end{equation*}
Similarly, for a $2$-morphism $\alpha\colon A \Rightarrow B$, $A, B \colon a \rightarrow b$, define
\begin{align*}
\Tr_L(\alpha)\,\,&:= \,\,i_{\bar{a}} \hc (\bar{a} \boxtimes \alpha) \hc e_a \\
\Tr_R(\alpha)\,\,&:=\,\, i_{a} \hc (\alpha \boxtimes \bar{a}) \hc e_{\bar{a}}.
\end{align*}
It turns out that $\Tr_L\colon \Cat_e \rightarrow \Cat_e$ is the identity functor and $\Tr_R$ is the functor given by the action of $a$ on $\Cat_e$. 

According to Definition $2.8$ of \cite{mackaay1999spherical}, $\D$ is called {\it spherical} if there is a natural isomorphism between $\Tr_L$ and $\Tr_R$. This means the identity functor and the $\lsup{a}{(\cdot)}$ functor on $\Cat_e$ are isomorphic for any group element $a$. In particular, it implies the action of $G$ fixes all objects of $\Cat_e$. This clearly does not hold in general. For instance, take any finite group $G$ acting nontrivially on an Abelian group $H$ and let $\rho\colon H \rightarrow G$ be the trivial morphism, then $(H,G, \rho)$ forms a crossed module by Section \ref{subsec:yetter}, and the corresponding \CatNameShort{$G$} from the crossed module has only one nontrivial sector, the sector indexed by the identity element which is equivalent to $\Vect_H$. Moreover, the $G$-action on $\Vect_H$ is given by the nontrivial $G$-action on $H$.

In \cite{mackaay1999spherical}, a monoidal $2$-category is defined to be finitely semisimple if there is a finite nonempty set of objects $S$, such that for any pair of objects $(a,b)$,
\begin{align}
\label{equ:semisimplicity}
\bigoplus\limits_{c \in S} \D(a,c) \boxtimes \D(c,b) \overset{\simeq}{\longrightarrow} \D(a,b).
\end{align}
We show $\D$ in general does not satisfy this property. Take $G = \{e\}$, then $\cross{\Cat}{G} = \Cat$ is just a ribbon fusion category, $\D$ has a single object $e$, and $\D(e,e) = \Cat$. If $\D$ were semisimple, then $S = \{e\}$, and $\Cat \boxtimes \Cat \simeq \Cat$, which fails to hold in general.

To summarize, starting from a unitary \CatNameShort{$G$}, a monoidal $2$-category with duals can be constructed, which in general does not satisfy the axioms of a semisimple spherical $2$-category. The requirement that the left trace functor be isomorphic to the right trace functor seems too strong.  Also, note that $\D(a,b)$ is a $\D(a,a)$-$\D(b,b)$ bi-module category.  For the semisimplicity condition, it seems more reasonable to require the tensor product in Equation \ref{equ:semisimplicity} to be taken over $\D(c,c)$ rather than over $\Vect$. 

On the other hand, spherical Gray categories are defined in \cite{barrett2012gray} where a Gray category is a semistrict tricategory. A semistrict monoidal $2$-category can be viewed as a Gray category with one object. We translate the definition of spherical Gray categories (Definitions $7.1$ and $7.2$ of \cite{barrett2012gray}) to a semistrict monoidal $2$-category $\D$. For any $2$-morphism $\alpha\colon A \Rightarrow A$, $A\colon a \rightarrow b$, the (new) left trace $\text{tr}_L(\alpha)\colon I_b \Rightarrow I_b$ and the (new) right trace $\text{tr}_R(\alpha)\colon I_a \Rightarrow I_a$ are defined by the following diagrams,
\begin{equation*}
\begin{tikzcd}
\text{tr}_L(\alpha)\,\,\colon\,\, I_b \arrow[r, Rightarrow, "e_A^{\dag}"] & A^* \hc A \arrow[r, Rightarrow, "A^* \hc \alpha"] & A^* \hc A \arrow[r, Rightarrow, "e_A"] & I_b,\\
\text{tr}_R(\alpha)\,\,\colon\,\, I_a \arrow[r, Rightarrow, "i_A"] & A \hc A^* \arrow[r, Rightarrow, "\alpha \hc A^*"] & A \hc A^* \arrow[r, Rightarrow, "i_A^{\dag}"] & I_a.
\end{tikzcd}
\end{equation*}
Note that the $2$-morphism $\lsup{\#}{\alpha}\colon \lsup{\#}{A} \Rightarrow \lsup{\#}{A}$ defined in this subsection is the same as the $2$-morphism $\#\alpha\colon \#A \Rightarrow \#A$ defined in \cite{barrett2012gray} (Section $4.2$). By Definition $7.2$ of \cite{barrett2012gray}, $\D$ is called {\it spherical} if the following identities hold.
\begin{alignat}{3}
\label{equ:new spherical}
(\text{tr}_L(\alpha) \boxtimes b^*)\hc e_{b^*}\quad &&= \quad(b \boxtimes \text{tr}_R(\lsup{\#}{\alpha})) \hc e_{b^*}, \quad\quad (\text{tr}_R(\alpha) \boxtimes a^*)\hc e_{a^*}\quad &&= \quad(a \boxtimes \text{tr}_L(\lsup{\#}{\alpha})) \hc e_{a^*}.  
\end{alignat}

For $\D = \D(\cross{\Cat}{G})$ for a unitary $\cross{\Cat}{G}$, it is direct to see that $\text{tr}_L(\alpha) = \text{tr}_R(\alpha) = \Tr(\alpha)$, where the last term means the trace of $\alpha$ in $\cross{\Cat}{G}$, and Equation \ref{equ:new spherical} reduces to the condition that $\Tr(\alpha) = \Tr(\lsup{b}{\alpha})$ which always holds in a \CatNameShort{$G$}. Hence $\D(\cross{\Cat}{G})$ is spherical monoidal $2$-category in the sense of \cite{barrett2012gray}.

Therefore, one may propose a spherical $2$-category as a spherical Gray category  with one object. Or more precisely, a spherical $2$-category is a semistrict monoidal $2$-category with duals that satisfy Equation \ref{equ:new spherical}. 

\begin{remark}
In \cite{barrett2012gray} (Definition $7.5$), semisimple Gray categories is defined. It can be checked that $\D = \D(\cross{\Cat}{G})$ is a semisimple Gray category. It should be noted that the notion of semisimple Gray categories is not sufficient for defining {\it fusion} 2-categories. This is because, according to that definition (adapted to monoidal $2$-categories), only semisimplicity at the level of $2$-morphisms is required. That is, for any object $a,b$, there is a set $J$ of simple 1-morphisms from $a$ to $b$, such that for any two 1-morphism $A, B \colon a \to b$, the composition gives a natural isomorphism:
\begin{align*}
\bigoplus\limits_{C \in J} \Hom(A, C) \otimes \Hom(C,B) \simeq \Hom(A,B).
\end{align*}
We speculate that in defining fusion 2-categories, \lq finite simplicity' at the level of 1-morphisms is also needed and this is a nontrivial issue.
\end{remark}

To this end, we add a comment on Gray diagrams which are used in \cite{barrett2012gray} to represent morphisms in a Gray category. Gray diagrams are generalizations of ribbon graphs for ribbon categories. Roughly, a Gray diagram is a three dimensional complex contained in a unit cube where the $i$-cells are labeled with $(3-i)$-morphisms with certain restrictions. Each Gray diagram can be evaluated to a $3$-morphism. We have shown that from a strict \CatNameShort{$G$} we can construct a Gray category with one object. On the other hand, morphisms in a \CatNameShort{$G$} can be represented by two dimensional graph diagrams (see Section \ref{subsec:stratification}). In fact, it can be shown that the graph diagrams for a \CatNameShort{$G$} correspond to the evaluation of the Gray diagrams for the Gray category constructed from that \CatNameShort{$G$}. In the language of \cite{barrett2012gray} (Definition $2.25$), evaluation of Gray diagrams means projecting the Gray diagrams in the cube with $(w,x,y)$-axis to the $(x,y)$-plane.

\section{Open Questions and Future Directions}
\label{sec:future}
In this paper, we constructed an invariant of $4$-manifolds out of a \CatName{$G$} category. The construction is a $4D$ analog of the Turaev-Viro invariants of $3$-manifolds. Here we point out a few directions$/$questions for future study.

\vspace{0.5cm}
\begin{enumerate}
\item What is the power of the invariant in terms of distinguishing $4$-manifolds? 

The Crane-Yetter invariant from a modular category is a classical invariant, which can be expressed in terms of the signature and the Euler characteristic \cite{Crane93evaluatingthe}. It would be interesting to also give the invariant from a \CatNameShort{$G$} an interpretation in terms of the intrinsic properties of the manifolds. We speculate that our invariant is related to the homotopy $3$-types with the exact relation to be studied in the future.

In \cite{kasprowski2018} it was shown that if an invariant of $4$-manifolds is multiplicative (up to a constant nonzero scalar) under connected sum and is invertible on $\mathbb{CP}^2$ and $\overline{\mathbb{CP}}^2$, then the invariant cannot distinguish smooth structures. On the other hand, for the invariant $Z_{\cross{\Cat}{G}}(\cdot)$, assuming it can be extended to a $\TQFT$\footnote{From the procedure of the construction, it is plausible to extend the invariant to a $\TQFT$. See \cite{williamson2016hamiltonian} for a Hamiltonian realization of the invariant.}, it is well known that,
\begin{align*}
Z_{\cross{\Cat}{G}}(M_1 \# M_2)Z_{\cross{\Cat}{G}}(\Sphere^4) &= Z_{\cross{\Cat}{G}}(M_1)Z_{\cross{\Cat}{G}}(M_2),
\end{align*}
for any two closed $4$-manifolds $M_1$ and $M_2$ if the vector space associated with $\Sphere^3$ is $1$-dimensional. Note that in Proposition \ref{prop:S^4_eval} we computed the invariant for $\Sphere^4$ which equals $D^2/|G|\neq 0$. It will be interesting to see if the invariants of $\mathbb{CP}^2$ and $\overline{\mathbb{CP}}^2$ are always nonzero. We have not found a simple triangulation for either of these two manifolds to make the calculation of their invariants feasible. Also, it needs to be checked if the vector space of $\Sphere^3$ is $1$-dimensional.

\vspace{0.3cm}
\item The definition of the current invariant is based on a triangulation of $4$-manifolds. In terms of calculations, this is not very efficient since the triangulations of most interesting $4$-manifolds contain a fairly large number of simplices, which makes it impractical to compute the invariant. It would be useful if the invariant can be defined in terms of other presentations of the manifolds such as Kirby diagrams. Note that the Crane-Yetter invariant indeed has such a formulation \cite{roberts1995skein} \cite{barenz2016dichromatic}.
\vspace{0.3cm}
\item It is appealing to generalize the invariant in the following three directions.

In \cite{yetter1993homologically} and \cite{roberts1996refined}, a refined version of the Crane-Yetter invariant was introduced, namely, the invariant was associated with a pair $(M, \omega)$, where $M$ is a closed oriented $4$-manifold and $\omega \in H_2(M, \Z_2)$, and the original Crane-Yetter invariant is a normalized sum of the refined invariant over all $\omega \in H_2(M, \Z_2)$. Moreover, the refined invariant gives a state sum formula of the second Stiefel-Whitney class and the Pontrjagin squares of second cohomology classes. It is interesting to see if our invariant also has a similar refinement.

Secondly, the extension of a braided spherical fusion category $\Cat$  to a \CatNameShort{$G$} $\cross{\Cat}{G}$ with $\Cat_e = \Cat$ depends on the vanishing condition of two obstructions, one in $O_3 \in H^3(G,\text{Inv}(\Cat))$ and other in $O_4 \in H^4(G,U(1))$ \cite{etingof2010fusion}. We conjecture that even if the $O_4$ obstruction does not vanish, one can still get a $4$-manifold invariant from some structures beyond \CatNameShort{$G$}s.

Lastly, we are interested in defining a \CatNameShort{$G$} where $G$ is allowed to be an infinite compact Lie group, such as $U(1)$, and in modifying the state sum model so that the partition function converges. It is expected that the resulting invariant would be stronger than the current one.

\vspace{0.3cm}
\item We already know that the invariant defined from a \CatNameShort{$G$} reduces to the Crane-Yetter invariant when $G$ is trivial. For any finite group $G$, by \cite{drinfeld2010braided}\cite{kirillov2002modular}\cite{kirillov2002modular2}, equivalence classes of \CatNameShort{$G$}s with a faithful $G$-grading are in one-to-one correspondence, by equivariantization and de-equivariantization, with equivalence classes of ribbon fusion categories containing $Rep(G)$ as a full subcategory. Thus we wonder in the general case whether there are any relations between the invariant from a $\GBSFC$ with a faithful grading and the Crane-Yetter invariant from the corresponding ribbon fusion category. On the other hand, in \cite{petit2008dichromatic} \cite{barenz2016dichromatic} an invariant of $4$-manifolds was defined for any pivotal functor $F: \Cat \longrightarrow \D$ where $\Cat$ is a spherical fusion category and $\D$ is a ribbon fusion category. It will be interesting to see if our invariant is related to theirs.

\vspace{0.3cm}
\item In $(2+1)$ dimensions, a nonsemisimple generalization of the Turaev-Viro invariant is the Kuperberg invariant \cite{kuperberg1996noninvolutory}. One may wonder if there is a four dimensional analog of the Kuperberg invariant, which generalizes the invariant in the current paper. The Kuperberg invariant is defined on a Heegaard diagram of a three manifold, and the relevant algebraic data is a finite dimensional Hopf algebra. A four dimensional version of a Heegaard diagram is a trisection \cite{gay2016trisecting}, which consists of three families of circles on a closed surface such that any two families of them form a Heegaard diagram for some $\#^{n} \Sphere^1 \times \Sphere^2$.  The question is what algebraic data can be used in this case. A Hopf algebra roughly consists of a product operation and a coproduct operation, or less rigorously two product operations. Hence a naive guess in dimension four would be some \lq algebra' with three product/coproduct operations satisfying certain coherent conditions. 
\end{enumerate}

\vspace{1cm}
\noindent\textbf{Acknowledgment} $\quad$ The author acknowledges the support from the Simons Foundation. The author would also like to thank Zhenghan Wang and Cesar Galindo for helpful discussions.

\vspace{1cm}
\noindent\textbf{Note added 11/25/2018}: while the paper is under review, we noticed that the second open question raised in Section \ref{sec:future} regarding a redefinition of the $G$-crossed invariant in terms of Kirby diagrams is claimed to be solved in \cite{barenz2018evaluating}. In fact, Kirby diagrams together with some information on 3-handles are required for such a redefinition.

\end{document}